%% file: main.tex
\documentclass[a4paper,times,3p]{elsarticle}
\usepackage[utf8]{inputenc}
\usepackage[T1]{fontenc}

\usepackage{microtype}

\usepackage{booktabs}

\usepackage{amsmath,amsfonts,amsthm,enumerate,graphicx,epstopdf,url,etoolbox,algorithm,algorithmic,listings}
\usepackage{subcaption}
\usepackage{caption}
\usepackage{diagbox}
\usepackage{graphbox}
\usepackage{enumitem}
\usepackage{dirtytalk}
\usepackage{tikz}
\usepackage[export]{adjustbox}
\usetikzlibrary{shapes, arrows}
\usepackage{tcolorbox}
\usepackage{mathtools}
\usepackage{wrapfig}
\usepackage{picins}

\definecolor{cmykcyan}{cmyk}{1,0,0,0}
\definecolor{cmykred}{cmyk}{0,1,1,0}
\definecolor{cmykblack}{cmyk}{0,0,0,1}

\newtheorem{theorem}{Theorem}

\newcommand{\br}{{\mathbf{r}}}
\newcommand{\bm}{{\mathbf{m}}}
\newcommand{\bx}{{\mathbf{x}}}
\newcommand{\bxi}{{\boldsymbol{\xi}}}
\newcommand{\bs}{{\mathbf{s}}}
\newcommand{\bmu}{{\boldsymbol{\mu}}}
\newcommand{\bnu}{{\boldsymbol{\nu}}}
\newcommand{\blambda}{{\boldsymbol{\lambda}}}
\newcommand{\hOm}{{\hat{\Omega}}}
\newcommand{\Om}{{\Omega}}
\newcommand{\hOmr}{{\hat{\Omega}^{\mathbf{r}}}}

\usepackage[font=footnotesize,labelfont=bf]{caption}
\usepackage{multirow}

\usepackage{pxfonts}

\usepackage{import}

\graphicspath{./graphics/}

\usepackage[textsize=scriptsize]{todonotes}
\usepackage{setspace}

\makeatletter
\def\ps@pprintTitle{%
 \let\@oddhead\@empty
 \let\@evenhead\@empty
 \def\@oddfoot{}%
 \let\@evenfoot\@oddfoot}
\makeatother

\journal{Computer Methods in Applied Mechanics and Engineering}

\begin{document}

\begin{frontmatter}

\title{PDE-Based Parameterisation Techniques for Planar Multipatch Domains}

\author[add1]{Jochen Hinz\corref{cor1}}
\ead{jochen.hinz@epfl.ch}
\author[add1]{Annalisa Buffa}
\ead{annalisa.buffa@epfl.ch}

\cortext[cor1]{Corresponding author}

\address[add1]{Institute of Mathematics, Ecole Polytechnique F\'ed\'erale de Lausanne, 1015 Lausanne, Switzerland.}

\begin{abstract}
This paper presents a PDE-based parameterisation framework for addressing the planar surface-to-volume (StV) problem of finding a valid description of the domain's interior given no more than a spline-based description of its boundary contours. The framework is geared towards isogeometric analysis (IGA) applications wherein the physical domain is comprised of more than four sides, hence requiring more than one patch. We adopt the concept of harmonic maps and propose several PDE-based problem formulations capable of finding a valid map between a convex parametric multipatch domain and the piecewise-smooth physical domain with an equal number of sides. In line with the \textit{isoparametric paradigm} of IGA, we treat the StV problem using techniques that are characteristic for the analysis step. As such, this study proposes several IGA-based numerical algorithms for the problem's governing equations that can be effortlessly integrated into a well-developed IGA software suite. \\
We augment the framework with mechanisms that enable controlling the parametric properties of the outcome. Parametric control is accomplished by, among other techniques, the introduction of a curvilinear coordinate system in the convex parametric domain that, depending on the application, builds desired features into the computed harmonic map, such as homogeneous cell sizes or boundary layers.
\end{abstract}

\end{frontmatter}

\input{sec-intro}
\input{sec-theory}
\input{sec-numerical_schemes}
\input{sec-control_mechanisms}
\input{sec-conclusion}

\bibliographystyle{elsarticle-num}
\bibliography{main}


\end{document}

%% file: sec-intro.tex
\section{Introduction}
\label{sect:Introduction}

Isogeometric analysis (IGA) \cite{hughes2005isogeometric, cottrell2009isogeometric} is a variant of the finite element method (FEM) that was conceived in an effort to bridge the gap between the geometrical and the numerical aspects of the computational science and engineering (CSE) workflow. In computer-aided design (CAD), the physical domain $\Om$ is represented by its bounding surface $\partial \Om$ using the field's de facto standard of NURBS / spline-based parametric descriptions. The analysis step, on the other hand, relies on a geometric format based on simplices and / or relatively basic polytopes (quadrilaterals, hexahedra, $\ldots$) which form the building blocks for finding a description of the domain $\Om_h$, where $\partial \Om_h$ is a (typically piecewise linear) collocation of $\partial \Om$. Most CSE workflows operate in the order $\partial \Om \rightarrow \partial \Om_h \rightarrow \Om_h$, wherein the surface to volume (StV) problem $\partial \Om_h \rightarrow \Om_h$ is referred to as the \textit{meshing} step. The conversion from a spline-based description of $\partial \Om$ to a simplistic representation of $\Om_h$ is regarded as a major robustness bottleneck \cite{cottrell2009isogeometric}. Furthermore, in many applications it is desirable to translate analysis results provided by, for instance, FEM back to appropriate changes in $\partial \Om$, which may be nontrivial, due to the differing geometrical formats. \\

\noindent To address these concerns, IGA employs the NURBS / spline-based modelling tools that are characteristic for CAD as a basis for both the geometrical modelling and the numerical analysis aspects of the CSE workflow. In IGA, the parametric description of $\partial \Om$ is immediately forwarded to a routine that solves the StV problem $\partial \Om \rightarrow \Om$, which becomes the IGA analogue of the classical meshing step. The operator that maps the parametric domain $\hat{\Om}$ onto $\Om$ is then utilised to perform a pullback of the governing equations into $\hat{\Om}$, where the same set of splines is employed as a basis for standard FEM techniques. Besides its promise of reducing the conversion overhead, numerical simulation based on IGA is showing promising results as spline-based FEM discretisations have been demonstrated to perform better than their classical (Lagrangian) counterparts on a range of benchmark problems \cite{beirao2011some}. \\

\noindent While IGA has matured into a prolific numerical method with encouraging applications within fields such as computational electromagnetism \cite{buffa2010isogeometric} and fluid dynamics \cite{bazilevs2008nurbs}, compared to their classical meshing counterparts, spline-based methods for addressing the StV step $\partial \Om \rightarrow \Om$ are still under-represented. On the one hand, this is partially explained by the relative novelty of IGA compared to, for instance, classical FEM techniques. On the other hand, while retaining the CAD-based representation in the StV step has culminated in an entirely novel class of approaches that exploit the higher-order continuity of spline basis functions, such approaches also come with novel challenges. For instance, verifying whether a spline-based map $\bx: \hat{\Om} \rightarrow \Om$ is indeed nondegenerate is a more complicated endeavour than in the piecewise linear case \cite{xu2010optimal, gravesen2012planar}. Furthermore, spline-based representations largely rule out generalisations of classical meshing algorithms that directly operate on the mesh's vertices, such as the advancing front method \cite{lohner1988generation}, as spline-based parameterisations do generally not cross their control points. \\

\noindent As a result, the majority of existing techniques are based on blending the (typically four) segments of $\partial \Om$ into the interior \cite{farin1999discrete}, (constrained and unconstrained) parameterisation quality optimisation \cite{gravesen2012planar, xu2010optimal, xu2011parameterization, falini2015planar, xu2013constructing} and PDE-based approaches \cite{hinz2018elliptic, hinz2021iga, hinz2020goal}. While methods from all three categories, depending on the type of the geometry, show promising results, the majority have only been studied in the singlepatch setting, i.e., when the parametric domain $\hat{\Om}$ is given by the unit quadrilateral. For complex domains $\Om$, a single quadrilateral may be too restrictive, which is why existing methods may have to be combined with segmentation algorithms that divide $\Om$ into smaller pieces, which are then parameterised from the unit quadrilateral one-by-one.\\
Another challenge is associated with computational differentiability: in order to form a closed design loop, the entire CSE pipeline, including the StV step, may have to be differentiated with respect to a set of design parameters. In the presence of segmentation, this may be challenging or impossible since segmentation may not be continuous in the provided boundary data. Furthermore, segmentation may take place with little regard to parameterisation quality metrics, which include the patch interfaces in the mutlipatch setting. \\

\noindent To address the limitations of the singlepatch setting, this paper introduces a PDE-based parameterisation framework that is compatible with multipatch domains $\Om \subset \mathbb{R}^2$. The idea is to introduce a multipatch covering of an appropriately-chosen convex, polygonal parametric domain $\hat{\Om} \subset \mathbb{R}^2$ and to construct a nondegenerate mapping operator $\bx: \hat{\Om} \rightarrow \Om \subset \mathbb{R}^2$ by approximately solving a PDE problem in $\hat{\Om}$ over a spline basis defined on the multipatch topology. The underlying PDE problem approximates a map whose inverse is comprised of a pair harmonic functions in $\Om$, wherein the boundary correspondence $\bx^{-1} \vert_{\partial \Om} = \partial \hat{\Om}$ becomes the Dirichlet boundary condition. We propose two different PDE-based formulations along with various IGA-based discretisations which are then studied in detail. \\

\noindent A major appeal of this framework is the fact that the patch interfaces establish themselves as part of the PDE solution and need not be strongly imposed using, for instance, segmentation. The parameterisation, including the interfaces, is continuous in the boundary data and straightforwardly differentiable. \\

\noindent For control over the parametric properties of the computed parameterisation, we augment the framework with a mechanism that changes the properties of $\bx: \hat{\Om} \rightarrow \Om$ by mapping inversely harmonically into a parametric domain with a curvilinear, instead of a Cartesian coordinate system. This coordinate transformation is accomplished by the introduction of a so-called controlmap $\mathbf{s}: \hat{\Om} \rightarrow \hat{\Om}$. We propose several techniques for constructing controlmaps for various desired parameterisation features, such as boundary layers and boundary orthogonality. As the controlmap is defined globally (i.e., over the entire parametric domain $\hat{\Om}$), control over the parametric properties includes the image of the patch interfaces under the mapping. \\

\noindent The choice to seek $\bx: \hat{\Om} \rightarrow \Om$ as the solution of a PDE problem is in line with the isoparametric paradigm of IGA: we handle both the geometrical as well as the analysis steps using IGA techniques. As a result, the proposed algorithms are straightforwardly integrated into a well-developed IGA software suite, reducing the code bloat resulting from relying on external tools in the StV step.

\subsection{Notation}
This paper denotes vectors in boldface. The $i-th$ entry of a vector is denoted by $x_i$. Similarly, the $ij$-th entry of a matrix is denoted by $A_{ij}$. Let $\mathbf{y}: \Om \rightarrow \mathbb{R}^m$ and $\bx: \Om \rightarrow \mathbb{R}^n$. Vectorial derivatives are taken along the second axis and we interchangeably employ the denotation
\begin{align*}
    \partial_{\bx} \mathbf{y} \equiv \frac{\partial \mathbf{y}}{\partial \bx}, \quad \text{with} \quad \left[\frac{\partial \mathbf{y}}{\partial \bx} \right]_{ij} = \frac{\partial y_i}{\partial x_j},
\end{align*}
where the vector derivative $\partial_{\bx} \mathbf{y}$ maps into $\mathbb{R}^{m \times n}$. The associated Nabla operator satisfies $\nabla_{\bx} \mathbf{y} = (\partial_{\bx} \mathbf{y})^T$. \\
Furthermore, we frequently work with vector spaces $\mathcal{V}$. By default, we employ the abuse of notation
\begin{align}
    \mathcal{V}^n = \underbrace{ \mathcal{V} \times \cdots \times \mathcal{V} }_{n \text{ terms}}
\end{align}
and similarly for tensorial spaces, i.e., $\mathcal{V}^{n \times n}$. Analogously, vectorial Sobolev spaces are denoted by $H^s(\Omega, \mathbb{R}^n)$, where $\Omega$ is the associated domain. For finite-dimensional spaces $\mathcal{V}_h$, $\left \{ \mathcal{V}_h \right \}$ refers to its canonical (spline) basis which we assume to be clear from context. \\
Let $\mathcal{V}$ be defined over the domain $\Om$. We define $\mathring{\mathcal{V}} \equiv \mathcal{V} \cap H^1_0(\Om)$ as the subspace of functions from $\mathcal{V}$ that have zero trace in $\partial \Om$. \\
By $\operatorname{Int}\left(\overline{\Om} \right)$, we denote the interior of a closed domain $\overline{\Om}$, while $\overline{\Om}$ denotes the closure of an open domain $\Om$.

\subsection{Problem Statement}
\label{subsect:problem_statement}
Let $\Om \subset \mathbb{R}^2$ be an open, simply connected Lipschitz domain whose boundary $\partial \Om$ is parameterised by 
an even number $K = 2n, \enskip n \in \mathbb{N}^{\geq 2}$ of open (spline) curves $C_k \subset \mathbb{R}^2$, oriented in counterclockwise direction. We have
$$\partial \Om = \bigcup \limits_{\mathclap{k \in \{1, \ldots, K\}}} \overline{C}_k, \quad \text{where} \quad i \neq j \implies C_i \cap C_j = \emptyset.$$
We assume that the $C_k$ are parameterised in the positive direction from the open unit interval by the spline maps $\mathbf{f}^k: (0, 1) \rightarrow \mathbb{R}^2$ with $\mathbf{f}^k \in C^1((0, 1), \, \mathbb{R}^2)$ and nonvanishing tangent. \\
Furthermore, let $\hat{\Om} \subset \mathbb{R}^2$ be a convex, polygonal parametric domain with $K$ sides $L_k \subset \mathbb{R}^2$ oriented in counterclockwise direction. We have
$$\partial \hat{\Om} = \bigcup \limits_{\mathclap{k \in \{1, \ldots, K\}}} \overline{L}_k, \quad \text{where} \quad i \neq j \implies \enskip L_i \cap L_j = \emptyset.$$
Each $L_k$ is parameterised in the positive direction on $\partial \hat{\Om}$ by an affine map $\mathbf{l}_k: (0, 1) \rightarrow \mathbb{R}^2$ of the form:
$$\mathbf{l}_k(s) = \bxi^k + \mathbf{t}^k s, \quad \text{with} \quad \{\bxi^k, \mathbf{t}^k\} \subset \mathbb{R}^2.$$
Assigning the $C_k$ to the $L_k$ in ascending order induces the boundary correspondence $\mathbf{F}: \partial \hat{\Om} \rightarrow \partial \Om$ that satisfies
\begin{align}
    \mathbf{F} \vert_{\overline{L}_k} = \overline{C}_k, \quad \text{or equivalently} \quad \mathbf{F} \circ \mathbf{l}_k = \mathbf{f}_k,
\end{align}
and we assume that $\mathbf{F}: \partial \hat{\Om} \rightarrow \partial \Om$ parameterises a Jordan curve in $\mathbb{R}^2$. \\
We assume that $\hat{\Om}$ is covered by a quadrangulation $\mathcal{Q}$ of a total of $N_p$ patches $\hat{\Om}_i$, i.e.,
\begin{align}
    \mathcal{Q} = \{ \hat{\Om}_1, \ldots, \hat{\Om}_{N_p} \}, \quad \text{with} \quad \hat{\Om} = \operatorname{Int} \left( \bigcup \limits_{\hat{\Om}_i \in \mathcal{Q}} \overline{\hat{\Om}}_i \right) \quad \text{and} \quad i \neq j \implies \hat{\Om}_i \cap \hat{\Om}_j = \emptyset.
\end{align}
Each $\hat{\Om}_i$ is the image of the reference patch $\Om^{\square} = (0, 1)^2$ under the diffeomorphic bilinear map $\bm^i: \Om^{\square} \rightarrow \hat{\Om}_i$. The facets of the quadrangulation are denoted by $\Gamma$, while boundary facets are denoted by $\Gamma^B := \{L_1, \ldots, L_K\}$ and interior facets by $\Gamma^I := \Gamma \setminus \Gamma^B$. \\
For boundary patches $\hat{\Om}_i$, the associated map $\bm^i$ restricted to the side of $\partial \Om^{\square}$ that maps onto $L_k$, is given either by $\mathbf{l}_k$ or $\mathbf{l}_k(1 - s)$, depending on the orientation along $\partial \hat{\Om}$. We denote the set of boundary patches by $\mathcal{Q}^B$.\\
The facets between pairs of neighbouring patches are denoted by $\gamma_{ij}$ and the collection of interior facets is given by
$$\Gamma^I = \bigcup \limits_{\mathclap{(i, j) \in F^I}} \gamma_{ij}, \quad \text{with} \quad F^I := \left \{(i, j) \, \, \vert \, \, \operatorname{Int} \left(\overline{\hat{\Om}_i} \cap \overline{\hat{\Om}}_j \right) \text{ is an open line segment in } \hat{\Om} \right \}.$$
\noindent Given no more than a boundary correspondence $\mathbf{F}: \partial \hOm \rightarrow \partial \Om$ that satisfies aforementioned assumptions, this paper deals with the spline-based StV problem $\partial \Om \rightarrow \Om$. More precisely, let $\mathcal{V}_h \subset H^1(\hat{\Om})$ be a finite-dimensional vector space and let 
$$\mathcal{U}_h^{\mathbf{F}} = \left \{\mathbf{v} \in \mathcal{V}_h^2 \enskip \vert \enskip \mathbf{v} = \mathbf{F} \enskip \text{on} \enskip \partial \hat{\Om} \right \}, \quad \text{with } \mathcal{V}_h \text{ such that} \quad \mathcal{U}^{\mathbf{F}}_h \neq \emptyset.$$
The purpose of this paper is providing a framework for finding a nondegenerate mapping operator $\bx_h: \hat{\Om} \rightarrow \Om$ with $\bx_h \in \mathcal{U}^{\mathbf{F}}_h$. Denoting the Cartesian coordinate functions in $\hat{\Om}$ by $\bxi = (\xi_1, \xi_2)^T$, we call a map $\bx: \hOm \rightarrow \Om$ \textit{nondegenerate} (NDG) if
\begin{align}
    0 \leq \inf \limits_{\bxi \in \hOm} \, \, \det J(\bx) \leq \sup \limits_{\bxi \in \hOm} \, \, \det J(\bx) \leq \infty, \quad \text{where} \quad J(\mathbf{x}) := \frac{\partial \mathbf{x}}{\partial \bxi}
\end{align}
denotes the Jacobian matrix of $\mathbf{x}: \hat{\Om} \rightarrow \Om$ in $\hat{\Om}$. Similarly, we call a map \textit{uniformly nondegenerate} (UNDG) if 
\begin{align}
    0 < c < \inf \limits_{\bxi \in \hOm} \, \, \det J(\bx) \leq \sup \limits_{\bxi \in \hOm} \, \, \det J(\bx) < C < \infty.
\end{align}
Clearly, uniform nondegeneracy of $\bx_h \in \mathcal{V}_h^2$ is favoured over nondegeneracy by most applications but imposes stronger requirements on $\mathbf{F}: \partial \hOm \rightarrow \partial \Om$ that are discussed in Section \ref{sect:boundary_correspondence_requirements}. \\
Letting $\mathcal{V}_h$ be spanned by $\{\phi_1, \ldots, \phi_N\}$, the mapping operator takes the form:
\begin{align}
\label{eq:mapping_general_form}
    \bx_h(\xi_1, \xi_2) = \sum_{i \in \mathcal{I}_I} \mathbf{c}^i \phi_i(\xi_1, \xi_2) + \sum_{j \in \mathcal{I}_B} \mathbf{c}^j \phi_j(\xi_1, \xi_2),
\end{align}
where $\mathbf{c}^k \in \mathbb{R}^2, \enskip \forall k \in \mathcal{I}_I \cup \mathcal{I}_B$, while $\mathcal{I}_I$ and $\mathcal{I}_B$ refer to the index-sets of vanishing and nonvanishing functions on $\partial \hat{\Om}$, respectively. With~\eqref{eq:mapping_general_form} in mind, the purpose of this paper is properly selecting the $\mathbf{c}^i$ while the $\mathbf{c}^j$ follow from the boundary correspondence and are therefore held fixed. \\
Besides nondegeneracy, this paper aims for mechanisms that allow for control over the parametric properties of $\bx_h: \hat{\Om} \rightarrow \Om$ while the framework should be implicitly differentiable, i.e., provide maps that are a continuous function of the supplied data, namely the boundary control points $\mathbf{c}^j, \enskip j \in \mathcal{I}_B$.

\subsection{Related Work}
\label{sect:related_work}
As stated in the introduction, existing techniques for the StV problem $\partial \Om \rightarrow \Om$ are predominantly based on blending the curves $C_k$ that make up $\partial \Om$ into the interior, selecting the $\mathbf{c}^i \in \mathcal{I}_I$ via an optimisation problem (with or without added constraints) and PDE-based methods. So far, most methods have only been studied in the singlepatch setting. \\

\noindent Interpolation-based methods, such as transfinite interpolation \cite{gordon1973transfinite}, are a class of approaches that had already been conceived before the onset of IGA. Such approaches attempt to parameterise the interior of $\partial \Om$ by taking the map $\bx_h: \hat{\Om} \rightarrow \mathbb{R}^2$ as a linear combination of the $\mathbf{f}_k: (0, 1) \rightarrow C_k$ times a set of (typically polynomial) blending functions defined in $\hat{\Om}$. In IGA, the most widely-used method is the bilinearly blended Coons' patch \cite{farin1999discrete}, a computationally inexpensive and often sufficiently powerful approach for singlepatch geometries. More advanced variants, such as Lagrange and Hermite interpolation, furthermore allow for control over the map's derivatives on $\partial \hat{\Om}$ using blending functions of polynomial degree $p \geq 2$. For an overview of interpolation techniques over the unit quadrilateral, see \cite[Chapter~5]{liseikin1999grid}. Generalisations to $n$-sided convex, polygonal domains have been made in \cite{varady2011transfinite} and \cite{salvi2014ribbon}, wherein the construction of appropriate blending functions becomes the main objective. As blending is based on polynomial constructions, interpolated surfaces can typically be equivalently expressed in local constructions based on splines. While computationally inexpensive and often highly effective in practical applications, interpolation-based methods provide no guarantee of nondegeneracy and the resulting maps are therefore often folded. \\

\noindent The second class of approaches minimises one or a positive sum of several quality cost functions over the $\mathbf{c}^i, \enskip i \in \mathcal{I}_I$. Compared to the classical literature, optimisation-based techniques have received more interest within the IGA-realm. Quality criteria are largely based on heuristics and optimisation typically seeks for maps with orthogonal isolines, homogeneous cell sizes or reduced cell skewness \cite{falini2015planar, buchegger2017planar}. Convex optimisation formulations are based on the length, Liao and uniformity functionals \cite{falini2015planar, hinz2020goal}, while nonconvex formulations are often based on a combination of the area, orthogonality and skewness functionals \cite{buchegger2017planar}. While computationally more demanding, nonconvex optimisation has a lower tendency to yield degenerate maps and allows for a wider range of quality criteria \cite[Chapter~6]{steinberg1993fundamentals}. Further examples include the \textit{Teichm\"uller} map \cite{nian2016planar} and the \textit{variational harmonic method} \cite{xu2013constructing}, which, given a sufficiently regular boundary correspondence, both approximate a bijection and are thus inherently less prone to yielding a degenerate map. In \cite{falini2015planar}, the most-commonly employed cost functions are studied in a THB-spline setting. \\

\noindent Penalisation methods enforce nondegeneracy by adding the Jacobian determinant $\det J(\bx_h)$ to the cost function's denominator thus creating a barrier that urges the optimiser to seek for local minima from within set of bijective maps. Examples include the \textit{Winslow} and \textit{modified Liao} functionals \cite[Chapter~8]{steinberg1993fundamentals} whose minimisation has to be initialised with a nondegenerate map to avoid division by zero. To relax this requirement, \cite{wang2021smooth} introduces a regularisation that enables degenerate initial iterates to converge to a valid map in the proximity of the original formulation's global minimiser. However, the radius of convergence remains small and a suitable initial iterate requires solving another optimisation problem first. \\

\noindent Another way of enforcing nondegeneracy is adding constraints that constitute a sufficient condition for bijectivity. A linear constraint is proposed in \cite{xu2011parameterization}. If convex cost functions are utilised, the problem remains convex. A nonconvex constraint is proposed in \cite{gravesen2012planar} wherein the map's (scalar) Jacobian determinant is expressed in a spline space $\mathcal{V}_h^J$ that contains it. If the determinant's weights with respect to $\{ \mathcal{V}_h^J \}$ are all positive, the map is valid and the iterate is deemed feasible. Expanding $\det J(\bx_h)$ over $\{ \mathcal{V}_h^J \}$ is furthermore a widely-used technique to test for nondegeneracy. As both constraints constitute sufficient but not necessary conditions for bijectivity, they may be too restrictive in practice. Furthermore, finding a feasible initial iterate may be nontrivial or impossible because for complex geometries, the feasible search space may be empty. \\

\noindent Optimisation-based approaches are readily generalised to the multipatch setting by minimising the same cost functions over the polygonal domain $\hat{\Om}$ rather than the unit square. Hereby, the patch interface control points become degrees of freedom in the formulation. Multipatch optimisation is employed in \cite{buchegger2017planar} where a suitable topology is chosen through a construction based on patch adjacency graphs. In \cite{bastl2021planar}, a time-dependent formulation is proposed which evolves the initial $\bx_h(t=0)$, that maps strictly into the interior of $\partial \Om$, to a map with the prescribed boundary correspondence. At each time-iteration, multipatch optimisation is utilised to warrant the parametric quality of the intermediate map. \\
To the best of our knowledge,  penalised or constrained optimisation problems have only been studied in the singlepatch setting. \\

\noindent The third class of approaches seeks the $\mathbf{c}^i, \enskip i \in \mathcal{I}_I$ by (approximately) solving a PDE problem. In the majority of cases, the PDE stems from the requirement that the mapping inverse $\bx_h^{-1}: \Om \rightarrow \hat{\Om}$ be a pair of harmonic functions in $\Om$. A justification for this is provided by the \textit{Rad\'o-Kneser-Choquet theorem} which states that a harmonic map $\bx_h^{-1}: \Om \rightarrow \hat{\Om}$ is diffeomorphic in $\Om$, provided $\hat{\Om}$ is convex. Thus, requiring the mapping inverse to be harmonic allows for treating a wider range of geometries as the parametric domain can be chosen freely. In \cite{falini2015planar} harmonic maps are approximated by employing a boundary element method (BEM) \cite{simpson2012two}. The method creates a large number of pairs $(\bxi_i, \bx_i)$, with $\bxi_i \in (0, 1)^2$ and $\bx_i \in \Om$ which are then utilised to fit a THB-spline map in the least-squares sense with added regularisation terms. The same BEM is adopted in \cite{falini2019thb} where multiply-connected domains are mapped inversely harmonically into punctured auxiliary domains, using a template segmentation approach to select an appropriate multipatch layout for fitting a map to the point pairs. \\
Another class of methods seeks (inversely) harmonic maps by approximately solving the equations of \textit{Elliptic Grid Generation} \cite[Chapter~5]{steinberg1993fundamentals} (EGG) in $\hat{\Om}$ using variational techniques. EGG stems from a pullback of the inverse harmonicity requirement into $\hat{\Om}$ where a (nonlinear) equation for the forward map is derived. In \cite{hinz2018elliptic} this formulation is employed to approximate harmonic maps on quadrilateral parametric domains. As the EGG equations are of second order, a spline space $\mathcal{V}_h \subset H^2(\hat{\Om})$ is employed. In \cite{hinz2020goal}, the equations are tackled with THB-splines combined with a posteriori refinement techniques to repair degeneracies stemming from insufficient numerical accuracy. The same publication proposes control mechanisms capable of tuning the map's parametric properties by introducing a suitable coordinate transformation in $\hat{\Om}$. An attempt to generalise IGA-EGG to the multipatch setting is made in \cite{hinz2021iga} where a mixed form that reduces the regularity requirement from $H^2(\hat{\Om})$ to $H^1(\hat{\Om})$ is proposed by introducing auxiliary variables for the map's Jacobian. While successful in practice, the formulation significantly increases the computational costs since the basis is likely subject to an inf-sup requirement \cite{bathe2001inf}, requiring the auxiliary space to be $p$- or $h$-refined with respect to the primal basis. \\
Finally, \cite{shamanskiy2020isogeometric} proposes a PDE-based approach that employs the equations of nonlinear elasticity with applications to both singlepatch and multipatch domains.

%% file: sec-theory.tex
\section{Theory}
This paper proposes a framework for computing parameterisations $\bx_h: \hOm \rightarrow \Om$ based on harmonic maps. In the following, we present an in-depth discourse on harmonic maps as well as finite element techniques for elliptic equations in nonvariational form, which shall be adopted to formulate discretisations in Section \ref{subsect:NDF_discretisations}.

\subsection{Harmonic Maps}
\label{subsect:harmonic_maps}
The motivation to seek the map $\bx_h: \hOm \rightarrow \Om$ as the inverse of a map that is harmonic in $\Om$ stems from the following famous result:
\begin{theorem}[Rad\'o-Kneser-Choquet]
\label{thrm:RKC}
    The harmonic extension of a homeomorphism from the boundary of a Jordan domain $\Om \subset \mathbb{R}^2$ onto the boundary of a convex domain $\hOm \subset \mathbb{R}^2$ is a diffeomorphism in $\Om$.
\end{theorem}
For proofs, we refer to \cite{kneser1926losung, choquet1945type, sauvigny1991embeddedness, castillo1991mathematical}. It should be noted that the convexity of $\hOm$ is a sufficient, but not a necessary condition. Furthermore, the same result is no longer true in $\mathbb{R}^3$ \cite{laugesen1996injectivity}. \\
Theorem \ref{thrm:RKC} has inspired many numerical approaches for finding a nondegenerate $\bx_h: \hOm \rightarrow \Om$ that approximates a map whose inverse is harmonic. Besides the nondegeneracy guarantee, this is furthermore explained by the regularity of harmonic maps which generally serve the map's quality from a numerical standpoint. \\
Numerical approaches go back to the pioneering works of Winslow \cite{winslow1966numerical}. Letting $\bx = (x_1, x_2)^T$, and defining the metric tensor 
$$G_{ij}(\bx) = g_{ij} \quad \text{with} \quad g_{ij} = \partial_{\xi_i} \bx \cdot \partial_{\xi_j} \bx,$$
Winslow's original approach seeks the map $\bx: \hOm \rightarrow \Om$ as result of the following minimisation problem
\begin{align}
    \frac{1}{2} \int \limits_{\Om}\operatorname{tr}\left( G^{-1} \right) \mathrm{d} \bx \rightarrow \min \limits_{\bx}, \quad \text{s.t.} \quad \bxi(\bx) = \mathbf{F}^{-1} \text{ on } \partial \Om.
\end{align}
Letting $\mathcal{U}^{\mathbf{F}} = \{ \mathbf{v} \in H^1(\hOm, \mathbb{R}^2) \,\, \vert \,\, \mathbf{v} = \mathbf{F} \enskip \text{on} \enskip \partial \hOm \}$, a pullback leads to 
\begin{align}
\label{eq:Winslow_pullback}
    \frac{1}{2} \int \limits_{\hOm} \frac{\operatorname{tr}(G)}{\det J} \mathrm{d} \bxi \rightarrow \min \limits_{\bx \in \mathcal{U}^\mathbf{F}},
\end{align}
while a discretisation replaces $\mathcal{U}^{\mathbf{F}} \rightarrow \mathcal{U}^{\mathbf{F}}_h$ in (\ref{eq:Winslow_pullback}). The minimisation of~\eqref{eq:Winslow_pullback} is highly impractical since the domain of the integrand is the set of all $\bx \in \mathcal{U}^{\mathbf{F}}$ that satisfy $\det J(\bx) > 0$ (almost everywhere). As such, minimisation has to be initialised with a nondegenerate initial map which is generally hard to find. \\
An alternative formulation is based on the harmonicity requirement's \textit{classical form}:
\begin{align}
    \Delta \bx^{-1} = 0 \quad \text{ in } \Om, \quad \text{s.t.} \quad \bx^{-1} = \mathbf{F}^{-1} \text{ on } \partial \Om,
\end{align}
where the Laplace operator is to be understood component-wise. A pullback leads to
\begin{align}
\label{eq:elliptic_classical_pullback}
    \Delta_{\bx} \bxi = 0 \quad \text{in } \hOm, \quad \text{s.t.} \quad \bx = \mathbf{F} \text{ on } \partial \hOm,
\end{align} 
where $\Delta_\bx$ denotes the Laplace-Beltrami operator. \\
The pullbacks from both~\eqref{eq:Winslow_pullback} and~\eqref{eq:elliptic_classical_pullback} inherently assume that $\bx^{-1}: \Om \rightarrow \hOm$ is invertible, thus potentially rendering the problems ill-posed in case $\hOm$ is not convex. However, assuming convexity of $\hOm$, we may multiply the two-component PDE from~\eqref{eq:elliptic_classical_pullback} by $T: \hOm \rightarrow \mathbb{R}^{2 \times 2}$, with $T = \left( \det J \right)^2 J(\bx)$ since $T$ does not vanish in the interior. The result is a two-component PDE for $\bx: \hOm \rightarrow \Om$ which can be classified as a quasilinear second-order elliptic PDE in nondivergence form \cite{hinz2018elliptic, hinz2020goal, thompson1998handbook}:
\begin{align}
\label{eq:EGG_classical}
    i \in \{1, 2\}: \quad A(\partial_{\bxi} \bx) \colon H(x_i) = 0, \quad \text{s.t.} \quad \bx = \mathbf{F} \text{ on } \partial \hOm,
\end{align}
where
$$H(y)_{ij} = \frac{\partial^2 y}{\partial \xi_i \partial \xi_j} \quad \text{denotes the Hessian in } \hOm, \quad \text{while} \quad A(\partial_{\bxi} \bx) := \begin{pmatrix} g_{22} & -g_{12} \\ -g_{12} & g_{11} \end{pmatrix}$$
and $A \colon B$ denotes the Frobenius inner product between two matrices. \\
The multiplication by $T: \hOm \rightarrow \mathbb{R}^{2 \times 2}$ removes $\det J$ from the original formulation's denominator, allowing schemes based on~\eqref{eq:EGG_classical} to be initialised with degenerate initial maps. On the other hand, minimisation based on~\eqref{eq:Winslow_pullback} most likely yields a nondegenerate map, while this may not hold for approaches based on~\eqref{eq:EGG_classical}, due the scheme's truncation error. \\
A further possibility is basing the scheme on the harmonicity requirement's \textit{weak form}. More precisely, with $\mathcal{V} = H^1(\Om)$:
\begin{align}
\label{eq:inverse_elliptic_weak}
    \text{find } \bx^{-1} \in \mathcal{V}^2, \quad \text{s.t.} \quad \int \limits_{\Om} \nabla \boldsymbol{\phi} \colon \nabla \bx^{-1} \mathrm{d} \bx = 0, \quad \forall \boldsymbol{\phi} \in \mathring{\mathcal{V}}^2 \quad \text{and} \quad \bx^{-1} = \mathbf{F}^{-1} \text{ on } \partial \Om,
\end{align}
which translates to an equation for $\bx: \hOm \rightarrow \Om$ via a pullback. \\
This paper presents algorithms for approximating $\bx_h \approx \bx$ based on formulations~\eqref{eq:EGG_classical} and~\eqref{eq:inverse_elliptic_weak}. As the former is in nondivergence form, in the following we give a brief summary on the finite element treatment of nondivergence form equations.

\subsection{Nondivergence form equations}
\label{subsect:NDF}
The finite element treatment of nondivergence form (NDF) equations is a relatively recent development with first contributions due to Lakkis and Pryer \cite{lakkis2011finite}. NDF-equations are of the form
\begin{align}
    B \colon H(u) + \text{ lower order terms } = f \quad \text{a.e. in } \hOm, \quad \text{s.t.} \quad u = g \text{ on } \partial \hOm.
\end{align}
Here, $f \in L^2(\hOm)$, while $B: \hOm \rightarrow L^\infty(\hOm, \mathbb{R}^{2 \times 2})$ is uniformly elliptic, i.e., there are constants $0 < c_1 \leq c_2 < \infty$ such that
\begin{align}
\label{eq:uniformly_elliptic}
    c_1 \leq \inf \limits_{\bxi \in \mathbb{R}^2, \| \bxi \| = 1} \bxi^T B \bxi \leq \sup \limits_{\bxi \in \mathbb{R}^2, \| \bxi \| = 1} \bxi^T B \bxi \leq c_2, \quad \text{ a.e. in } \hOm.
\end{align}
The set of all symmetric and uniformly elliptic $B: \hOm \rightarrow L^\infty(\hOm, \mathbb{R}^{2 \times 2})$ is referred to as $\text{SPD}^{2 \times 2}(\hOm)$. 
In the following, we take $g = 0$ and disregard lower order terms for convenience. \\
For $g = 0$ and $\hOm \subset \mathbb{R}^2$ convex, it can be shown that $u \in H^2(\hOm) \cap H^1_0(\hOm)$ as long as $B$ satisfies the so-called Cord\'es condition \cite{maugeri2000elliptic}. In $\mathbb{R}^2$ the Cord\'es condition is implied by~\eqref{eq:uniformly_elliptic}. Defining 
\begin{align}
\label{eq:NDF_scaling}
    \gamma(B) := \operatorname{tr}(B) / \| B \|_F^2,
\end{align} 
where $\| \, \cdot \, \|_F$ denotes the Frobenius norm $\sqrt{A \colon A}$ of a matrix, finite element discretisations are based on the following Petrov-Galerkin formulation of the problem's \textit{strong form}:
\begin{align}
\label{eq:NDF_weak_discrete}
    \text{find } u \in H^2(\hOm) \cap H^1_0(\hOm) \quad \text{s.t.} \quad  \int \limits_{\hOm} \gamma(B) \tau(\phi) \left( B \colon H(u)  - f \right) \mathrm{d} \bxi = 0, \quad \forall \phi \in \mathcal{V},
\end{align}
for some suitably-chosen test space $\mathcal{V}$. \\
Here, $\tau: \mathcal{V} \rightarrow L^2(\hOm)$ is a suitably-chosen operator that warrants coercivity of the associated bilinear form over finite-dimensional subspaces $\mathcal{V}_h \subset \mathcal{V}$. The (optional) scaling $\gamma(\, \cdot \,)$ guarantees that $\gamma(B) \, B$ resembles the identity matrix $\mathcal{I}^{2 \times 2}$ and simplifies the analysis of numerical schemes based on~\eqref{eq:NDF_weak_discrete}. The choices of $\tau: \mathcal{V} \rightarrow L^2(\hOm)$ for $\mathcal{V} = H^2(\hOm) \cap H^1_0(\hOm)$ are $\tau_{\text{NS}}(v) = \Delta v$ and $\tau_{\text{LS}}(v) = B \colon H(v)$ \cite{blechschmidt2021error, gallistl2017variational, smears2013discontinuous}, while for $\mathcal{V} = H^1_0(\hOm)$,  $\tau_{\text{ID}}(v) = v$ \cite{feng2017finite, lakkis2011finite}. To enable discretisations over finite element spaces $\mathcal{V}_h \subset H^1(\hOm)$, mixed-FEM formulations of~\eqref{eq:NDF_weak_discrete} are introduced in \cite{lakkis2011finite, gallistl2017variational} while \cite{smears2013discontinuous, blechschmidt2021error} propose $C^0$ discontinuous Galerkin schemes that introduce interior penalty terms over the facets of the FEM mesh, acting on the discrete solution's normal gradient. \\
As spaces $\mathcal{V}_h$ resulting from local spline-based constructions over multipatch topologies are generally only in $H^1(\hOm)$, this paper adopts the mixed formulations based on Gallistl \cite{gallistl2017variational} and Lakkis-Pryer \cite{lakkis2011finite} as well as the $C^0$-DG formulation from \cite{blechschmidt2021error} and applies them to linearisations of~\eqref{eq:EGG_classical}. In the case of $C^0$-DG, penalty terms can be restricted to the interior patch interfaces $\gamma_{ij} \in \Gamma^I$ of $\hOm$.

%% file: sec-numerical_schemes.tex
\section{Numerical Schemes}
\label{sect:numerical_schemes}
In this section we propose several numerical schemes for finding approximate solutions of the inverse harmonicity formulations based on both~\eqref{eq:EGG_classical} and~\eqref{eq:inverse_elliptic_weak}. Given $\hOm$ and a suitable multipatch covering (see Subsection \ref{subsect:problem_statement}), we denote by $\boldsymbol{\Xi}_i = (\Xi_{i, 1}, \Xi_{i, 2})$ the pair of local (open) knotvectors associated with the $i$-th patch along with the associated canonical spline space \smash{$\widehat{\mathcal{V}}_{h, i} \subset H^2(\Om^{\square})$}. We denote by \smash{$\mathcal{V}_h^{\text{disc}} \subset L^2(\hOm)$} the space that results from a push-forward of the local spaces \smash{$\widehat{\mathcal{V}}_{h, i}$}, i.e,
$$\mathcal{V}_h^{\text{disc}} = \bigcup \limits_{i \in \{1, \ldots, N_p \}} \left\{ v_h \circ (\mathbf{m}^i)^{-1} \, \, \vert \, \, v_h \in \widehat{\mathcal{V}}_{h, i} \right \}.$$
Then, we define the subspace $\mathcal{V}_h := \mathcal{V}_h^{\text{disc}} \cap C^0(\hOm)$ and assume that the local spline bases' knotvector tuples $\boldsymbol{\Xi}_i$ are selected in such a way that the canonical basis $\{ \mathcal{V}_h \}$ of $\mathcal{V}_h$ forms a partition of unity on $\hOm$ that is compatible with $\mathbf{F}: \partial \hOm \rightarrow \partial \Om$ in the sense that the set $\mathcal{U}^{\mathbf{F}}_h = \{ \mathbf{v} \in \mathcal{V}_h^2 \, \, \vert \, \, \mathbf{v} = \mathbf{F} \text{ on } \partial \hOm \} \neq \emptyset$.

\subsection{NDF discretisations}
\label{subsect:NDF_discretisations}
The discretisations of~\eqref{eq:EGG_classical} are based on a variation of the Petrov-Galerkin formulation from~\eqref{eq:NDF_weak_discrete}. Here, we restrict ourselves to the choices $\tau \in \{\tau_{\text{NS}}, \tau_{ID} \}$.
In the following, we propose iterative solution strategies targeting a variational form of~\eqref{eq:EGG_classical}. For the sake of a unified presentation, we let $\mathcal{U}^{\mathbf{F}} := \{\mathbf{v} \in H^2(\hOm, \mathbb{R}^2) \, \vert \, \mathbf{v} = \mathbf{F} \text{ on } \partial \hOm \}$ and we introduce the form $\mathcal{L}: \text{SPD}^{2 \times 2}(\hOm) \times \mathcal{U}^{\mathbf{F}} \times \mathcal{U}^{\mathbf{0}} \rightarrow \mathbb{R}$ with

\begin{align}
\label{eq:NDF_strong_operator}
    \mathcal{L}(B, \bx, \boldsymbol{\phi}) := \int \limits_{\hOm} \tau(\phi_i) B \, \colon \, H(x_i) \, \mathrm{d} \bxi,
\end{align}
where we sum over repeated indices. For the time being, we assume that the data is sufficiently regular for the problem and its linearisations to be well-posed over $\mathcal{U}^{\mathbf{F}}$ for the test space $\mathcal{U}_{\text{test}} = \mathcal{V}_{\text{test}}^2$ with $\mathcal{V}_{\text{test}} = H^1_0(\hOm)$ ($\tau = \tau_{\text{ID}}$) and $\mathcal{V}_{\text{test}} = H^2(\hOm) \cap H^1_0(\hOm)$ ($\tau = \tau_{\text{NS}}$), respectively. The linearisations are then modified for compatibility with the $C^0(\hOm)$-nature of spline spaces over multipatch topologies while discretisations follow readily from replacing vector spaces by their finite-dimensional counterparts.  \\
In what follows, we shall substitute various flavours of $A(\, \cdot \,)$ (cf. Subsection \ref{subsect:harmonic_maps}) scaled by $\gamma(\, \cdot \,)$ (cf. Subsection \ref{subsect:NDF}) into~\eqref{eq:NDF_strong_operator}. Besides being customary in NDF-discretisations, we have noticed the scaling to have a positive effect on the iterative schemes' radii of convergence and the required number of iterations.

\subsubsection{Fixed-Point Iteration}
\label{subsect:NDF_fixed}
The most elementary linearisation is based on a fixed-point iteration which freezes $A( \, \cdot \,)$ of~\eqref{eq:EGG_classical} in the previous iterate $\bx^k$ and seeks $\bx: \hOm \rightarrow \Om$ as the limit $k \rightarrow \infty$ of the recursion
\begin{align}
    i \in \{1, 2\}: \quad A(\partial_{\bxi} \bx^k) \, \colon \, H(x^{k+1}_i) = 0, \quad \text{s.t.} \quad \bx^{k+1} = \mathbf{F} \text{ on } \hOm.
\end{align}
We note that $A(\partial_{\bxi} \bx)$ may equivalently be written in the form
\begin{align}
\label{eq:A_CC}
    A(\partial_{\bxi} \bx) = C^T C, \quad \text{with} \quad C(\partial_{\bxi} \bx) = \begin{pmatrix} \tfrac{\partial x_2}{\partial \xi_2} & -\tfrac{\partial x_2}{\partial \xi_1} \\ - \tfrac{\partial x_1}{\partial \xi_2} & \tfrac{\partial x_1}{\partial \xi_1} \end{pmatrix}.
\end{align}
Since $C(\, \cdot \,)$ has the same characteristic polynomial as $J = \partial_{\bxi} \bx$, we conclude that $A(\, \cdot \,) \in \text{SPD}^{2 \times 2}(\hOm)$ whenever $\bx: \hOm \rightarrow \mathbb{R}^2$ is UNDG. As such, the uniform ellipticity requirement~\eqref{eq:uniformly_elliptic} is violated for degenerate intermediate maps. To circumvent this issue, we introduce the stabilisation $A_\mu(\, \cdot \,) := A( \, \cdot \,) + \mu \mathcal{I}^{2 \times 2}$, with $0 < \mu < 1$ and base a numerical scheme on the following linearised classical form:
\begin{align}
    i \in \{1, 2\}: \quad A_{\mu}(\partial_{\bxi} \bx^k) \, \colon \, H(x^{k+1}_i) - \mu \Delta x^k_i = 0, \quad \text{s.t.} \quad \bx^{k+1} = \mathbf{F} \text{ on } \hOm.
\end{align}
With $A_{\mu}^k := A_{\mu}(\partial_{\bxi} \bx^k)$ and $\gamma_\mu^k := \gamma(A^k_\mu)$ (cf. equation~\eqref{eq:NDF_scaling}), a variational formulation seeks $\bx: \hOm \rightarrow \Om$ as the limit $k \rightarrow \infty$ of the recursion
\begin{align}
    \text{find } \bx^{k+1} \in \mathcal{U}^{\mathbf{F}} \quad \text{s.t.} \quad \mathcal{F}_{\mu}( \bx^{k+1}, \bx^{k}, \boldsymbol{\phi}) = 0, \quad \forall \boldsymbol{\phi} \in \mathcal{U}^{\mathbf{0}},
\end{align}
where
\begin{align}
    \mathcal{F}_{\mu}(\bx^{k+1}, \bx^k, \boldsymbol{\phi}) := \mathcal{L}(\gamma_\mu^k \, A_\mu^k, \bx^{k+1},  \boldsymbol{\phi}) - \mu \mathcal{L}(\gamma_\mu^k \, \mathcal{I}^{2 \times 2}, \bx^{k}, \boldsymbol{\phi}).
\end{align}
In practice, we take $\mu = 10^{-4}$. Here, a reasonable stopping criterion terminates the recursion as soon as $\| \bx^{k+1} - \bx^{k} \| / \| \bx^k \| < \varepsilon$, where a suitable norm depends on the augmented scheme (with $C^0$-support). 

\subsubsection{Newton Approach}
\label{subsect:NDF_newton}
As in the fixed-point iteration, a linearisation based on Newton's method needs to be adjusted for the possibility of encountering iterates $\bx^{k}$ with $A(\partial_{\bxi} \bx^{k}) \notin \text{SPD}^{2 \times 2}$. As such, we again employ the eigenspectrum shift $A( \, \cdot \,) \rightarrow A_{\mu}(\, \cdot \,)$  and base a Newton scheme on the residual form $\mathcal{N}_{\mu}: \mathcal{U}^{\mathbf{F}} \times \mathcal{U}^{\mathbf{0}} \rightarrow \mathbb{R}$, with
\begin{align}
    \mathcal{N}_{\mu}(\bx, \boldsymbol{\phi}) := \mathcal{L}(\gamma_{\mu}^{\bx} A_{\mu}^{\bx}, \bx, \boldsymbol{\phi}), \quad \text{where} \quad A_{\mu}^{\bx} := A_{\mu}(\partial_{\bxi} \bx) \quad \text{ and } \quad \gamma_{\mu}^{\bx} := \gamma(A_\mu^{\bx}).
\end{align}
Given some intermediate iterate $\bx^{k} \in \mathcal{U}^{\mathbf{F}}$, the Newton scheme computes the increment $\partial \bx^k \in \mathcal{U}^{\mathbf{0}}$ from:
\begin{align}
\label{eq:EGG_newton_continuous}
    \text{find } \partial \bx^k \in \mathcal{U}^{\mathbf{0}}, \quad \text{s.t.} \quad \mathcal{N}_{\mu}^{\prime}(\bx^k, \boldsymbol{\phi}, \partial \bx^k) = - \mathcal{N}_{\mu}(\bx^k, \boldsymbol{\phi}), \quad \forall \boldsymbol{\phi} \in \mathcal{U}^{\mathbf{0}},
\end{align}
wherein $\mathcal{N}_{\mu}^{\prime}(\, \cdot \,, \, \cdot \,, \mathbf{v})$ denotes the Gateaux derivative of $\mathcal{N}_{\mu}(\, \cdot \,, \, \cdot \,)$ with respect to its first argument in the direction of $\mathbf{v} \in \mathcal{U}^{\mathbf{0}}$. The new iterate becomes $\bx^{k+1} = \bx^k + \kappa \partial \bx^k$, where the optimal value of $\kappa \in (0, 1]$ is estimated using a line search routine. \\
Contrary to the fixed-point iteration, for $\mu > 0$ the root of $\mathcal{N}_{\mu}(\, \cdot \,, \, \cdot \,)$ generally differs from that of $\mu = 0$. As such, the eigenspectrum shift constitutes a regularisation rather than a stabilisation. Therefore, $\mu$ needs to be taken small and we utilise $\mu = 10^{-5}$ in practice. While in discretisations based on~\eqref{eq:EGG_newton_continuous}, the value of $\mu$ can be reduced to $\mu = 0$ in an outer loop, in practice this is usually not necessary. In fact, schemes based on~\eqref{eq:EGG_newton_continuous} converge in the vast majority of cases even for $\mu = 0$ and the stabilisation with $\mu > 0$ merely improves convergence behaviour for severely folded initial iterates.

\subsubsection{Hessian recovery approach}
\label{subsubsect:lakkispryer}
Having discussed the linearisations of the continuous variational formulation of~\eqref{eq:EGG_classical}, we can proceed to concrete discretisations. Clearly, for bases \smash{$\mathcal{V}_h \subset H^2(\hOm)$}, discretisations follow from replacing \smash{$\mathcal{U}^{\mathbf{F}} \rightarrow \mathcal{U}_h^{\mathbf{F}}$}. For compatibility with spaces \smash{$\mathcal{V}_h \subset H^1(\hOm)$}, in the following, we extend the linearisations from subsections \ref{subsect:NDF_fixed} and \ref{subsect:NDF_newton} with the weak Hessian recovery approach proposed in \cite{lakkis2011finite}. In what follows, we assume that $\tau(\, \cdot\, ) = \tau_{\text{ID}}(\, \cdot \,)$. \\
Assuming sufficient regularity of $u: \hOm \rightarrow \mathbb{R}$ and $\Phi: \hOm \rightarrow \mathbb{R}^{2 \times 2}(\hOm)$, the Hessian recovery approach is based on the following integration by parts formula:
$$\int \limits_{\hOm} H(u) \, \colon \, \Phi \, \mathrm{d} \bxi = - \int \limits_{\hOm} \nabla u \, \cdot \, \left(\nabla \cdot \Phi \right) \mathrm{d} \bxi + \int \limits_{\partial \hOm} \nabla u \cdot \left( \Phi \mathbf{n} \right) \mathrm{d} \Gamma,$$
wherein $\mathbf{n}: \partial \hOm \rightarrow \mathbb{R}^2$ denotes the outward normal vector on $\partial \hOm$, while the divergence $\nabla \cdot (\, \cdot \,)$ applied to $\Phi: \hOm \rightarrow \mathbb{R}^{2 \times 2}$ is taken row-wise. We introduce $\mathcal{U}^{\mathbf{F}} := \{ \mathbf{v} \in H^1(\hOm, \mathbb{R}^2)\, \vert \, \mathbf{v} = \mathbf{F} \text{ on } \partial \hOm \}$ and $\mathcal{W} := H^1(\hOm, \mathbb{R}^{2 \times 2})$ as well as $X := (\bx, \widehat{H}) \in \mathcal{U}^{\mathbf{F}} \times \mathcal{W}^2$ and $\Sigma := (\boldsymbol{\phi}, \Phi) \in \mathcal{U}^{\mathbf{0}} \times \mathcal{W}^2$. Analogous to~\eqref{eq:NDF_strong_operator}, we base a numerical scheme on the form $\mathcal{L}^H: \text{SPD}^{2 \times 2} \times \left(\mathcal{U}^{\mathbf{F}} \times \mathcal{W}^2 \right) \times \left(\mathcal{U}^{\mathbf{0}} \times \mathcal{W}^2 \right)$, with
\begin{align}
\label{eq:EGG_operator_Hessian}
    \mathcal{L}^H(B, X, \Sigma) & = \int \limits_{\hOm} \phi_i B \, \colon \, \widehat{H}_i \, \mathrm{d} \bxi + \int \limits_{\hOm} \left( \widehat{H}_i \, \colon \, \Phi_i + \nabla x_i \, \cdot \, \left(\nabla \cdot \Phi_i \right) \right) \mathrm{d} \bxi - \int \limits_{\partial \hOm} \nabla x_i \cdot (\Phi_i \mathbf{n}) \, \mathrm{d} \Gamma,
\end{align}
wherein we sum over repeated indices. Note that here, elements $Q \in \mathcal{W}^2$ are of the form $Q = (Q_1, Q_2) \in \mathcal{W} \times \mathcal{W}$ (and are therefore indexed in the same way as vectors). Letting, again, $A_\mu^k := A(\partial_{\bxi} \bx^k)$, the fixed-point iteration is based on
\begin{align}
    \text{find } X^{k+1} \in \mathcal{U}^{\mathbf{F}} \times \mathcal{W}^2 \quad \text{s.t.} \quad \mathcal{F}^H_{\mu}(X^{k+1}, X^k, \Sigma) = 0 \quad \forall \Sigma \in \mathcal{U}^{\mathbf{0}} \times \mathcal{W}^2,
\end{align}
where
\begin{align}
    \mathcal{F}^H_{\mu}(X^{k+1}, X^k, \Sigma) := \mathcal{L}^H(\gamma^k_\mu A_\mu^k, X^{k+1}, \Sigma) - \mu \mathcal{L}^H(\gamma^k_\mu \mathcal{I}^{2 \times 2}, X^{k}, \Sigma).
\end{align}
Following \cite{lakkis2011finite}, a discretisation replaces $\mathcal{U}^{\mathbf{F}} \rightarrow \mathcal{U}_h^{\mathbf{F}}$ and $\mathcal{W} \rightarrow \mathcal{W}_h := \mathcal{V}_h^{2 \times 2}$, with $\mathcal{V}_h \subset H^1(\hOm)$. \\
Similarly, the Newton approach is based on
\begin{align}
    \mathcal{N}^H_{\mu}(X, \Sigma) := \mathcal{L}^H(\gamma_{\mu}^{\bx} A_{\mu}^{\bx}, X, \Sigma)
\end{align}
and seeks the increment $\partial X^k \in \mathcal{U}^{\mathbf{0}} \times \mathcal{W}^2$ as in~\eqref{eq:EGG_newton_continuous} by taking the Gateaux derivative of $\mathcal{N}_{\mu}^H(\, \cdot \, , \, \cdot \,)$ with respect to its first argument. We discretise in the same way as in the fixed-point iteration. \\
The Hessian recovery approach increases the problem's cardinality from $\sim 2 \operatorname{dim}(\mathcal{V}_h )$ to $\sim 10 \operatorname{dim}(\mathcal{V}_h)$. However, we note that $\mathcal{L}^H( A(\partial_{\bxi} \bx), \, \cdot \,, \, \cdot \, )$ is nonlinear only in the first term on the right hand side of~\eqref{eq:EGG_operator_Hessian}. As such, the linearisation's bilinear form only needs to be reassembled partially and an efficient implementation can, in fact, operate on the Schur complement of the matrix's constant blocks, making the cardinality increase manageable in practice. For a more in-depth discourse on an efficient implementation, we refer to \cite{lakkis2011finite}.

\subsubsection{Rotation-free approach}
\label{subsubsect:gallistl}
This approach adopts the formulation proposed by Gallistl et al. in \cite{gallistl2017variational}. Here, we restrict ourselves to the choice $\tau(\, \cdot \,) = \tau_{\text{NS}}(\, \cdot \,)$. Furthermore, for reasons that shall become apparent shortly, we focus exclusively on the fixed-point linearisation. \\
Rather than directly solving for $\bx: \hOm \rightarrow \Om$, the rotation-free approach, in simple terms, is based on a formulation which seeks the map's Jacobian $J := \partial_{\xi} \bx$. In order for $J \in H^1(\hOm, \mathbb{R}^{2 \times 2})$ to be the gradient of a two-component function, it requires the rows of $J$ to be rotation-free for which it introduces a suitable Lagrange multiplier. We introduce the space 
$$\mathcal{W}^{\mathbf{F}} := \left \{ (\mathbf{v}_1, \mathbf{v}_2) \in H^1(\hOm, \mathbb{R}^2) \times H^1(\hOm, \mathbb{R}^2) \, \, \vert \, \, \text{the tangential trace of } \mathbf{v}_i \text{ equals } \partial_{\mathbf{t}} \mathbf{F}_i \text{ in } \hOm \,  \right \},$$
where $\partial_{\mathbf{t}} (\, \cdot \,)$ denotes the tangential derivative along $\partial \hOm$. Furthermore, we introduce 
$$Q := \left \{q \in L^2(\hOm) \, \vert \, \smallint \limits_{\partial\hOm} q \, \mathrm{d} \bxi = 0 \right \}.$$
With \smash{$Y := (\widehat{J}, \mathbf{p}) \in W^{\mathbf{F}} \times Q^2$} and $\Sigma := (\Phi, \mathbf{q}) \in \mathcal{W}^{\mathbf{0}} \times Q^2$ , the continuous formulation is based on the operator $\mathcal{L}^{\text{rot}}: \text{SPD}^{2 \times 2} \times (\mathcal{W}^{\mathbf{F}} \times Q^2) \times (\mathcal{W}^{\mathbf{0}} \times Q^2) \rightarrow \mathbb{R}$, with
\begin{align}
\label{eq:Gallistl_operator}
    \mathcal{L}^{\text{rot}}(B, Y, \Sigma) = \int \limits_{\hOm} \left( \nabla \cdot \Phi_i \right) \, B \, \colon \, \partial_{\xi} \widehat{J}_i \, \mathrm{d} \bxi + \int \limits_{\hOm} \left(\nabla \times \Phi_i \right) \mathbf{p}_i \, \mathrm{d} \bxi + \int \limits_{\hOm} \left(\nabla \times \widehat{J}_i \right) \mathbf{q}_i \, \mathrm{d} \bxi
\end{align}
and $\nabla \times \mathbf{v} := \partial_{\xi_2} \mathbf{v}_1 - \partial_{\xi_1} \mathbf{v}_2$. \\
The fixed-point linearisation reads:
\begin{align}
\label{eq:Gallistl_fixed_point}
    \text{find } Y^{k+1} \in \mathcal{W}^{\mathbf{F}} \times Q^2 \quad \text{s.t.} \quad \mathcal{F}_\mu^{\text{rot}}(Y^{k+1}, Y^k, \Sigma) = 0, \quad \forall \Sigma \in \mathcal{W}^{\mathbf{0}} \times Q^2,
\end{align}
with
\begin{align}
    \mathcal{F}_\mu^{\text{rot}}(Y^{k+1}, Y^k, \Sigma) := \mathcal{L}^{\text{rot}}(\gamma^k_\mu A^k_\mu, Y^{k+1}, \Sigma) - \mu \mathcal{L}^{\text{rot}}(\gamma^k_\mu \mathcal{I}^{2 \times 2}, Y^k, \Sigma),
\end{align}
where now $A^k_\mu := A(J^k) + \mu \mathcal{I}^{2 \times 2}$, with the $i$-th row of $J^k$ given by $\widehat{J}_i^k$ and, as before, $\gamma_\mu^k := \gamma(A^k_\mu)$. \\

\noindent For the operator utilised in the discretisation, we have to make two adjustments. Firstly, depending on $\mathbf{F}: \partial \hOm \rightarrow \partial \Om$, the discrete space $\mathcal{W}_h^{\mathbf{F}}$ will generally be empty. As such, we implement the boundary condition \textit{weakly}. Secondly, as the discrete root \smash{$Y^{k+1}_h \in \mathcal{W}_h^{\mathbf{F}} \times Q_h^2$} of~\eqref{eq:Gallistl_fixed_point} generally does not satisfy \smash{$\nabla \times \widehat{J}_{h, i} = 0$} (pointwise) for $i \in \{1, 2\}$, the discrete problem requires stabilisation. Here, we follow the stabilisation proposed in \cite{gallistl2017variational}. Defining
$$\epsilon(B) := \min \limits_{\hOm} \frac{\operatorname{tr}(B)^2}{\left \| B \right \|^2_F} - 1, \quad \lambda(\, \cdot \,) := \frac{2 + \sqrt{\alpha \epsilon(\, \cdot \,)}}{2} \quad \text{and} \quad \sigma(\, \cdot \,) := \sqrt{1 - \tfrac{\lambda(\, \cdot \,)}{2}},$$
with $0 < \alpha < 1$, the stabilisation consists of adding
$$ \mathcal{K}_{\text{stab}}(B, \widehat{J}, \Phi) := \sigma(B) \int \limits_{\hOm} \left(\nabla \times \Phi_i \right) \left(\nabla \times \widehat{J}_i \right) \quad \text{(summing over repeated indices)}$$
to the right-hand-side of~\eqref{eq:Gallistl_operator}. Here, we use $\alpha = 0.9$ and in practice, we approximate $\epsilon(B) \approx \epsilon_h(B)$ by taking the minimum over all evaluations in the abscissae of the quadrature scheme used to compute the integrals. Finally, the stabilised operator $\mathcal{L}^{\text{rot}, \text{stab}}: \text{SYM}^{2 \times 2} \times \mathcal{W} \times Q^2 \rightarrow \mathbb{R}$, along with the weak imposition of the boundary data reads:
\begin{align}
    \mathcal{L}_{\eta}^{\text{rot}, \text{stab}}(B, Y, \Sigma) := \mathcal{L}^{\text{rot}}(B, Y, \Sigma) + \mathcal{K}_{\text{stab}}(B, \widehat{J}, \Phi) + \eta \sum \limits_{L_j \in \Gamma^B} \frac{1}{h_j} \int \limits_{L_j} \left(\partial_{t} \mathbf{F} - \widehat{J} \, \hat{\mathbf{t}} \right) \cdot \left( \Phi \hat{\mathbf{t}} \right) \mathrm{d} \Gamma,
\end{align}
where $\mathcal{W} := H^1(\hOm, \mathbb{R}^{2 \times 2})$ while $h_i$ denotes the average diameter of all knot spans on $L_i \subset \partial \hOm$ and $\hat{\mathbf{t}}$ denotes the unit tangent along $\partial \hOm$. The factor $\eta > 0$ needs to be taken sufficiently large and in practice, we utilise $\eta = 10^3$. \\
The discrete problem is subject to the same inf-sup condition as the Stokes problem \cite{gallistl2017variational}. Here, we utilise the subgrid space pair \cite{bressan2018inf}. If $\mathcal{W}_h = \mathcal{V}_h^{2 \times 2}$ and $Q_h$ is constructed by a modification of some finite-dimensional $\mathcal{U}_h \subset H^1(\hOm)$ (to incorporate the zero average condition), this implies that $\mathcal{V}_h$ is uniformly $h$-refined with respect to $\mathcal{U}_h$. In practice, the space $\mathcal{U}_h$ results from removing every other knot in the knotvectors utilised to construct the primal space $\mathcal{V}_h$. As such, the cardinality of the problem is $\sim 4.5 \times \operatorname{dim}(\mathcal{V}_h)$ and computational efficiency can be greatly improved by only reassembling the nonlinear part of the fixed-point iteration's global matrix equation. The scheme, in its current form, is not compatible with Newton's method due to the required stabilisation.\\
It remains to be said that the map $\bx: \hOm \rightarrow \Om$ can be recovered by solving:
$$\text{find } \bx \in \mathcal{U}^{\mathbf{F}}, \quad \text{s.t.} \quad \int \limits_{\hOm} \left(\partial_{\bxi} \bx - \widehat{J} \, \right) \, \colon \, \partial_{\bxi} \boldsymbol{\phi} \, \mathrm{d} \bxi = 0, \quad \forall \boldsymbol{\phi} \in \mathcal{U}^{\mathbf{0}},$$
and similarly for the discrete counterpart.

\subsubsection{$C^0$-DG approach}
\label{subsubsect:blechschmidt}
Having presented two approaches in mixed form, we now proceed to an approach based on the $C^0$-DG formulation from \cite{blechschmidt2021error}. A $C^0$-DG formulation is particularly appealing as it completely avoids auxiliary variables. Furthermore the discrete basis is (by assumption) sufficiently regular in the interior of patches for penalisation to be restricted to the interior patch facets $\gamma_{ij} \in \Gamma^I$. As opposed to mixed formulations, the $C^0$-DG approach employs the patchwise exact Hessian while weakly imposing continuity of the map's Jacobian across interior interfaces. Here, we restrict ourselves to the choice $\tau( \, \cdot \,) = \tau_{\text{NS}}(\, \cdot \,)$. The operator from~\eqref{eq:NDF_strong_operator} is adjusted as follows: $\mathcal{L} \rightarrow \mathcal{L}^{\text{DG}}_{\eta}$, with $\mathcal{L}_{\eta}^{\text{DG}}: \text{SPD}^{2 \times 2} \times \mathcal{U}^{\mathbf{F}} \times \mathcal{U}^{\mathbf{0}} \rightarrow \mathbb{R}$ satisfying
\begin{align}
    \mathcal{L}^{\text{DG}}_\eta(B, \bx, \boldsymbol{\phi}) := \sum \limits_{k=1}^{N_p} \int \limits_{\hOm_k} \Delta \phi_i B \, \colon \, H(x_i) \, \mathrm{d} \bxi + \eta \sum_{\gamma_{jl} \in \Gamma^I} \frac{1}{h(\gamma_{jl})} \int \limits_{\gamma_{jl}} [\![ \nabla x_i ]\!] \, \colon \, [\![ \nabla \phi_i ] \! ] \, \mathrm{d} \Gamma.
\end{align}
Here, $[ \! [ \, \mathbf{v} \, ] \! ]$ denotes the (entry-wise) jump term of $\mathbf{v} \otimes \mathbf{n} \in L^2(\gamma_{ij}, \mathbb{R}^{2 \times 2})$, with $\mathbf{n}$ the unit outer normal on $\gamma_{ij}$ in arbitrary but fixed direction while $h(\gamma_{ij})$ denotes the average diameter of all knot spans on the facet $\gamma_{ij}$. The penalisation parameter $\eta > 0$ has to be chosen sufficiently large. In practice, facing geometries with characteristic length scales of $\mathcal{O}(1)$, we utilise $\eta = 10$. \\
The fixed-point iteration as well as the Newton approach are adapted to this formulation simply by replacing $\mathcal{L} \rightarrow \mathcal{L}^{\text{DG}}_\eta$ in $\mathcal{F}_{\mu}(\, \cdot \,, \, \cdot \,, \, \cdot \,)$ and $\mathcal{N}_{\mu}(\, \cdot \,, \, \cdot \,)$, respectively, which is not repeated here for the sake of brevity.

\subsection{Regularised weak form discretisation}
\label{subsect:weak_form_discretisation}
This scheme is based on the weak inverse harmonicity requirement from~\eqref{eq:inverse_elliptic_weak}. In what follows, we let $\mathcal{U}^{\mathbf{F}} = \{ \mathbf{v} \in H^1(\hOm, \mathbb{R}^2) \, \vert \, \mathbf{v} = \mathbf{F} \text{ on } \hOm \}$ while \smash{$\mathcal{U}^{\mathbf{F}}_{\text{bij}} := \{\mathbf{v} \in \mathcal{U}^{\mathbf{F}} \, \vert \, \mathbf{v} \text{ is uniformly nondegenerate} \}$}. Noting that $\bx^{-1}(\bxi) = \bxi$ in $\hOm$, a pullback of the weak inverse harmonicity requirement leads to:
\begin{align}
\label{eq:inverse_elliptic_weak_pullback}
    \text{find } \bx \in \mathcal{U}^{\mathbf{F}}_{\text{bij}} \quad \text{s.t.} \quad \mathcal{L}^{\text{W}}(\bx, \boldsymbol{\phi}) = 0, \quad \forall \boldsymbol{\phi} \in \mathcal{U}^{\mathbf{0}},
\end{align}
with $\mathcal{L}^{\text{W}}: \mathcal{U}^{\mathbf{F}}_{\text{bij}} \times \mathcal{U}^{\mathbf{0}} \rightarrow \mathbb{R}$ given by
\begin{align}
\label{eq:inverse_elliptic_weak_pullback_operator}
    \mathcal{L}^{\text{W}}(\bx, \boldsymbol{\phi}) := \int \limits_{\hOm}\frac{\partial_{\bxi} \boldsymbol{\phi} \, \colon \, A(\partial_{\bxi} \bx)}{\det J} \mathrm{d} \bxi
\end{align}
and $A(\, \cdot \,)$ as in~\eqref{eq:A_CC}. The formulation based on~\eqref{eq:inverse_elliptic_weak_pullback_operator} can be regarded as the Galerkin method of the weak inverse harmonicity formulation while Winslow's original approach constitutes the associated Ritz-Galerkin method. \\
The appearance of $\det J$ in the denominator, as in Winslow's original approach, prohibits the substitution of degenerate maps, hence the requirement to restrict the domain of $\mathcal{L}^W(\cdot, \, \boldsymbol{\phi})$ to $\mathcal{U}^{\mathbf{F}}_{\text{bij}}$ instead of $\mathcal{U}^{\mathbf{F}}$. However, this limits the scope of algorithms based on~\eqref{eq:inverse_elliptic_weak_pullback} to improving the parametric quality of an already (uniformly) nondegenerate map. In order to attenuate this harsh requirement, we employ the regularisation proposed in \cite{ji2022penalty, wang2021smooth}, whose original purpose was regularising the Winslow function by replacing
\begin{align}
\label{eq:jacdet_regularisation}
    \det J \rightarrow \mathcal{R}_{\varepsilon}(\det J), \quad \text{where} \quad \mathcal{R}_{\varepsilon}(x) := \frac{x + \sqrt{4 \varepsilon^2 + x^2}}{2}.
\end{align}
We denote the regularised operator by $\mathcal{L}_\varepsilon^{\text{W}}(\cdot, \, \cdot)$, whose domain is restored to $\mathcal{U}^{\mathbf{F}} \times \mathcal{U}^{\mathbf{0}}$. The asymptotic behaviour of~\eqref{eq:jacdet_regularisation} reads
$$\lim \limits_{x \rightarrow -\infty} \mathcal{R}_\varepsilon(x) = 0 \quad \text{and} \quad \lim \limits_{x \rightarrow \infty} \mathcal{R}_{\varepsilon}(x) = x, \quad \text{with} \quad \mathcal{R}_{\varepsilon}(0) = \varepsilon.$$
For $\varepsilon >0$, $\mathcal{R}_{\varepsilon} \in C^1(\mathbb{R})$ and $\mathcal{R}_\varepsilon(x) > 0 \, \forall x \in \mathbb{R}$. As such, the regularisation can be combined with a gradient-based algorithm acting on~\eqref{eq:inverse_elliptic_weak_pullback}, such as Newton's method, again replacing $\mathcal{U}^{\mathbf{F}} \rightarrow \mathcal{U}^{\mathbf{F}}_h$ for a discretisation. Heuristically, $\varepsilon = 10^{-4}$ is a reliable choice as it dramatically increases the radius of convergence and even extends it into the set of degenerate initial iterates in practice. The regularisation is a convenient means to urge a globalised Newton-based root finder to decrease the size of the Newton step in case the updated iterate accidentally leaves the set of NDG maps. In the absence of regularisation, the division by zero typically causes a numerical algorithm to diverge, even when initialised with an NDG initial map. \\
The value of $\varepsilon > 0$ can be reduced to $\varepsilon = 0$ in an outer loop, in which case it is almost guaranteed that the resulting map is nondegenerate since for $\varepsilon \rightarrow 0$, $\mathcal{R}_{\varepsilon}^{-1}(\, \cdot \,)$ acts as a barrier term, as in Winslow's original approach.
While the discrete root of~\eqref{eq:inverse_elliptic_weak_pullback} substituted into the Winslow functional typically yields a value slightly greater than Winslow's global minimiser over $\mathcal{U}_h^{\mathbf{F}}$, we have noticed~\eqref{eq:inverse_elliptic_weak_pullback} to converge faster and more reliably than Winslow's original (regularised) formulation. Furthermore, it is plausible to assume that for $\varepsilon \rightarrow 0$,~\eqref{eq:inverse_elliptic_weak_pullback} has a unique root, while the discretisation of Winslow's approach, may produce local minima.\\
Compared to the NDF discretisations (cf. Subsection \ref{subsect:NDF_discretisations}), the radius of convergence is small. As such, the method is best initialised with one of the NDF discretisations' solutions, for which it typically converges in no more than $5$ Newton iterations. Furthermore, it provides a convenient way of untangling a degenerate map produced by an NDF discretisation without the need to recompute the map over a refined space. Convergence failure of~\eqref{eq:inverse_elliptic_weak_pullback} may furthermore indicate that the set \smash{$\mathcal{U}_{h, \text{bij}}^{\mathbf{F}}$} is empty, thus making refinement mandatory.

\subsection{Boundary correspondence requirements}
\label{sect:boundary_correspondence_requirements}
\begin{figure}[h!]
\centering
    \includegraphics[scale=0.8]{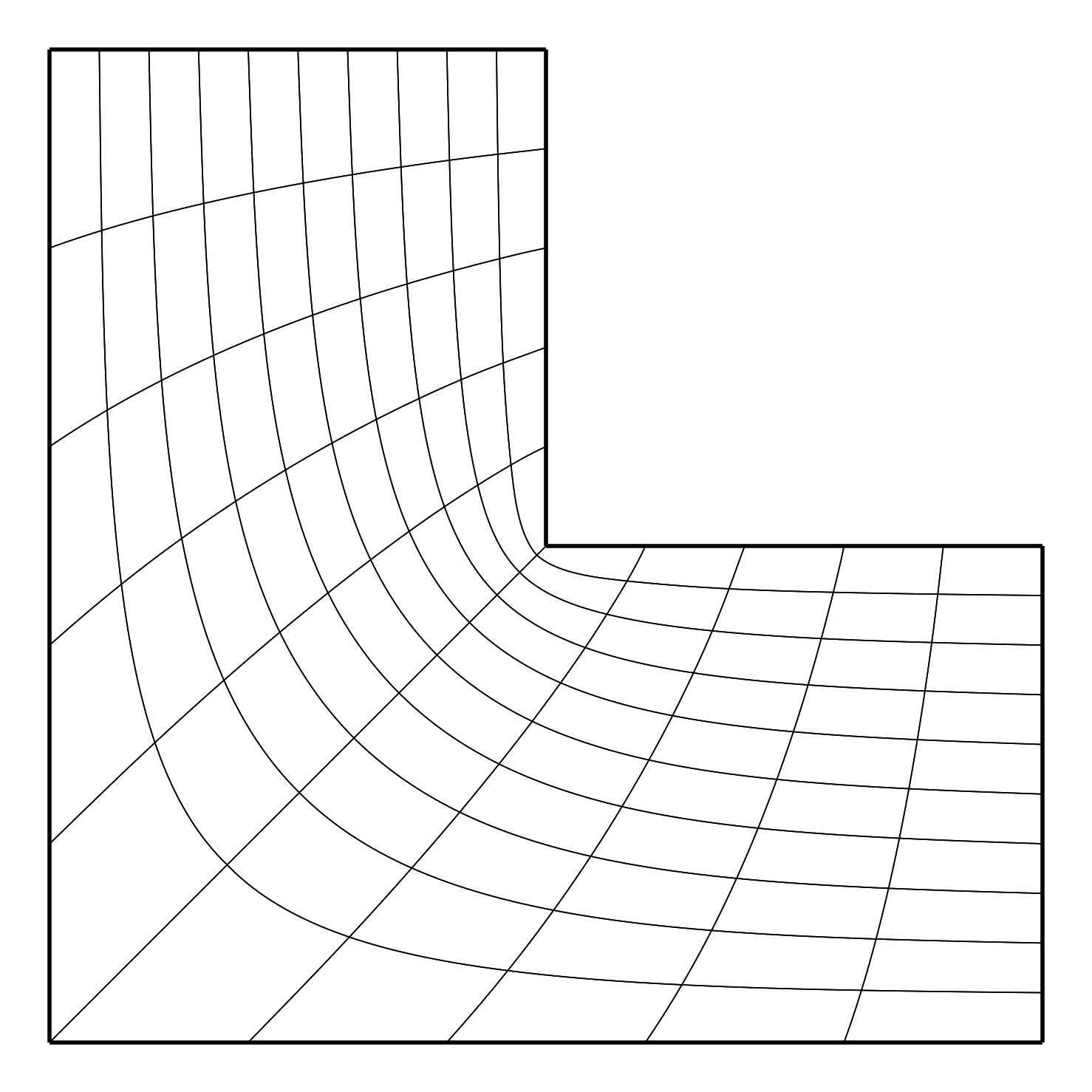}
    \caption{The L-bend geometry computed from a two-patch covering of the unit interval $\hOm = (0, 1)^2$. This geometry provides an example for singularities created by concave corners. The eigenvalue of the map's metric tensor that corresponds to the direction transversal to the concave corner is given for various refinement levels $h: 0.0118$, $h/2: 0.00729$, $h/4: 0.00453$ $h/8: 0.00284$. Each map has been computed by minimising the Winslow function~\eqref{eq:Winslow_pullback}.}
    \label{fig:L_bend_singular}
\end{figure}

Without aspirations to provide formal proofs, this section discusses the requirements that $\hOm$, $\Om$ and the boundary correspondence $\mathbf{F}: \hOm \rightarrow \Om$ have to satisfy in order for the variational formulations of Section \ref{sect:numerical_schemes} to be well-posed. \\
As $A(\partial_{\bxi} \bx)$, as defined in~\eqref{eq:EGG_classical}, has the same characteristic polynomial as the map's metric tensor $G_{ij} = \partial_{\bxi_i} \bx \cdot \partial_{\bxi_j} \bx$, it is plausible to assume that a necessary condition for well-posedness of the NDF discretisations is that the harmonic map $\bx^{-1}: \Om \rightarrow \hOm$ satisfies
\begin{align}
\label{eq:requirement_well_posedness}
    0 < \inf \limits_{\bx \in \overline{\Omega}} \, \det J(\bx^{-1}) \leq \sup \limits_{\bx \in \overline{\Omega}} \, \det J(\bx^{-1}) < \infty.
\end{align}
While Theorem \ref{thrm:RKC} guarantees that $\bx^{-1}: \Om \rightarrow \hOm$ is diffeomorphic in $\Om$, it provides no guarantee that the map is differentiable on the closure $\overline{\Om}$ of $\Om$. Failure to satisfy~\eqref{eq:requirement_well_posedness} will cause $A(\, \cdot \,)$ to no longer be uniformly elliptic in the exact solution, which may render the problem ill-posed in this case. While this section's algorithms may nevertheless succeed in finding the discrete problem's root, we may expect to encounter conditioning issues in the linearisation's bilinear forms in a refinement study. For $\det J(\bx^{-1})$ to stay uniformly bounded on the closure, we require that $\mathbf{F}^{-1}: \partial \Om \rightarrow \partial \hOm$ maps the convex corners of $\partial \Om$ onto the convex corners of $\partial \hOm$ while the smooth segments of $\partial \Om$ are mapped onto straight line segments of $\partial \hOm$ where we furthermore require that $\partial_t \mathbf{F}^{-1}$ be continuous in the vertices that are mapped onto vertices $\mathbf{v}_{ij}\in \partial \hOm$ connecting two sides $L_i$ and $L_j$ of $\partial \hOm$ without creating a corner (and similarly for $\mathbf{F}: \partial \hOm \rightarrow \partial \Om$). \\
Clearly, mapping a straight line segment of $\overline{L}_i \cup \overline{L}_j \subset \partial \hOm$ onto two sides $\overline{C}_i \cup \overline{C}_j \subset \partial \Om$ with a convex corner in the shared vertex will cause $\det J(\bx) \rightarrow \infty$ in the vertex connecting $L_i$ and $L_j$. Similarly, mapping the same vertex onto a vertex of $\partial \Om$ that creates a concave corner will cause a singularity, i.e., $\det J(\bx) \rightarrow 0$. While discrete approximations typically remain UNDG, this behaviour will be observable in a refinement study, see Figure \ref{fig:L_bend_singular}. Since the weak-form operator $\mathcal{L}^W_{\varepsilon}: \mathcal{U}^{\mathbf{F}} \times \mathcal{U}^{\mathbf{0}} \rightarrow \mathbb{R}$ may map into $\mathbb{R}$ despite division by zero on $\partial \hOm$, even for $\varepsilon \rightarrow 0$, it may allow for boundary correspondences that exhibit less regularity. Notwithstanding the well-posedness of the formulation, we regard discrete maps $\bx_h: \hOm \rightarrow \Om$ that approximate a map $\bx: \hOm \rightarrow \Om$ with a singularity on $\partial \hOm$ as undesirable from a numerical perspective. \\
In what follows, we refer to a boundary correspondence $\mathbf{F}: \partial \hOm \rightarrow \partial \Om$ that satisfies above requirements as a \textit{diffeomorphic boundary correspondence}. Section \ref{sect:control_mechanisms} augments this section's formulations with a mechanism that allows for the creation of diffeomorphic boundary correspondences even when $\partial \Om$ has no corners.

\subsection{Choosing an initial guess}
\label{subsect:initial_guess}
Both the NDF discretisations (cf. Subsection \ref{subsect:NDF_discretisations}) and the weak form discretisation (cf. Subsection \ref{subsect:weak_form_discretisation}) are nonlinear and need to be initialised with a suitable initial map $\bx_h^0$. On a single patch, NDF-type discretisations are typically initialised using the bilinearly-blended Coons' patch approach. However, since blending leads to potentially complicated polynomial constructions in the multipatch case, we initialise the NDF discretisations by a map $\bx^0$ whose components are harmonic in $\hOm$ (with the prescribed boundary correspondence). The map $\bx^0$ is then approximated by discretising the Laplace equation in the usual way. As $\Om$ is generally nonconvex, the initial map is typically degenerate. With this initial guess, the Newton schemes reliably converge after typically $5$ iterations, while the fixed-point iteration requires $\sim 15$ iterations. \\
The weak form discretisation has a smaller radius of convergence and using a harmonic map (in $\hOm$) as initialisation typically leads to convergence failure. As such, this scheme is initialised with the solution of an NDF scheme, for which it reliably converges after typically $5$ Newton iterations (using $\epsilon = 10^{-4}$ in the regularisation). The value of $\varepsilon$ can then be gradually reduced in an outer loop. In practice, this is rarely necessary.

\subsection{Numerical Experiments}
\begin{figure}[h!]
\centering
\begin{subfigure}[t]{0.25 \textwidth}
    \includegraphics[width=\textwidth]{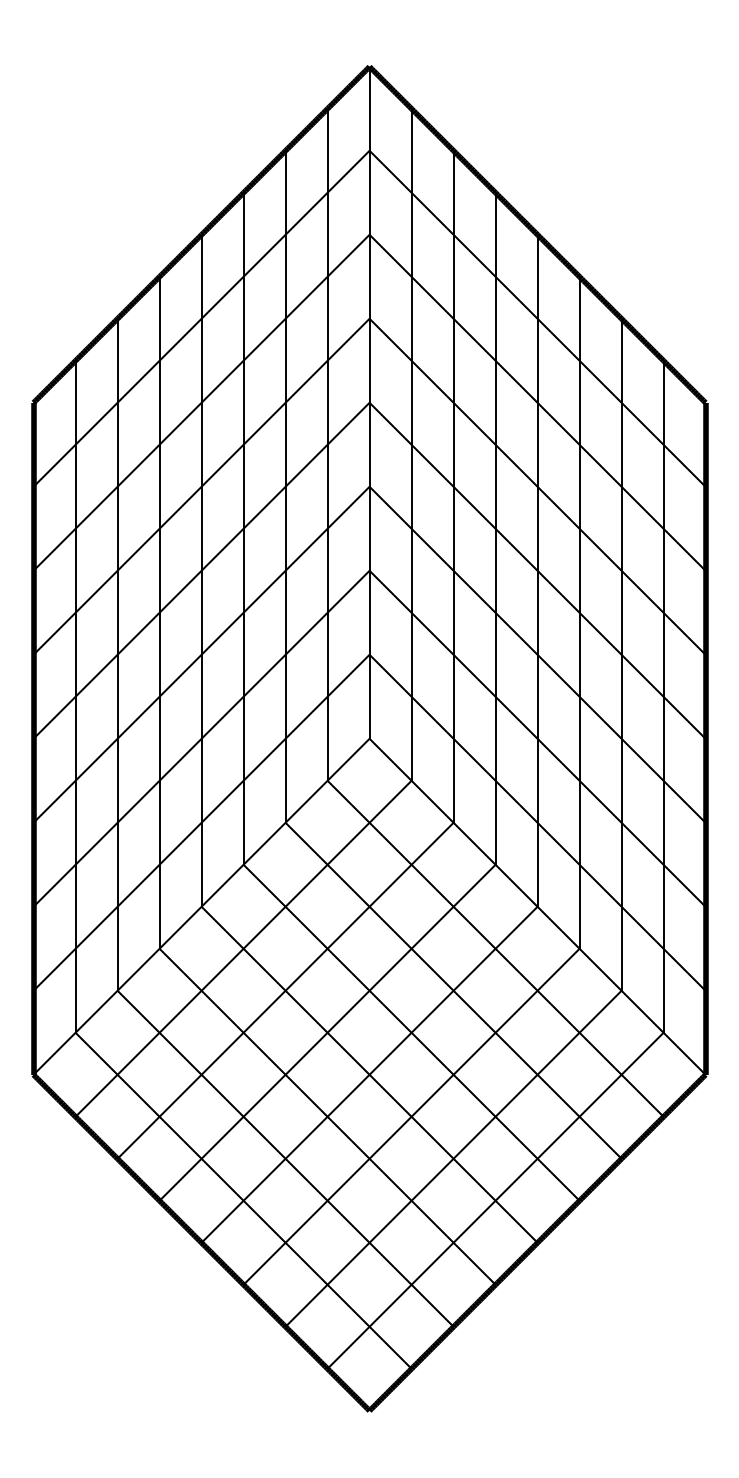}
    \caption{The multipatch covering of $\hOm$ along with the uniform knotspans.}
    \label{fig:bat_benchmark_param}
\end{subfigure}%
\hspace*{2cm}
\begin{subfigure}[t]{0.26 \textwidth}
    \includegraphics[width=\textwidth]{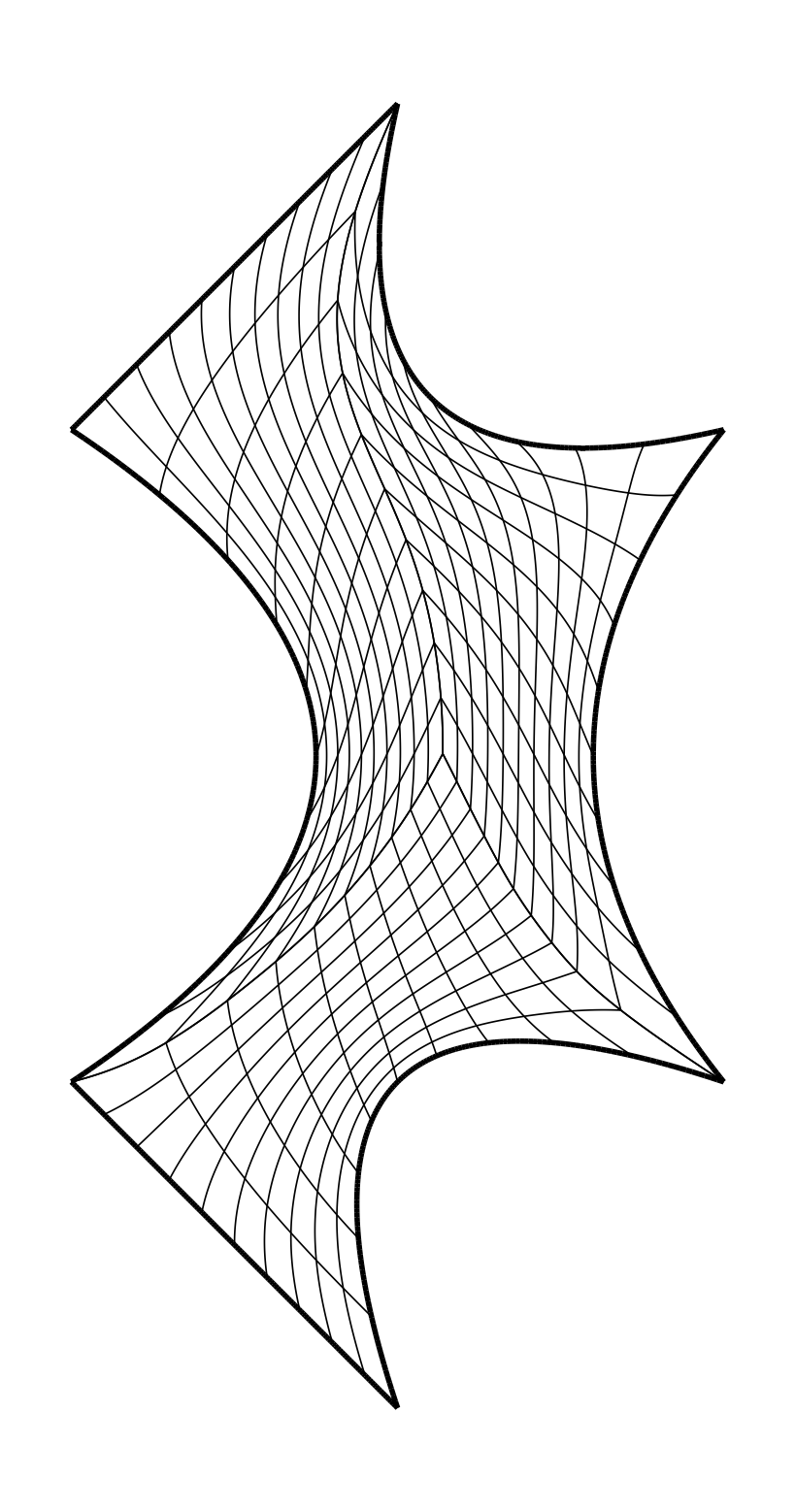}
   \caption{Uniform plot of the harmonic map at refinement level $0$ computed with the $C^0$-DG approach from Section \ref{subsubsect:blechschmidt}.}
   \label{fig:bat_benchmark_geom}
\end{subfigure}
\caption{The parametric domain and the approximation of a harmonic map for a bat-shaped geometry.}
\label{fig:bat_benchmark}
\end{figure}

In this section, we apply the algorithms from Section \ref{sect:numerical_schemes} to a number of benchmark test cases to experimentally determine the scheme's convergence rates. All schemes, as well as all control mechanisms from Section \ref{sect:control_mechanisms} have been implemented in the open-source finite element library \textit{Nutils} \cite{nutils7}. \\
As a first experiment we are considering the bat shaped geometry from Figure \ref{fig:bat_benchmark_geom} along with the parametric domain from Figure \ref{fig:bat_benchmark_param}.
\noindent The geometry is a piecewise $C^\infty$ curvilinear polygon whose sides $C_i$ are quadratic polynomials. We therefore expect the harmonic map $\bx^{-1}: \Om \rightarrow \hOm$ to satisfy $\bx^{-1} \in H^s(\Om, \mathbb{R}^2)$ with $2 \leq s \leq 3$. The boundary correspondence \smash{$\mathbf{F}: \partial \hOm \rightarrow \partial \Om$} (suitably extended into the interior) can be expressed exactly in any finite-dimensional space $\mathcal{V}_h$ with polynomial degree $p \geq 2$, thanks to its piecewise quadratic nature. We are estimating the convergence rate for all three NDF discretisations (cf. Subsection \ref{subsect:NDF_discretisations}) as well as the weak form discretisation from Subsection \ref{subsect:weak_form_discretisation} and Winslow's original approach (cf. Subsection \ref{subsect:harmonic_maps}). The NDF discretisations are initialised with the forward-Laplace initial guess (cf. Subsection \ref{subsect:initial_guess}) while the weak form discretisation is initialised with the $C^0$-DG approach's solution and Winslow's minimisation with the weak form discretisation's solution. This is repeated for several levels of $h$-refinement, each splitting each univariate knotspan in half (without changing the boundary correspondence). For the rotation-free approach, we perform a fixed-point iteration while all other linearisations are based on Newton's method. The fixed point linearisation converges after typically $12$ iterations while the Newton approach requires typically $5$ for the NDF discretisations and another $3$ for the weak form discretisation. The convergence rate is estimated in the $H^1(\hOm)$-norm. Denoting three consecutive solutions by $\bx_h$, $\bx_{h/2}$ and $\bx_{h/4}$, respectively, the convergence rate is estimated as
$$
\kappa \approx \log_2 \left( \frac{\left \| \bx_h - \bx_{h/2} \right \|_{H^1(\hOm)}}{\left \| \bx_{h/2} - \bx_{h/4} \right \|_{H^1(\hOm)}} \right),
$$
and we utilise the last three levels of refinement to estimate $\kappa$ in the above. We perform a refinement study assigning a uniform knotvector without internal knot repetitions to each interior and boundary facet whereby the coarsest knotvector contains three interior knots.
\begin{figure}[h!]
\centering \includegraphics[width=0.35 \textwidth]{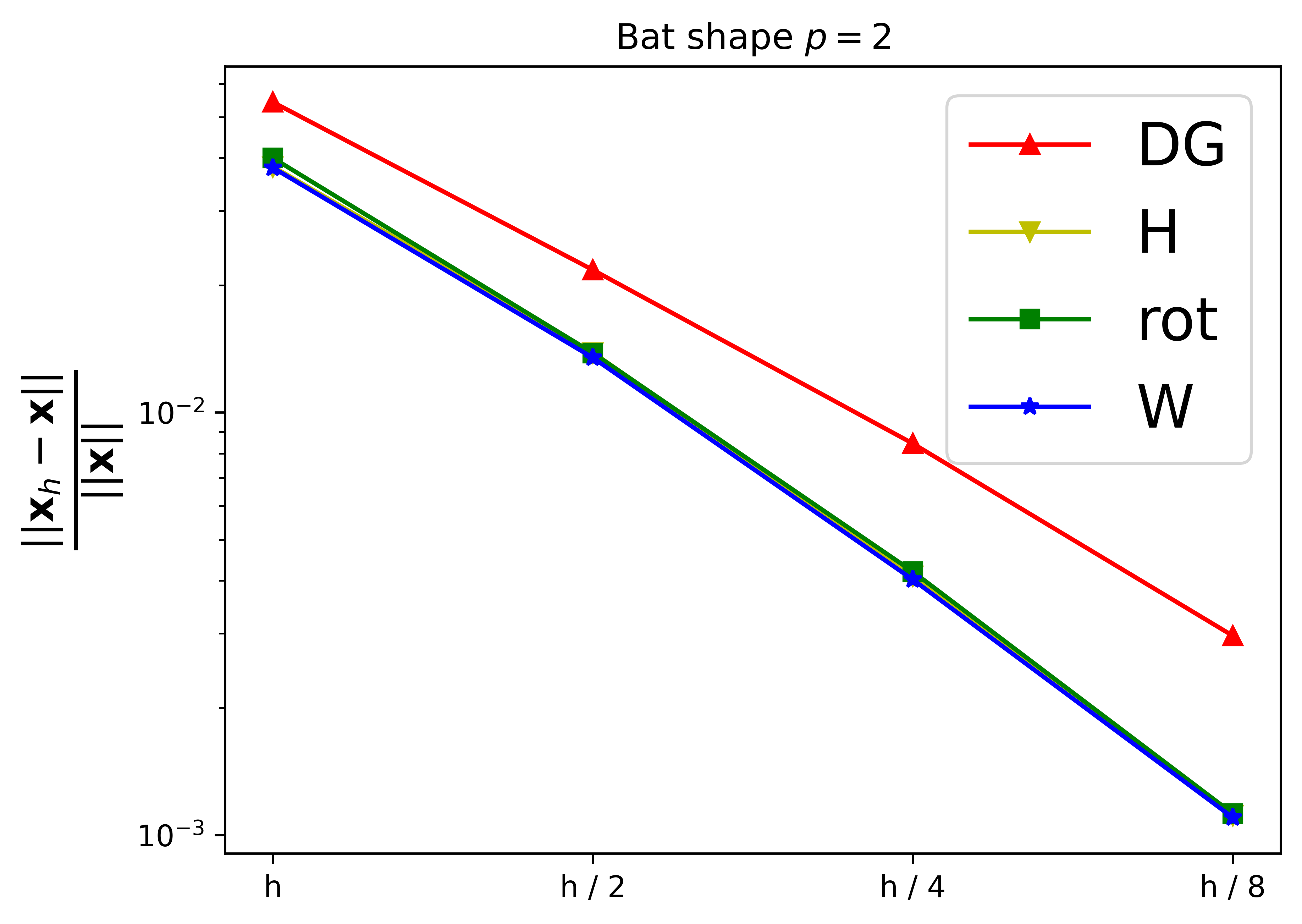} \hspace*{1cm}\includegraphics[width=0.35 \textwidth]{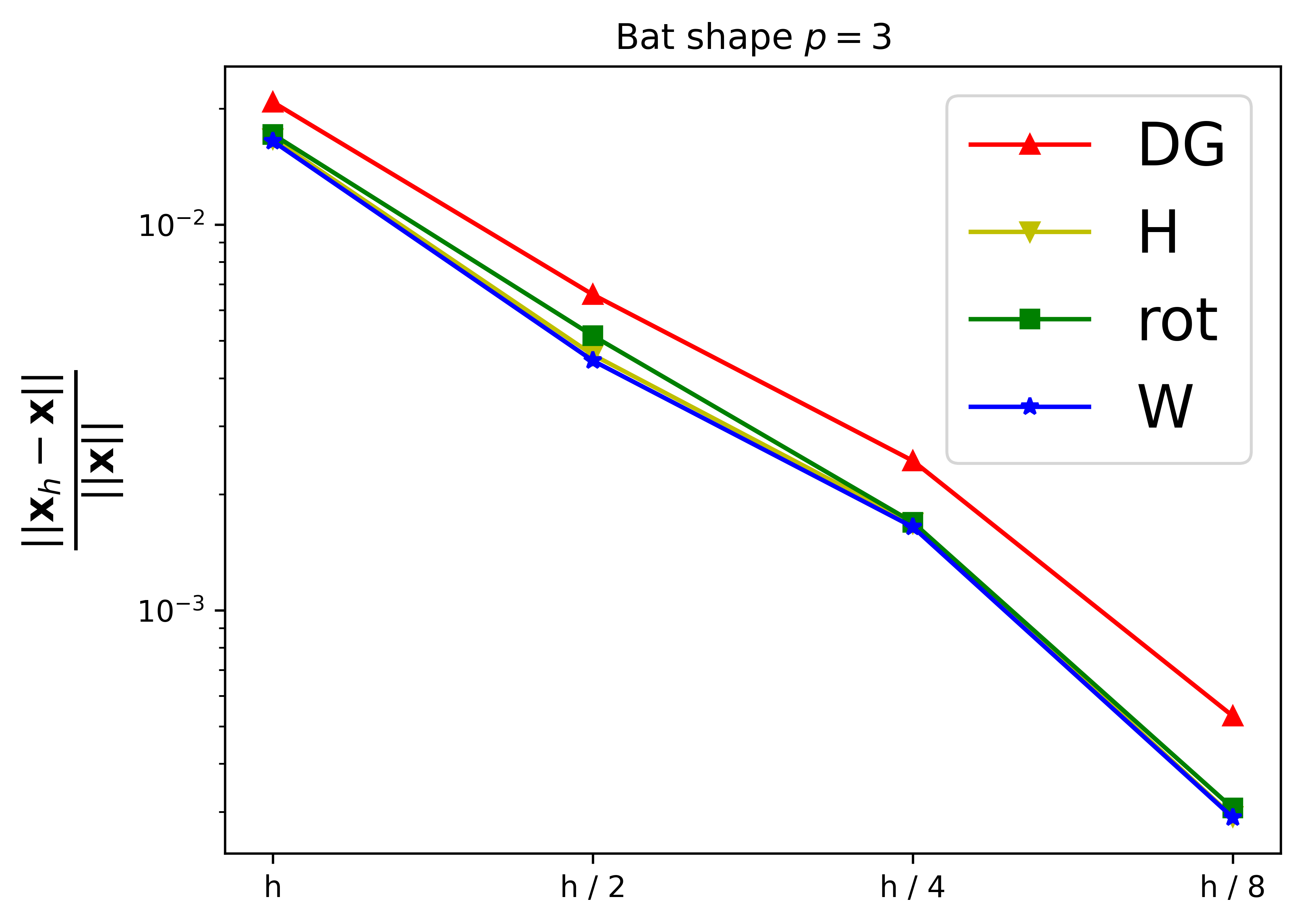} \\[0.5cm]
\includegraphics[width=0.35 \textwidth]{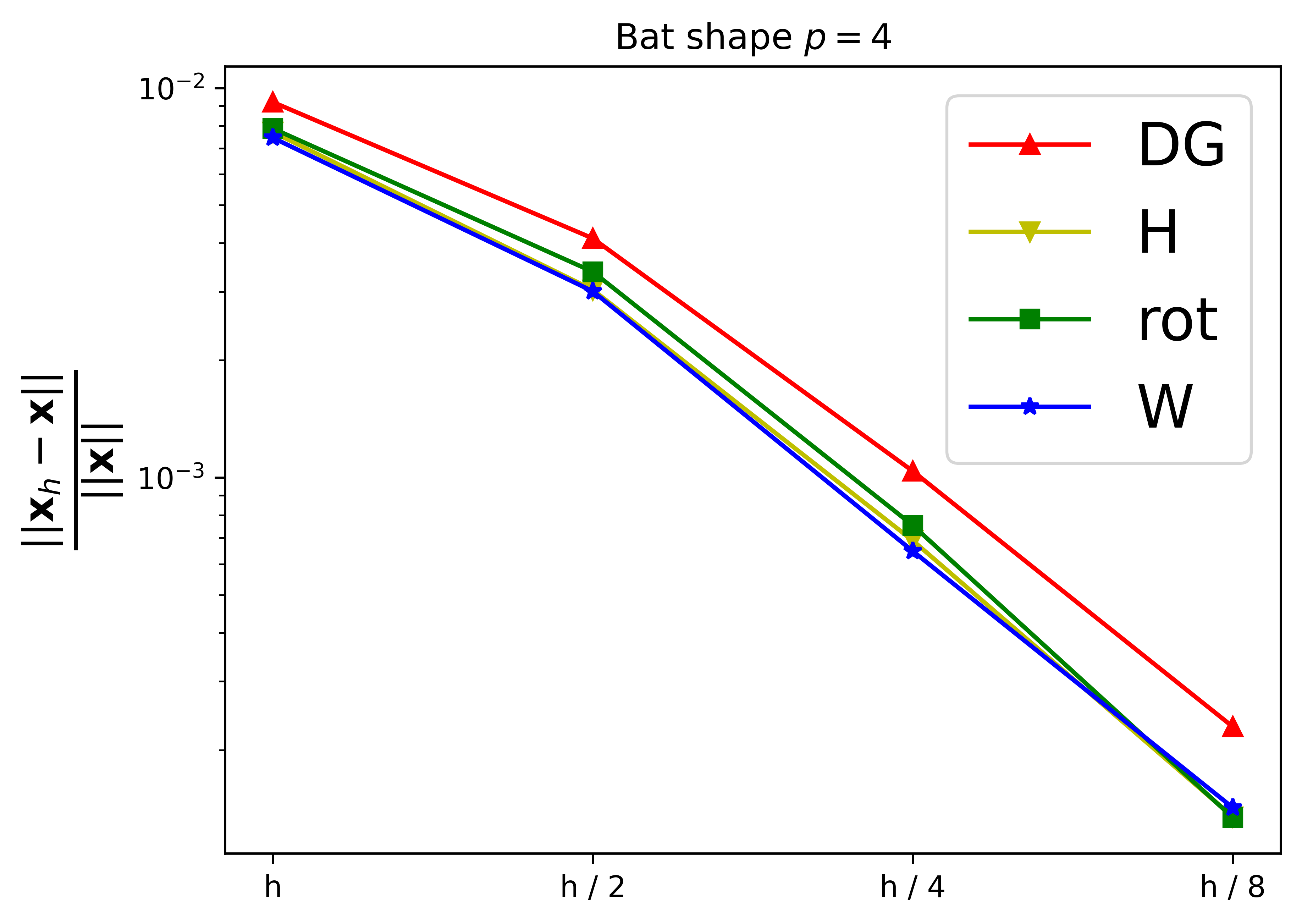} \hspace*{1cm} \includegraphics[width=0.35 \textwidth]{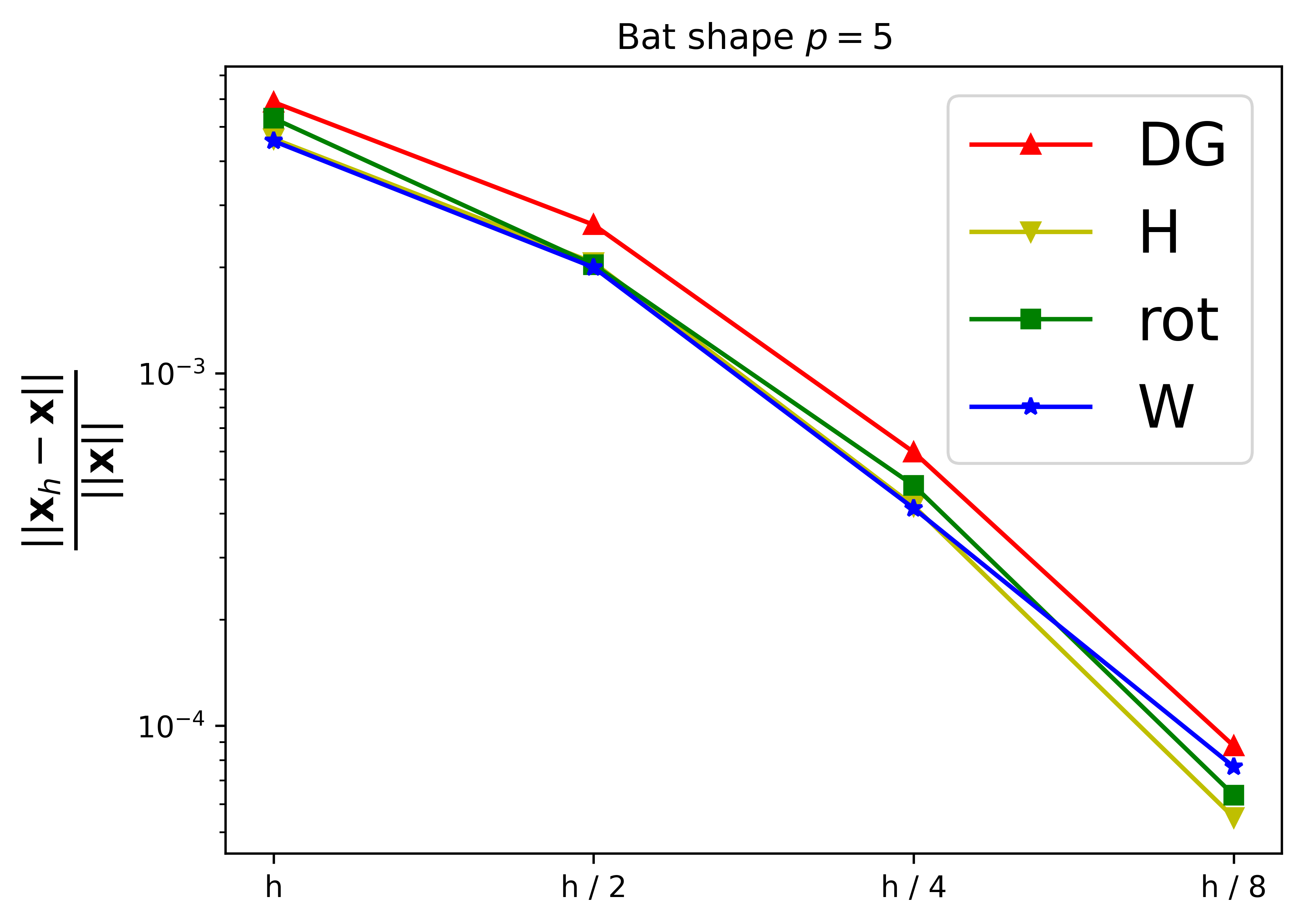}
\caption{Convergence behaviour of the various discretisations applied to the bat shaped geometry from Figure \ref{fig:bat_benchmark} for different values of the polynomial degree $p$. Here the labels, in the order of appearance, refer to the $C^0$-DG (Section \ref{subsubsect:blechschmidt}), Hessian-recovery (Section \ref{subsubsect:lakkispryer}) and the rotation-free approach (Section \ref{subsubsect:gallistl}) as well as the regularised weak-form discretisation (Section \ref{subsect:weak_form_discretisation}).}
\label{fig:convergence_plot_bat_shape}
\end{figure}

\begin{figure}[h!]
\centering \includegraphics[width=0.4 \textwidth]{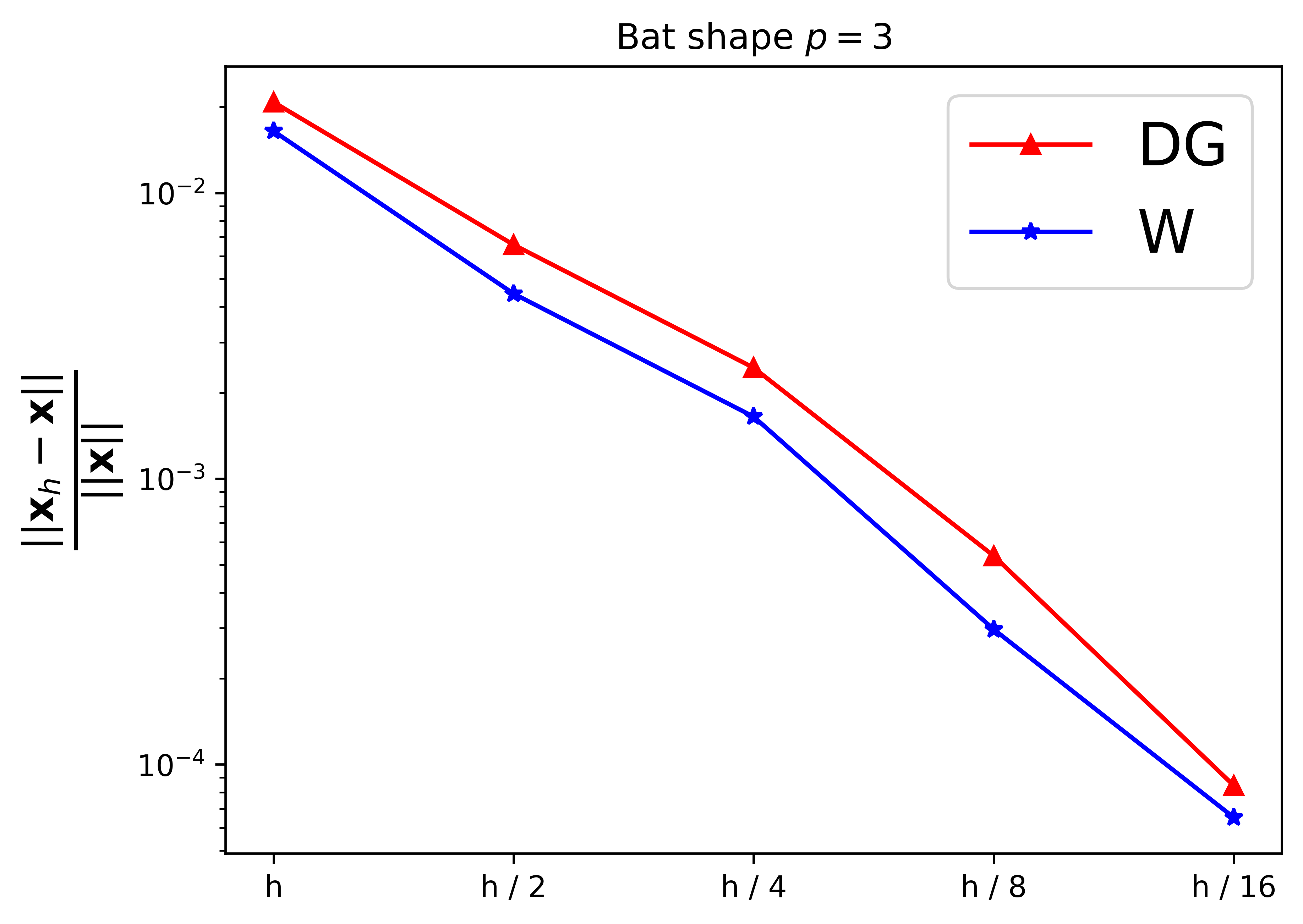}
\caption{Convergence behaviour of the $C^0$-DG and the weak form discretisation with one additional level of refinement compared to Figure \ref{fig:convergence_plot_bat_shape}.}
\label{fig:convergence_plot_bat_shape_extended}
\end{figure}

\begin{table}[h!]
\renewcommand{\arraystretch}{1.5}
\centering
\begin{tabular}{c|c|c|c|c}
 $p$ & 2 & 3 & 4 & 5 \\
 \hline
$\kappa \left(\mathcal{L}^{\text{DG}}\right)$ & 1.37 & 1.36 & 2.08 & 2.04 \\
\hline
$\kappa\left(\mathcal{L}^H \right)$ & 1.68 & 1.32 & 2.09 & 2.22  \\
\hline
$\kappa \left( \mathcal{L}^{\text{rot}} \right)$ & 1.66 & 1.61 & 2.13 & 2.01  \\
\hline
$\kappa \left(\mathcal{L}^{\text{W}} \right)$ & 1.68 & 1.31 & 2.19 & 2.22 \\
\hline
$\kappa \left(\text{Winslow} \right)$ & 1.75 & 1.32 & 2.16 & 2.31 \\
\end{tabular}
\caption{Approximate convergence rates for the bat shape geometry for various values of the polynomial degree $p$.}
\label{tab:refinement_study_bat_shape_p2to5}
\end{table}
Table \ref{tab:refinement_study_bat_shape_p2to5} contains the approximate convergence rates of this section's discretisations for various values of the polynomial degree $p \geq 2$ while Figure \ref{fig:convergence_plot_bat_shape} plots the relative $H^1$-norm discrepancy between the various refinement level's solutions and the exact solution which is approximated by minimising the Winslow function over a spline space that is one level of refinement ahead of the maximum refinement level in the plots (here: $h/16$). The table and plots suggest that all discretisations perform similarly well, a notable exception being the $C^0$-DG approach from Section \ref{subsubsect:blechschmidt} which consistently produces the largest relative $H^1$-norm discrepancy. We note however that the $C^0$-DG approach is the computationally least expensive method that can be initialised with a degenerate map. The approximate convergence rates improve for larger values of $p \geq 2$, eventually reaching saturation for $p = 5$. This is not surprising given that it is plausible to assume that the largest attainable $H^1(\hOm, \mathbb{R}^2)$ convergence rate is bounded by the value of $2 \leq s \leq 3$ associated with the harmonic map $\bx^{-1} \in H^s(\Om, \mathbb{R}^2)$. A notable exception is the outcome for $p = 3$, which consistently ranks below all other choices. Figure \ref{fig:convergence_plot_bat_shape_extended} depicts the convergence behaviour of the $C^0$-DG and weak form discretisations with one additional level of refinement. Applying the convergence rate estimator to the last three consecutive levels of refinement yields $\kappa(\mathcal{L}^{\text{DG}}) \approx 2.18$, $\kappa(\mathcal{L}^{\text{W}}) \approx 2.43$ and $\kappa(\text{Winslow}) \approx 2.48$. 
This suggests that the results of Table \ref{tab:refinement_study_bat_shape_p2to5}, which correspond to refinement levels that are of practical interest, are in the non-asymptotic convergence regime while the asymptotic convergence rate (which corresponds to a practically undesirably large number of DOFs) is slightly better than the table suggests. However, given that the choice $p=2$ produces a convergence log-plot comprised of nearly straight lines, the results of Table \ref{tab:refinement_study_bat_shape_p2to5} suggest that the convergence rate is slightly below the for $p=2$ maximally attainable rate of $\kappa(\, \cdot \,) = 2$ in this case, most likely as a result of the nonlinearity. Overall, the results suggest that the two NDF discretisations in mixed form perform similarly well while slightly outperforming the $C^0$-DG approach at the expense of higher computational costs. Of all the discretisations presented in Section \ref{sect:numerical_schemes}, the regularised weak form discretisation consistently produces the best results. In fact, it performs only marginally worse than Winslow's original approach, yet converging significantly faster and more reliably while furthermore avoiding local minima in practice. In our practical experience the $C^0$-DG approach, despite being outperformed by the discretisations in mixed form, suffices for the purpose of finding a nondegenerate or nearly nondegenerate initial iterate for initialising the regularised weak form discretisation in the vast majority of cases. As such, a combination of the two methods constitutes the best trade-off between robustness, solution quality and computational costs. \\
All NDF discretisations converge very reliably in typically $~5$ iterations when initialised with the forward Laplace initial guess from Section \ref{subsect:initial_guess}, making them suitable for the use in autonomously operating workflows. In autonomous workflows, combining an NDF discretisation with a posteriori refinement in case nondegeneracy does not carry over to the numerical approximation constitutes the most robust choice. Here, a mixed-form discretisation becomes a viable choice thanks to the better convergence rate. While the computational costs are higher, they remain manageable when operating on the Schur complement of the bilinear form's constant blocks. Overall, the Hessian recovery approach tends to be the better choice in this case, despite the problem's larger cardinality compared the the rotation-free approach, since it can be combined with Newton's method. Accelerating the convergence of the rotation-free approach with constitutes a topic for future research. For this, it may be possible to adopt a multipatch generalisation of the preconditioned Anderson acceleration approach from \cite{ji2023improved} which is highly effective in the singlepatch case. \\

\begin{figure}[h!]
\centering
\begin{subfigure}[b]{0.4 \textwidth}
    \includegraphics[width=\textwidth]{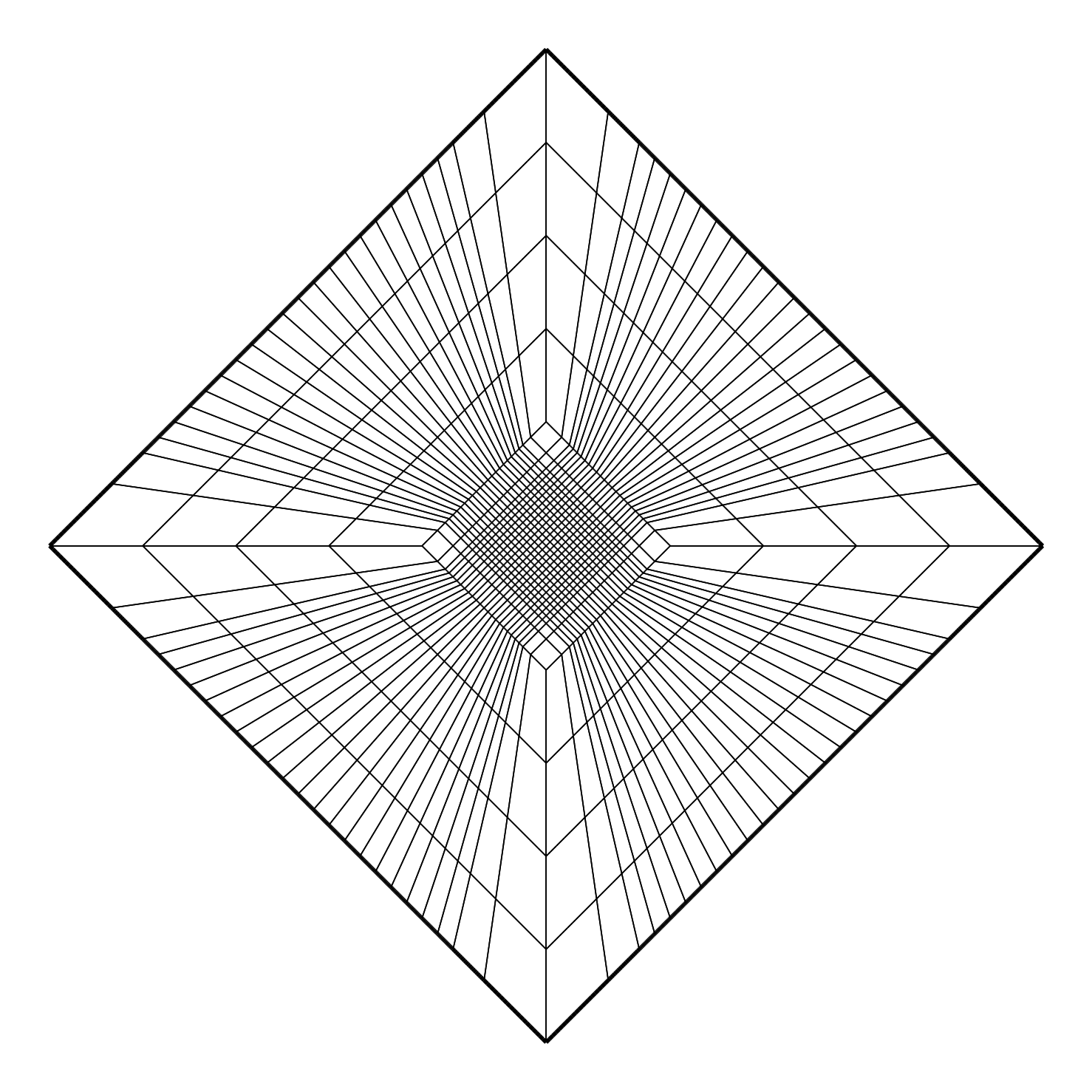}
    \caption{The parametric domain multipatch covering showing the knotspans at the zeroth refinement level.}
    \label{fig:male_screw_benchmark_param}
\end{subfigure}%
\hspace*{1cm}
\begin{subfigure}[b]{0.4 \textwidth}
    \includegraphics[width=\textwidth]{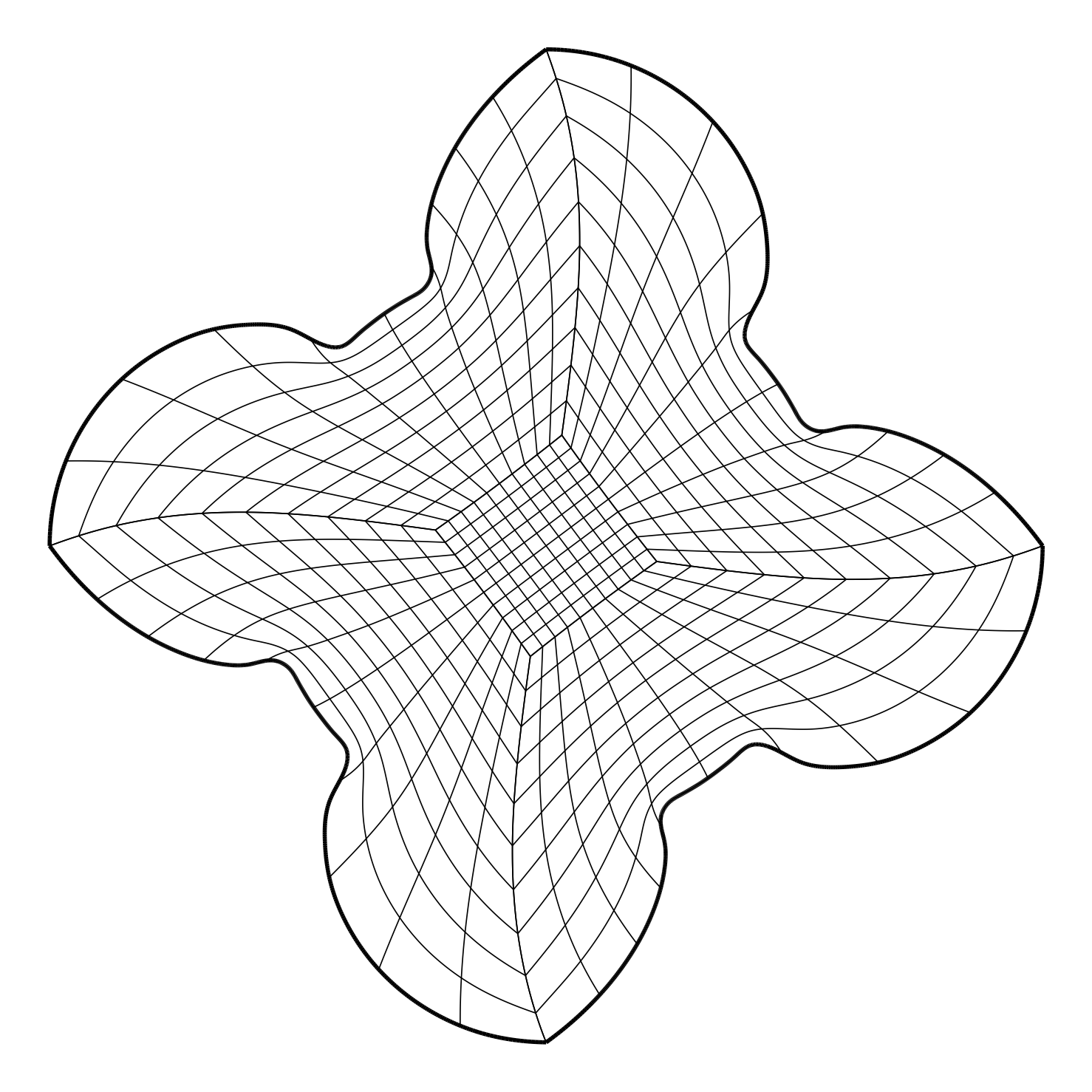}
    \caption{Uniform plot of the harmonic map at refinement level $0$ computed with the $C^0$-DG approach from Section \ref{subsubsect:blechschmidt}.}
    \label{fig:male_screw_benchmark_geom}
\end{subfigure}
\caption{The parametric domain and the approximation of a harmonic map for a screw geometry.}
\label{fig:male_screw_benchmark}
\end{figure}

\begin{figure}[h!]
\centering \includegraphics[width=0.5 \textwidth]{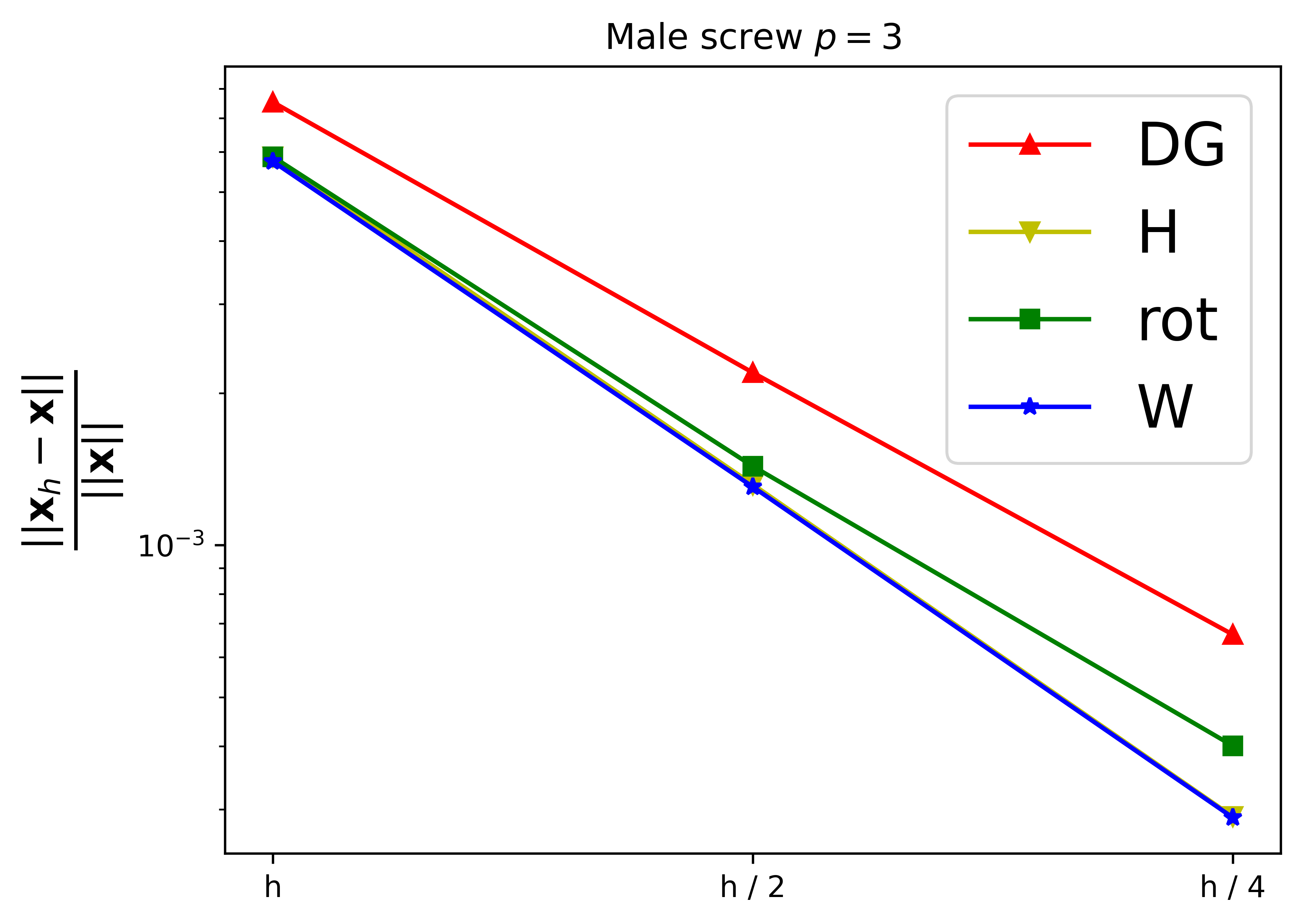}
\caption{Convergence plot of the various methods applied to the male screw geometry from Figure \ref{fig:male_screw_benchmark} with $p=3$.}
\label{fig:male_screw_convergence_plot}
\end{figure}

As a final example, we are considering the screw geometry from Figure \ref{fig:male_screw_benchmark_geom} which is mapped inversely harmonically into the parametetric domain from Figure \ref{fig:male_screw_benchmark_param}, which also shows the knotspans of the bicubic knotvectors with maximum regularity on each individual patch. Since the boundary correspondence is itself a piecewise bicubic spline with maximum regularity, we are only considering the choice $p=3$ here. The convex corners of $\hOm$ are mapped onto the convex corners of $\Om$ in the counterclockwise direction. Table \ref{tab:refinement_study_male_screw_p3} contains the approximate convergence rates for the various numerical schemes while Figure \ref{fig:male_screw_convergence_plot} shows the approximate $H^1$-norm distance to the exact solution, as before.

\begin{table}[h!]
\renewcommand{\arraystretch}{1.5}
\centering
\begin{tabular}{c|c|c|c|c|c}
 $p = 3$ & $\kappa\left(\mathcal{L}^H \right)$ & $\kappa \left( \mathcal{L}^{\text{rot}} \right)$ & $\kappa \left(\mathcal{L}^{\text{DG}}\right)$ & $\kappa \left(\mathcal{L}^{\text{W}}\right)$ & $\kappa \left(\text{Winslow} \right)$ \\
 \hline
  & 2.12 & 1.97 & 1.78 & 2.13 & 2.12
\end{tabular}
\caption{Approximate convergence rates for the male screw geometry using the various numerical schemes.}
\label{tab:refinement_study_male_screw_p3}
\end{table}

From the table and plot we may largely draw the same conclusions as in the previous example with the $C^0$-DG approach being outperformed by the other methods while the weak form discretisation fares the best, even slightly outperforming Winslow's method in this example. 

%% file: sec-control_mechanisms.tex
\section{Control Mechanisms}
\label{sect:control_mechanisms}

\subsection{Techniques for Parametric Control}
The parameterisations generated by Winslow's original approach~\eqref{eq:Winslow_pullback} or its PDE-based counterparts perform well on a wide range of benchmark problems \cite{gravesen2012planar}. However, as individual applications may require parameterisations with specific features, such as boundary layers in flow problems or homogeneous cell-sizes in problems subject to a CFL-condition, the techniques from Section \ref{sect:numerical_schemes} may be too rigid. Clearly, choosing the multipatch covering $\mathcal{Q}$ based on the application's specific needs may provide relief. However, in practice, this may prove too restrictive since the covering remains bilinear, which does not, for instance, allow for the creation of boundary layers. \\
Parametric control can be achieved in two main ways:
\begin{enumerate}
    \item Augmenting the standard inverse Laplace problem with a nonhomogeneous diffusivity.
    \item Mapping inversely harmonically into a parametric domain with a curvilinear rather than a Cartesian coordinate system.
\end{enumerate}
Let $\boldsymbol{\phi}: \Om_1 \rightarrow \mathbb{R}^2$ satisfy:
\begin{align}
\label{eq:classical_inverse_laplace_diffusivity}
    i \in \{1, 2\}: \quad \nabla \cdot \left( D \nabla \boldsymbol{\phi}_i \right) = 0 \quad \text{in } \Om_1, \quad \text{s.t. } \boldsymbol{\phi} = \mathbf{F} \text{ on } \partial \Om_1.
\end{align}
For point 1., we state the following theorem \cite{bauman2001univalent}:
\begin{theorem}[Divergence-form equations]
\label{thrm:nondegeneracy_div_form_equations}
   Let $D \in \text{SPD}^{2 \times 2}(\Om_1)$ be uniformly elliptic and let $\boldsymbol{\phi} \in H^1(\Om_1, \mathbb{R}^{2}) \cap C^0(\Om_1, \mathbb{R}^2)$ be the weak solution of~\eqref{eq:classical_inverse_laplace_diffusivity}. If $\mathbf{F}$ is diffeomorphic between $\partial \Om_1$ and $\partial \Om_2$ and $\Om_2$ is convex, then $\boldsymbol{\phi}: \Om_1 \rightarrow \Om_2$ satisfies $\det \partial_{\bx} \boldsymbol{\phi} \geq 0$ a.e. in $\Om_1$.
\end{theorem}
Under stronger regularity requirements on $D$, $\Om_1$ and $\Om_2$, Theorem \ref{thrm:nondegeneracy_div_form_equations} can be extended to uniform nondegeneracy $\det \partial_{\bx} \boldsymbol{\phi} > 0$ (a.e. in $\Om_1$). For details we refer to \cite{bauman2001univalent}. This means in particular that for merely essentially bounded $D \in \text{SPD}^{2 \times 2}(\Om)$, we need to account for the possibility of $\det \partial_{\bxi} \bx \rightarrow 0$ or $\det \partial_{\bxi} \bx \rightarrow \infty$ in the interior of $\hOm$, which may require stabilisation. Taking $\Om_1 = \Om$ and $\Om_2 = \hOm$, it is reasonable to assume that Theorem \ref{thrm:nondegeneracy_div_form_equations} also applies to, for instance, the weak-form approach from Section \ref{subsect:weak_form_discretisation}, even though it exchanges the dependencies, i.e.,  $\bxi(\bx) \rightarrow \bx(\bxi)$. A limitation of introducing a nonhomogeneous diffusivity is that it is currently unknown whether the inverted problem can be cast into a form that does not contain the Jacobian determinant in the denominator, as in the NDF-discretisation from Section \ref{subsect:NDF_discretisations}. However, the NDF-discretisations remain highly practical since they can compute a nondegenerate reference solution to initialise an iterative scheme with $D \neq \mathcal{I}^{2 \times 2}$ based on the weak-form discretisation. \\

\noindent For point 2., a coordinate transformation is conveniently accomplished by introducing a controlmap $\br: \hOm \rightarrow \hOmr$. As such, we now allow the target domain of $\bx^{-1}: \Om \rightarrow \hOmr$ to be a parametric surface, too. In what follows, differential operators receive a subscript to indicate differentiation w.r.t. various coordinate systems. For instance, $\nabla \rightarrow \nabla_{\br}$ to indicate differentiation w.r.t. the entries of $\br: \hOm \rightarrow \hOmr$. \\
The introduction of $\br: \hOm \rightarrow \hOmr$ furthermore enables creating boundary correspondences $\mathbf{F}^{\br \rightarrow \bx}: \partial \hOmr \rightarrow \partial \Om$ that are diffeomorphic between $\hOmr$ and $\Om$ when $\Om$ has no corners by, for instance, choosing $\hOmr$ to be the unit disc. For Theorem \ref{thrm:RKC} to apply to the pair $(\hOmr, \Om)$, we require $\hOmr$ to be convex. We denote the map that maps inversely harmonically into the domain $\Om^{\br} = \br(\hOm)$ by $\bx^{\br}(\br): \hOmr \rightarrow \Om$. The same map can be converted to the original coordinate system via a pullback. We employ the abuse of notation $\bx^{\br}(\br)$ instead of $\bx(\br(\bxi))$ to indicate a change of coordinate system and assume that the reader is aware of the compositions involved. \\
\pichskip{15pt}
\parpic[r][b]{%
  \begin{minipage}{100mm}
  \centering
    \includegraphics[width=1 \textwidth]{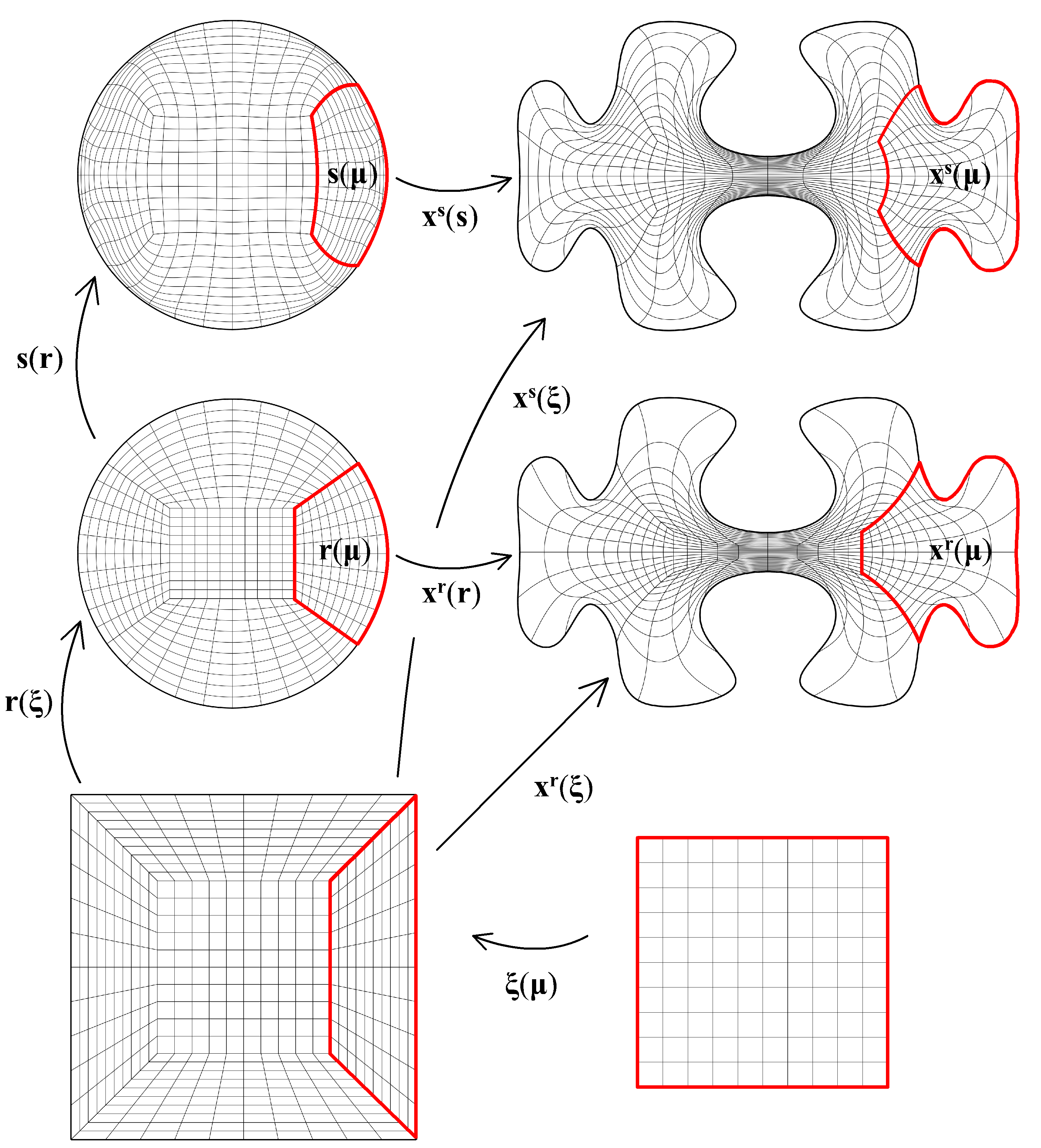}%
    \captionof{figure}{Figure summarising the dependencies between $\bx, \bs, \br, \bxi$ and $\bmu$.}
    \label{fig:dependency_summary}
  \end{minipage}
}


\noindent In the IGA-setting, parametric control via $\br: \hOm \rightarrow \hOmr$ is conveniently achieved by reinterpreting the PDE-based formulations over the Cartesian coordinate system $\bxi = (\xi_1, \xi_2)^T$ as problems posed over $\hat{\Om}^{\br}$, with the curvilinear coordinate system induced by $\br(\bxi) = (r_1(\xi_1, \xi_2), r_2(\xi_1, \xi_2))^T$. We may then use basic differential geometry identities to express the associated integrals in the original coordinate system via a pullback. As such, the operators from Section \ref{sect:numerical_schemes} now receive an additional $\br$-dependence. However, this does not change the nature of the equations as long as the $\br \circ \bm^i$ remain diffeomorphisms. \\

As in point 1., a broad class of reparameterisation methods follows from seeking the controlmap as the solution of~\eqref{eq:classical_inverse_laplace_diffusivity}. More precisely, given a reference controlmap $\br: \hOm \rightarrow \hOmr$, we may find a new controlmap $\bs(\br): \hOmr \rightarrow \hOmr$ by solving~\eqref{eq:classical_inverse_laplace_diffusivity} with $\Om_1 = \Om_2 = \hOmr$ while selecting a diffusivity $D$ that builds desired features into the solution. The map $\bs(\bxi): \hOm \rightarrow \Om^{\br}$ then follows from a pullback and a diffeomorphic boundary correspondence $\mathbf{F}^{\br \rightarrow \bs}: \partial \hOmr \rightarrow \hOmr$ is given by the identity.
\noindent If $\bs(\br): \hOmr \rightarrow \hOmr$ is the identity on $\partial \hOmr$, the effect of the coordinate transformation, induced by $\bs(\bxi): \hOm \rightarrow \hOmr$, on $\bx^{\bs}(\bxi): \hOm \rightarrow \Om$ can be predicted by noting that $\bx^{\bs}(\bxi) = \bx^{\br} \circ \bs(\bxi)$. All the dependencies between $\bx, \bs, \br, \bxi$ and the local coordinate system $\bmu$ in $\Omega^{\square}$ are summarised in Figure \ref{fig:dependency_summary}. \\

\noindent As before, depending on the regularity of the diffusivity, the solutions may contain singularities $\det J \rightarrow 0$ or unbounded growth, i.e., $\det J \rightarrow \infty$. While the vanishing or diverging of $\det J$ is typically avoided by discrete approximations, this behaviour will be observable in a refinement study. For the $\bs(\bxi) \circ \bm^i$ to be diffeomorphisms, this means that for $\bs \in H^1(\hOmr, \mathbb{R}^{2})$, jumps in the Jacobian $\partial_{\br} \bs(\br)$ may only occur on the $\br(\partial \hOm_i)$. As such, we require the diffusivity to be patchwise continuous. \\

\noindent Given a reference controlmap $\br: \hOm \rightarrow \hOmr$, the most general approach combines methods $1.$ and $2.$, leading to the coupled system
\begin{align}
\label{eq:domain_optimisation_coupled_equation}
   \text{find } (\bx(\br), \bs(\br)) \quad \text{ s.t. for } i \in \{1, 2 \}: \left \{ \begin{array}{c} \nabla_{\bx} \cdot \left(D^{\bx} \nabla_{\bx} \bs_i \right) = 0 \\ \nabla_{\br} \cdot \left(D^{\bs} \nabla_{\br} \bs_i \right) = 0 \end{array} \right. \text{ in } \hOmr \quad \text{ and } \quad \left \{ \begin{array}{l} \bx = \mathbf{F}^{\br \rightarrow \bx} \\ \bs = \mathbf{F}^{\br \rightarrow \bs} \end{array} \right. \text{ on } \partial \hOmr,
\end{align}
where, typically, $\mathbf{F}^{\br \rightarrow \bs}(\br) = \br$. Here, the practically useful dependencies are $D^{\bx} = D^{\bx}(\bs, \bx)$ and $D^{\bs} = D^{\bs}(\br, \bx)$. Note that the first equation is inverted and the unknown becomes the differential operator $\nabla_{\bx}$ as a function of $\bx(\br): \hOmr \rightarrow \Om$. The system is associated with a global operator comprised of two separate operators, one for each equation
\begin{align}
\label{eq:domain_optimization_coupled}
    \mathcal{L}^{\bx, \bs}(\bx, \bs, \boldsymbol{\phi}_1, \boldsymbol{\phi}_2, D^{\bx}, D^{\bs}) = \mathcal{L}^{\bx}(\bx, \boldsymbol{\phi}_1, D^{\bx}, \bs) + \mathcal{L}^{\bs}(\bs, \boldsymbol{\phi}_2, D^{\bs}, \bx).
\end{align}
For $D^{\bx} = \mathcal{I}^{2 \times 2}$, the operator $\mathcal{L}^{\bx}(\cdot, \, \cdot, \, \cdot, \, \cdot)$ can be based on any of the operators from Section \ref{sect:numerical_schemes}. For $D^{\mathbf{x}} \neq \mathcal{I}^{2 \times2}$ we restrict ourselves to the weak-form discretisation. With $\bs = \bs(\br)$, the regularised weak-form operator becomes (c.f. equation~\eqref{eq:inverse_elliptic_weak_pullback_operator}):
\begin{align}
\label{eq:weak_form_operator_controlmap}
    \mathcal{L}^{\text{W}}_{\varepsilon}(\bx, \boldsymbol{\phi}, D^{\bx}, \bs) = \int \limits_{\hOmr}\frac{\left(C(\partial_{\br} \bs) \nabla_{\br} \boldsymbol{\phi} \right) \, \colon \left(Q^T(\partial_{\br} \bx, \partial_{\br} \bs) D^{\bx}(\bs, \bx) \, Q(\partial_{\br} \bx, \partial_{\br} \bs) \right)}{\mathcal{R}_{\varepsilon} \left( \det Q(\partial_{\br} \bx, \partial_{\br} \bs) \right)} \, \mathrm{d} \br,
\end{align}
with
$$
Q(\partial_{\br} \bx, \partial_{\br} \bs) := C(\partial_{\br} \bx) \left(\nabla_{\br} \bs \right) \quad \text{and } C(\, \cdot \,) \text{ as in~\eqref{eq:A_CC}.}
$$
The operator corresponding to the second part $\nabla_{\br} \cdot \left(D^{\bs} \nabla_{\br} \bs \right) = 0$ reads
\begin{align}
    \mathcal{L}^{\bs}(\bs, \boldsymbol{\phi}, D^{\bs}, \bx) = \int \limits_{\hOmr} \partial_{\br} \boldsymbol{\phi} \, \colon \, ( \partial_{\br} \bs \, D^{\bs}(\br, \bx)) \, \mathrm{d} \br.
\end{align}
We note that the coordinate transformation $\partial_{\bxi} \rightarrow \partial_{\bs(\bxi)}$ in the NDF operators from Section \ref{subsect:NDF_discretisations} reintroduces the Jacobian determinant $\det \partial_{\bxi} \bs$ in the denominator upon pullback of the equations from $\hOmr$ into $\hOm$. As such, an iterative algorithm has to be initialised with a nondegenerate controlmap $\bs^0: \hOmr \rightarrow \hOmr$ when $\bx^{\bs}: \hOm \rightarrow \Om$ and $\bs: \hOmr \rightarrow \hOmr$ are coupled via $D^{\bx}$ or $D^{\bs}$. In this case, we recommend basing the scheme on~\eqref{eq:weak_form_operator_controlmap} instead. In the classical literature, parametric control via $\bs$, rather than through a pullback, is accomplished by introducing additional terms in~\eqref{eq:EGG_classical} \cite[Chapter~4]{thompson1998handbook}. While this formulation enables removing $\partial_{\bxi} \bs$ from the denominator, it is not applicable if $\bs \notin H^2(\hOm, \mathbb{R}^2)$ since it requires second-order derivative information of $\bs(\bxi)$, making it unsuited for this paper's use-cases. While it may be possible to reduce the regularity requirements of $\bs: \hOmr \rightarrow \hOmr$ in a way anologous to Section \ref{subsect:NDF_discretisations}, this is beyond the scope of this paper.\\

Given a reference controlmap $\br: \hOm \rightarrow \hOmr$, the controlmap $\bs: \hOmr \rightarrow \hOmr$ is conveniently built from a push-forward of the same finite-dimensional space $\mathcal{V}_h$ used to represent the map $\bx_h: \hOm \rightarrow \Om$. As mentioned before, substituting a degenerate intermediate controlmap $\bs: \hOmr \rightarrow \hOmr$, produced by an iterative root-finding algorithm applied to the coupled system, may cause problems due to division by zero. In practice this is avoided by initialising the scheme by the tuple $(\bx^{\br}, \br)$ (i.e., the solution for $D^{\bx} = D^{\bs} = \mathcal{I}^{2 \times 2}$ over the reference controlmap $\br: \hOm \rightarrow \hOmr$) which is computed using one of the NDF-discretisation from Section \ref{subsect:NDF_discretisations}. The barrier term in~\eqref{eq:weak_form_operator_controlmap} then prevents intermediate iterates $(\bx^{\bs}, \bs)^i$ from leaving the set of nondegenerate maps. As before, a discretisation takes the test functions $(\boldsymbol{\phi}_1, \boldsymbol{\phi}_2)$ from the finite-dimensional space $\mathcal{U}_h^{\mathbf{0}} \times \mathcal{U}_h^{\mathbf{0}}$ and finds the root using Newton's method. In practice, the coupled scheme converges reliably for a wide range of diffusivities $D^\bx, D^{\bs}$ when initialised with $(\bx^{\br}, \br)$. \\

\noindent Depending on the choice of $D^{\bx}$ and $D^{\bs}$, the solution of the coupled system may no longer be uniformly nondegenerate (even for boundary correspondences that lead to UNDG maps for $D = \mathcal{I}^{2\times2}$). To avoid singularities, we shall often introduce a stabilisation on the patch vertices. As such, let $\Gamma^v = \{ \mathbf{v}^1, \ldots, \mathbf{v}^{N_v} \} \subset \overline{\hOm}$ be the set of patch vertices shared by at least two patches, i.e.,
$$\Gamma^v := \left \{ \mathbf{v} \in \overline{\hOm} \, \, \big| \, \, \exists (i, j) \in \{1, \ldots, N_p\} \times \{1, \ldots, N_p\} \text{ s.t. } \overline{\hOm}_i \cap \overline{\hOm}_j = \{ \mathbf{v} \} \right \}.$$
The singularities are avoided by introducing an appropriate regularisation. For this purpose, we introduce the Gaussian blending functions
\begin{align}
\label{eq:exponential_decay_regularisation}
    g_i^{\kappa}(\br) := A_i \exp{\left(-\left(\frac{\kappa}{d_i^{\text{min}}} \| \br - \br(\mathbf{v}^i) \| \right)^2 \right)}, \quad \text{with} \quad d_i^{\min} := \min \left \{ \|\mathbf{v}^i - \mathbf{v}^j \| \, \, \big| \, \, j \in \{1, \ldots, N_p\} \setminus \{i\} \right \}
\end{align}
and $\kappa > 0$. Here, the $A_i$ are chosen such that
$$\forall \mathbf{v}^i \in \Gamma^v: \quad \sum_{i=1}^{N_v} g_i^{\kappa}(\br(\mathbf{v}^i)) = 1.$$
Let $D$ be the diffusivity in question and let $\mathcal{D}_i = \{D_i^1, \ldots, D_i^q \}$ be the set containing the limits
$$D_i^j := \lim \limits_{\bxi \rightarrow \mathbf{v}^i} D(\br(\bxi)) \quad \text{ s.t. } \bxi \in \hOm_j \quad \text{for each patch } \hOm_j \text{ with } \overline{\hOm}_j \cap \{\mathbf{v}^i\} = \{ \mathbf{v}^i \}.$$
We define $\overline{D}_i$ as the average of the $D_i^j \in \mathcal{D}_i$, i.e.,
\begin{align}
\label{eq:heuristical_D_scaling}
    \overline{D}_i := \sum \limits_{D_i^j \in \mathcal{D}_i} D_i^j.
\end{align}
The regularisation ensures that the regularised $\overline{D}^{\kappa}(D) \in \text{SPD}^{2 \times 2}(\hOmr)$ is single-valued in the $\br(\mathbf{v}^i)$ by replacing
\begin{align}
\label{eq:patch_vertex_stabilisation}
    D \rightarrow \overline{D}^{\kappa}(D) := \left( 1 - \sum \limits_{i=1}^{N_v} g_i^{\kappa} \right) D + \sum_{i=1}^{N_v} g_i^{\kappa} \overline{D}_i.
\end{align}
The decay rate $\kappa > 0$ in~\eqref{eq:exponential_decay_regularisation} tunes the degree of regularisation and is relatively insensitive to the characteristic length-scale of $\Omega^{\br}$ thanks to the scaling by $d_i^{\min}$. It should be noted that other regularisations exist besides~\eqref{eq:heuristical_D_scaling}.

To better see what the effect of reparameterising under a nonhomogenous diffusivity $D$ is, we note that~\eqref{eq:classical_inverse_laplace_diffusivity} is the Euler-Lagrange equation of
\begin{align}
\label{eq:laplace_reparam_minimisation}
    \min \limits_{\boldsymbol{\phi}: \Om_1 \rightarrow \Om_2} \frac{1}{2} \int \limits_{\Om_1} \operatorname{tr} \left(\nabla_{\bx} \boldsymbol{\phi}^T D \nabla_{\bx} \boldsymbol{\phi} \right) \mathrm{d} \bx 
    \quad \text{s.t.} \quad \boldsymbol{\phi} = \mathbf{F} \text{ on } \partial \Om_1.
\end{align}
The interpretation as a minimisation problem is helpful in predicting the effect of $D$ on $\boldsymbol{\phi}: \Om_1 \rightarrow \Om_2$.

\subsection{Patch interface removal}
\label{subsect:interface_removal}
\begin{figure}[h!]
\centering
    \includegraphics[width=0.7 \textwidth, align=c]{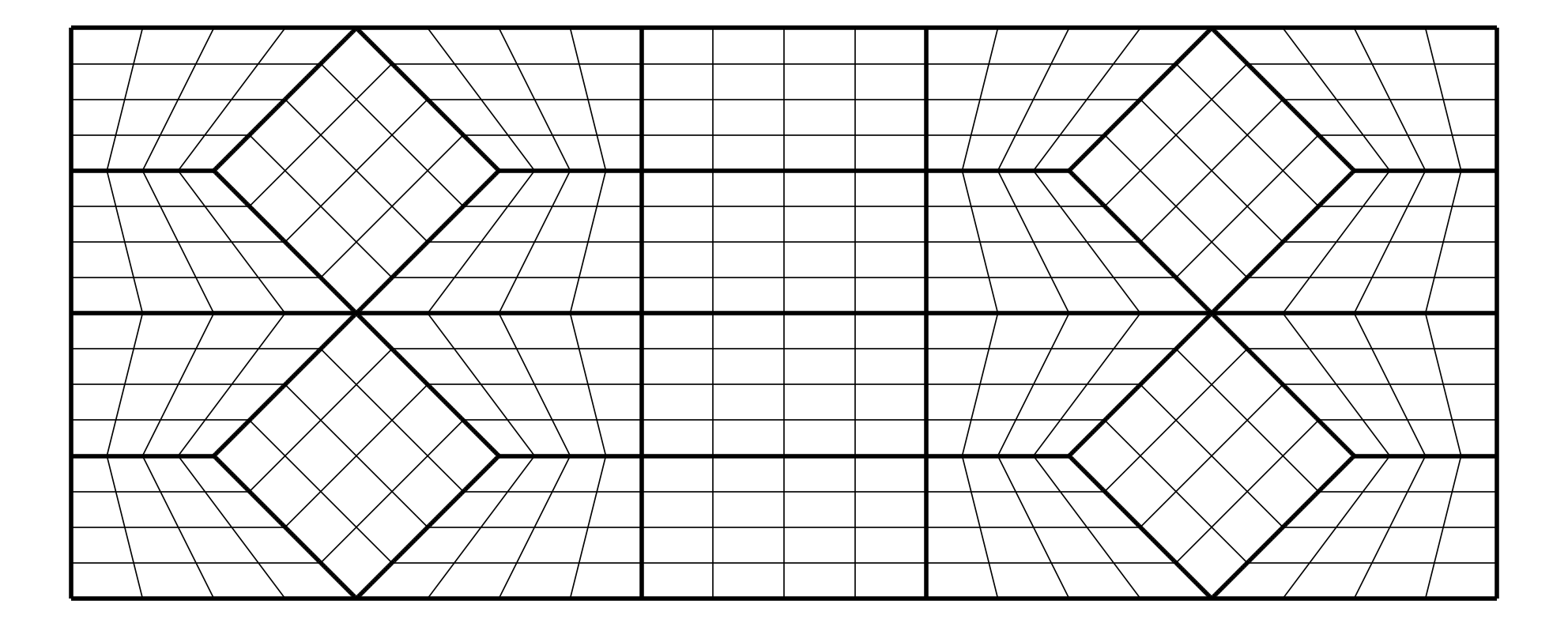}
\caption{Bilinear multipatch covering of a rectangular domain comprised of $24$ patches.}
\label{fig:interface_removal_regularised_ref}
\end{figure}

\begin{figure}[h!]
\centering
    \includegraphics[height=3.0cm, align=c]{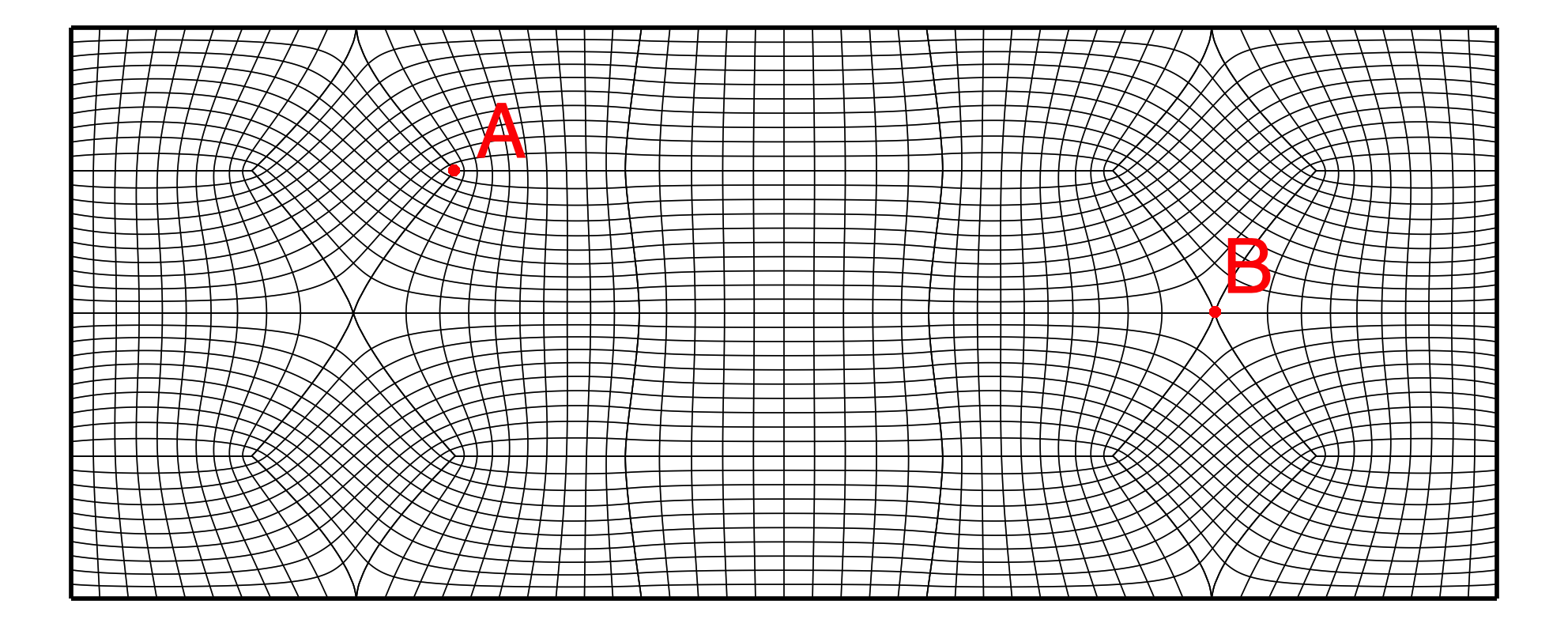} \includegraphics[height=3.0cm, align=c]{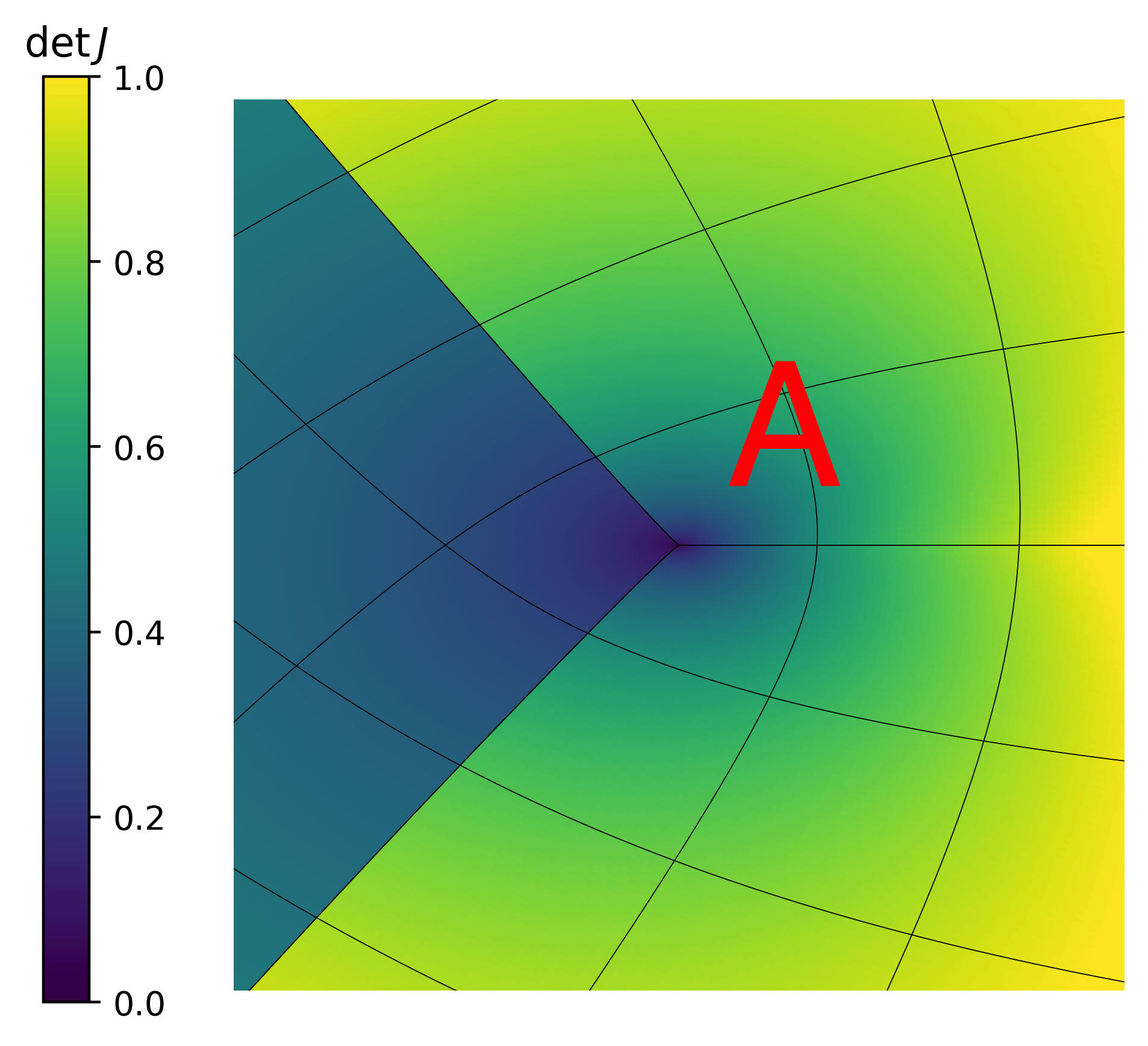} \includegraphics[height=3.0cm, align=c]{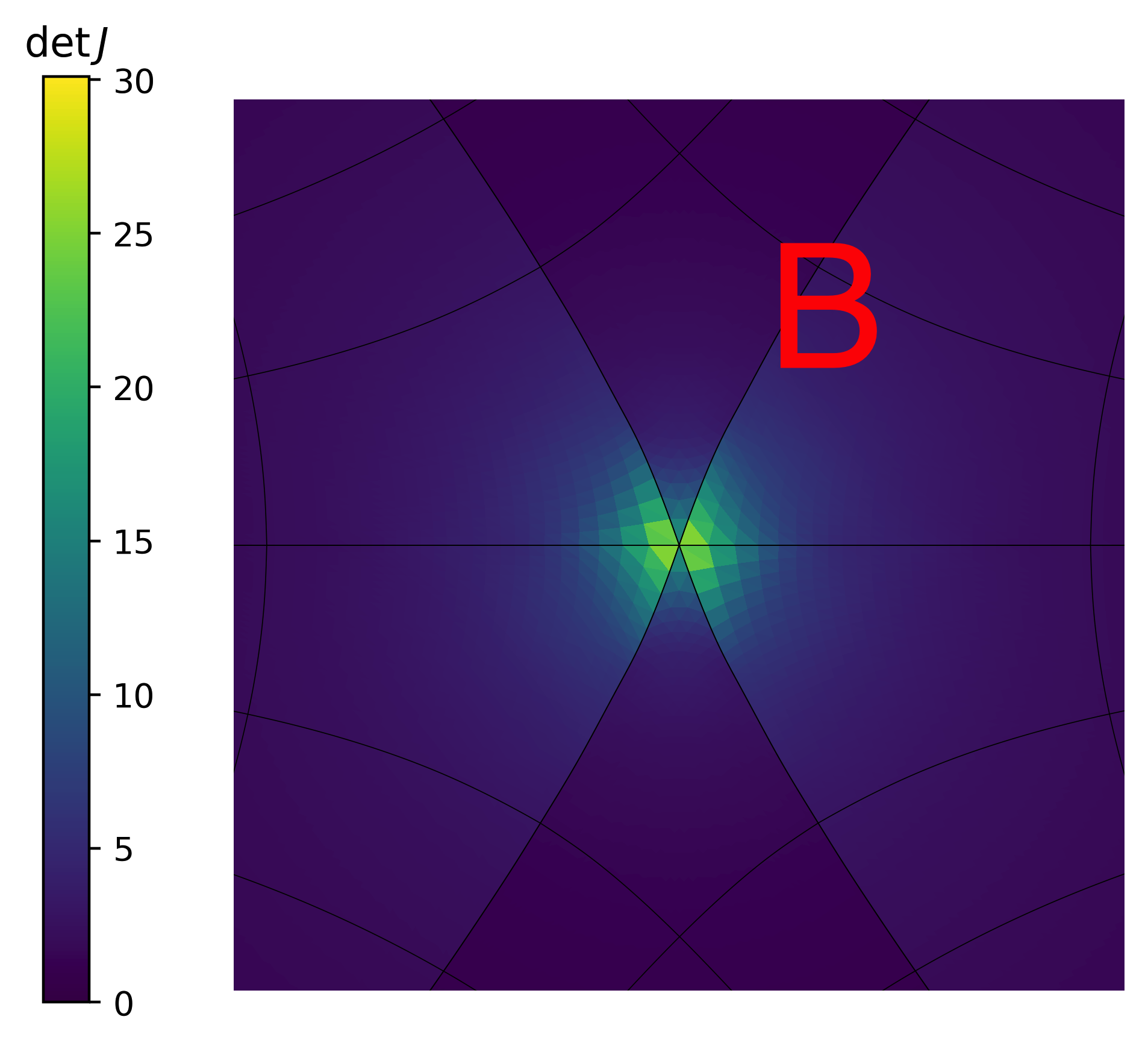} \\
    \includegraphics[height=3.0cm, align=c]{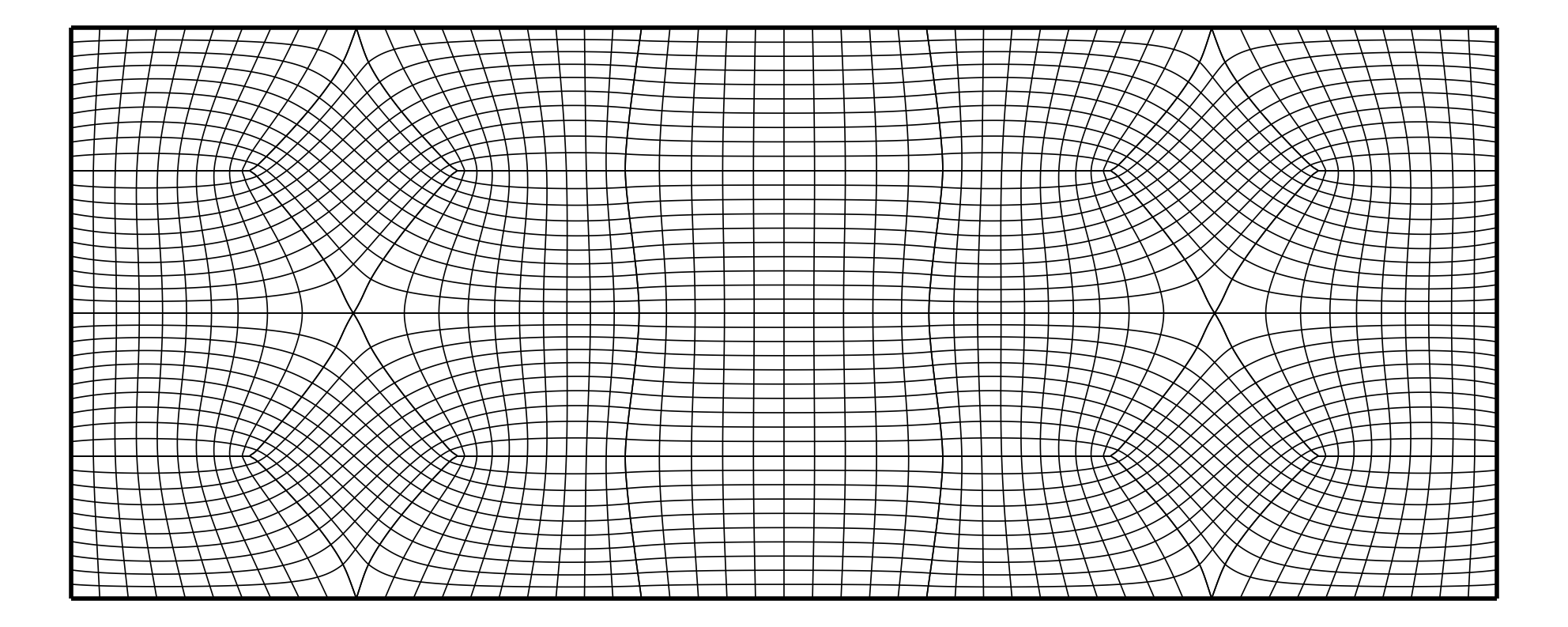}
    \includegraphics[height=3.0cm, align=c]{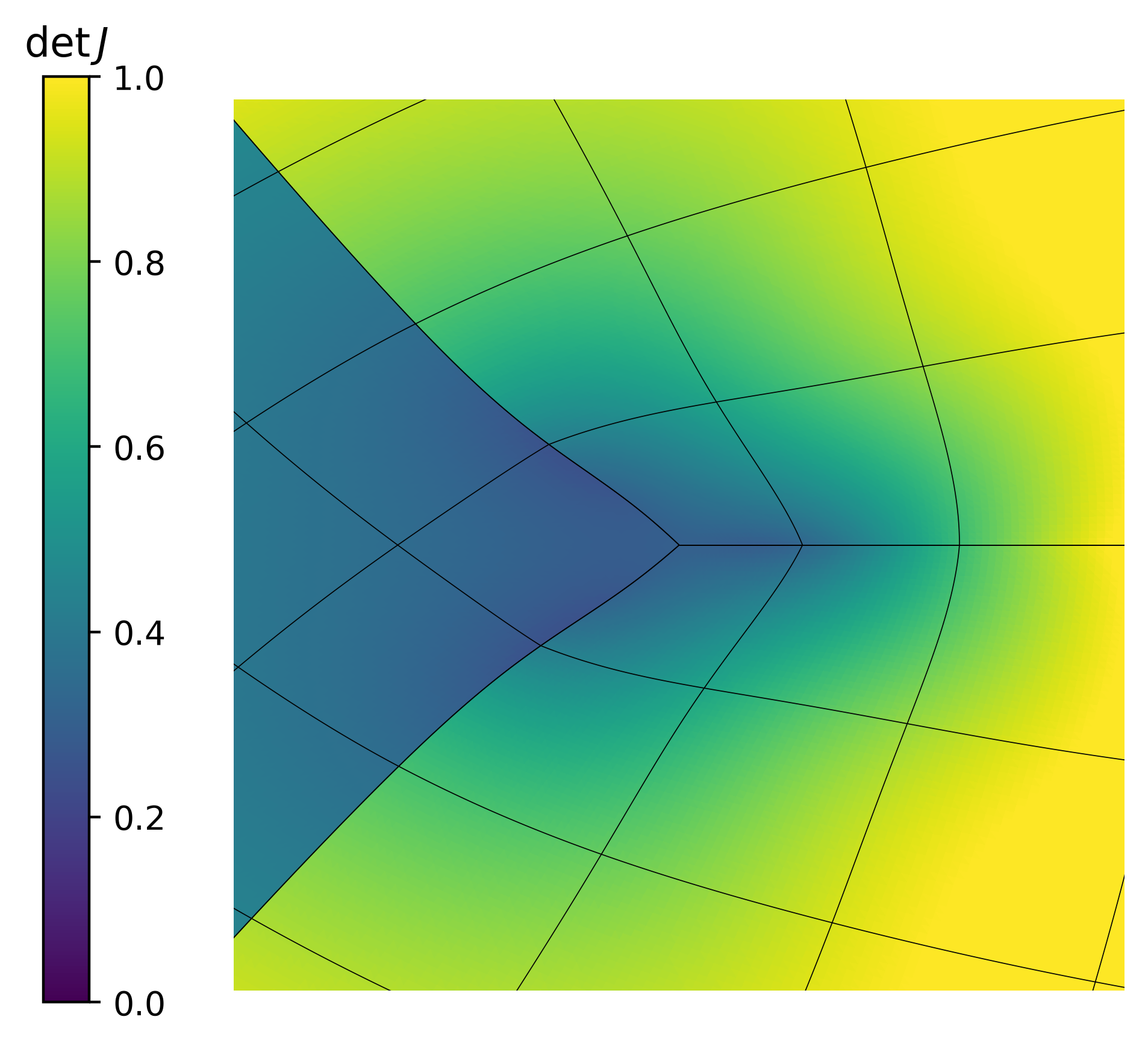} \includegraphics[height=3.0cm, align=c]{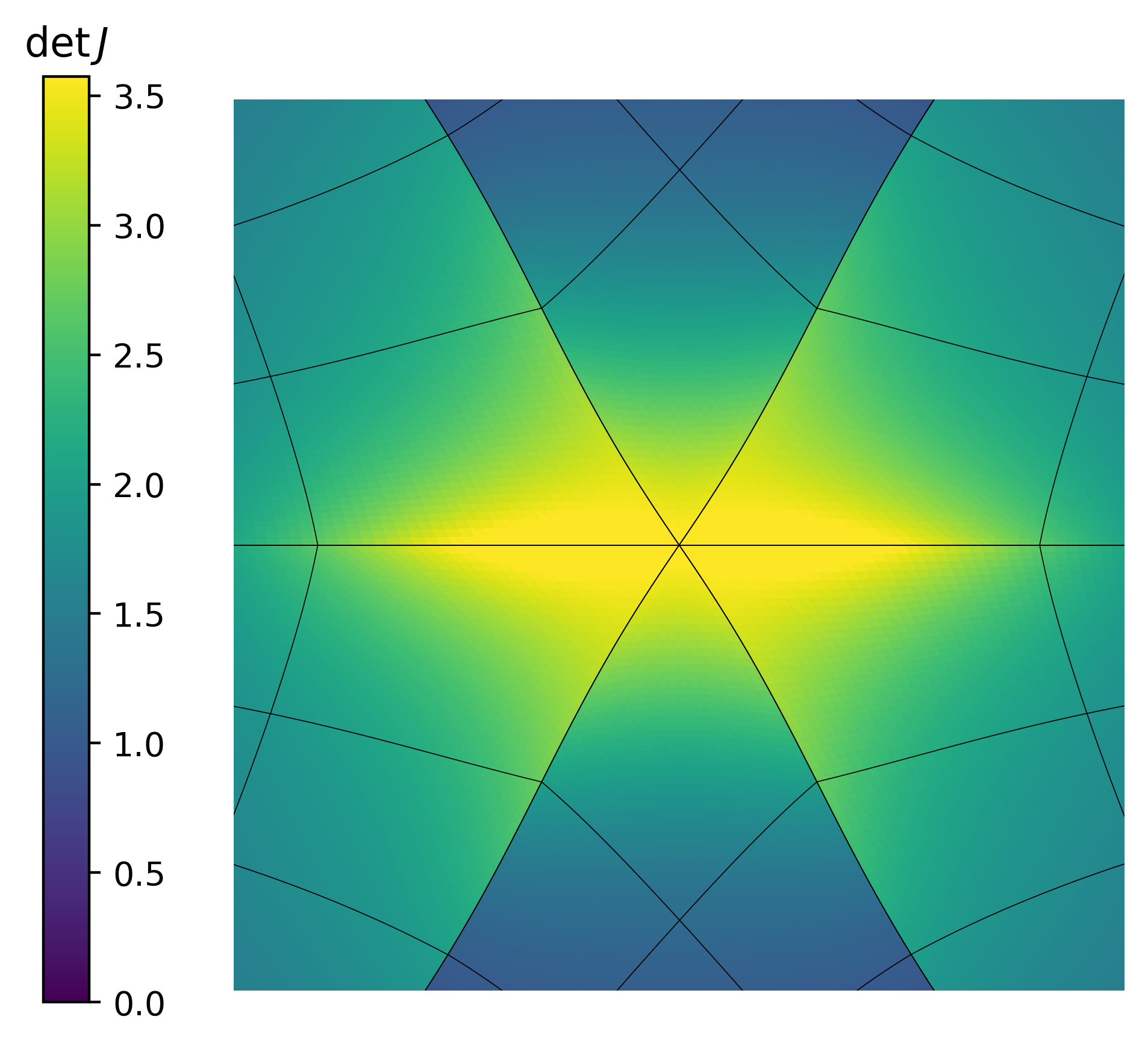}
\caption{Plot of the controlmaps $\bs: \hOm \rightarrow \hOm$ (top) and $\bs^{\text{reg}}: \hOm \rightarrow \hOm$ (bottom) along with a zoom-in onto vertices in which the Jacobian determinant grows / shrinks unboundedly in the absence of regularisation.}
\label{fig:interface_removal_regularised}
\end{figure}

The images of local (in $\Om^{\square}$) isolines under the $\bx^{\br}(\bxi) \circ \bm^i$ will generally form a (possibly steep) angle across patch interfaces when joined together on $\Om$. In certain applications, it may be desirable to decrease or largely remove the steep interface angles. As $\bx: \hOmr \rightarrow \Om$ is diffeomorphic in $\hOmr$, a controlmap $\bs: \hOmr \rightarrow \hOmr$ that removes steep angles will, by extension, remove them in the recomputed map $\bx^{\bs}: \hOm \rightarrow \Om$. As such, interface removal can be regarded as an \textit{a priori} step since it requires no prior knowledge of $\bx: \hOm \rightarrow \Om$. \\
We would like to accomplish
\begin{align}
\label{eq:local_jump_0}
    \forall \gamma_{jk} \in \Gamma^I: \quad  [\![ (\partial_{\bmu^{\perp}} \bs(\br))]\!] = 0  \quad \text{on } \br(\gamma_{jk}),
\end{align}
wherein $[\![ \, \cdot \, ]\!]$ now denotes the ordinary entry-wise jump term while $\partial_{\bmu^{\perp}}$ denotes the directional derivative transversal to $\br(\gamma_{ij})$ (i.e., either $\partial_{\mu_1}$ or $\partial_{\mu_2}$ on $\br(\gamma_{jk}^+)$ and $\br(\gamma_{jk}^-)$). Requirement~\eqref{eq:local_jump_0} can be weakly enforced by utilising in~\eqref{eq:domain_optimisation_coupled_equation} the diffusivity
\begin{align}
\label{eq:laplace_local_trace_penalty}
    D^{\bs}(\br, \bx) = D^{\bs}_{\Gamma^I}(\br) := \partial_{\mu_1} \br \otimes \partial_{\mu_1} \br + \partial_{\mu_2} \br \otimes \partial_{\mu_2} \br \quad \text{on} \quad \br(\hat{\Omega}_i).
\end{align}
Meanwhile, if $D^{\bx} = \mathcal{I}^{2 \times 2}$,~\eqref{eq:domain_optimisation_coupled_equation} is decoupled and the map $\bx(\bxi): \hOm \rightarrow \Om$ can be computed from a degenerate initial guess using an NDF-discretisation. The diffusivity from~\eqref{eq:laplace_local_trace_penalty} urges $\bs: \hOmr \rightarrow \hOmr$ to map the patchwise isolines $\bs(\bmu)$ smoothly across patch interfaces.
As the magnitude of $\partial_{\mu_i} \br$ depends on $\br \circ \bm^i: \Om^{\square} \rightarrow \br(\hOm_i)$ (and may therefore be subject to considerable changes between patches), the best results are obtained by a normalisation, i.e.,
\begin{align}
\label{eq:laplace_local_normalised_trace_penalty}
    D^{\bs}(\br, \bx) = \widehat{D}^{\bs}_{\Gamma^I}(\br) := \widehat{\partial}_{\mu_1} \br \otimes \widehat{\partial}_{\mu_1} \br + \widehat{\partial}_{\mu_2} \br \otimes \widehat{\partial}_{\mu_2} \br \quad \text{on} \quad \br(\hat{\Omega}_i), \quad \text{with } \widehat{\partial}_{\mu_i} \br := \frac{1}{\left \| \partial_{\boldsymbol{\mu_i}} \br \right \|} \partial_{\boldsymbol{\mu_i}} \br.
\end{align}
Note that 
$$\operatorname{tr}(\mathcal{I}^{2 \times 2}) = \operatorname{tr}\left(\widehat{D}^{\bs}_{\Gamma^I}(\br) \right) = 2.$$
The normalisation has a similar effect as minimising~\eqref{eq:local_jump_0} while now suppressing jumps in the normalised transversal component of the $\br(\gamma_{jk})$ with $\gamma_{jk} \in \Gamma^I$, i.e., we are penalising jumps in the transverse direction but not the direction's magnitudes.
\begin{figure}[h!]
\centering
    \includegraphics[width=0.4 \textwidth]{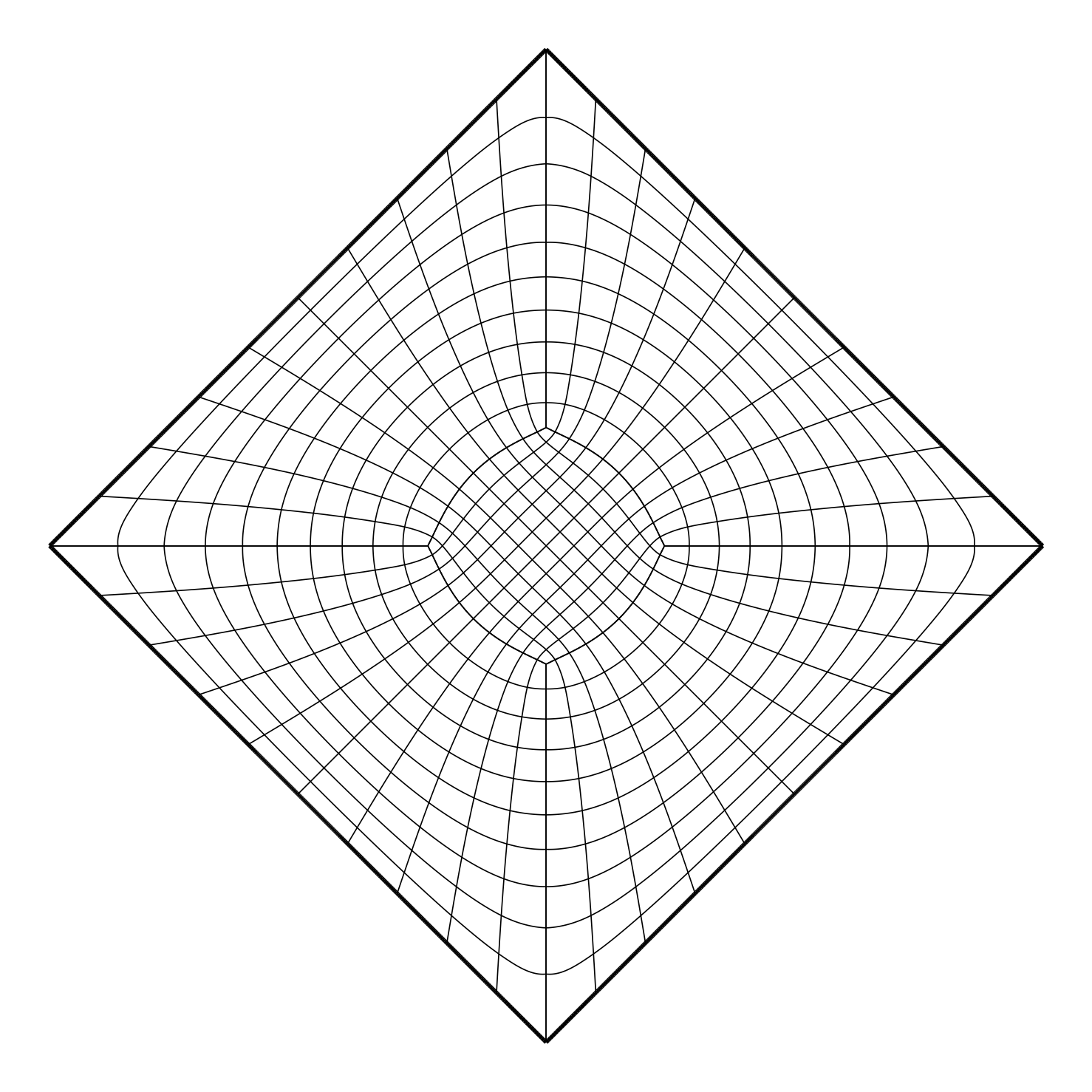} \hspace*{0.5cm}
    \includegraphics[width=0.4 \textwidth]{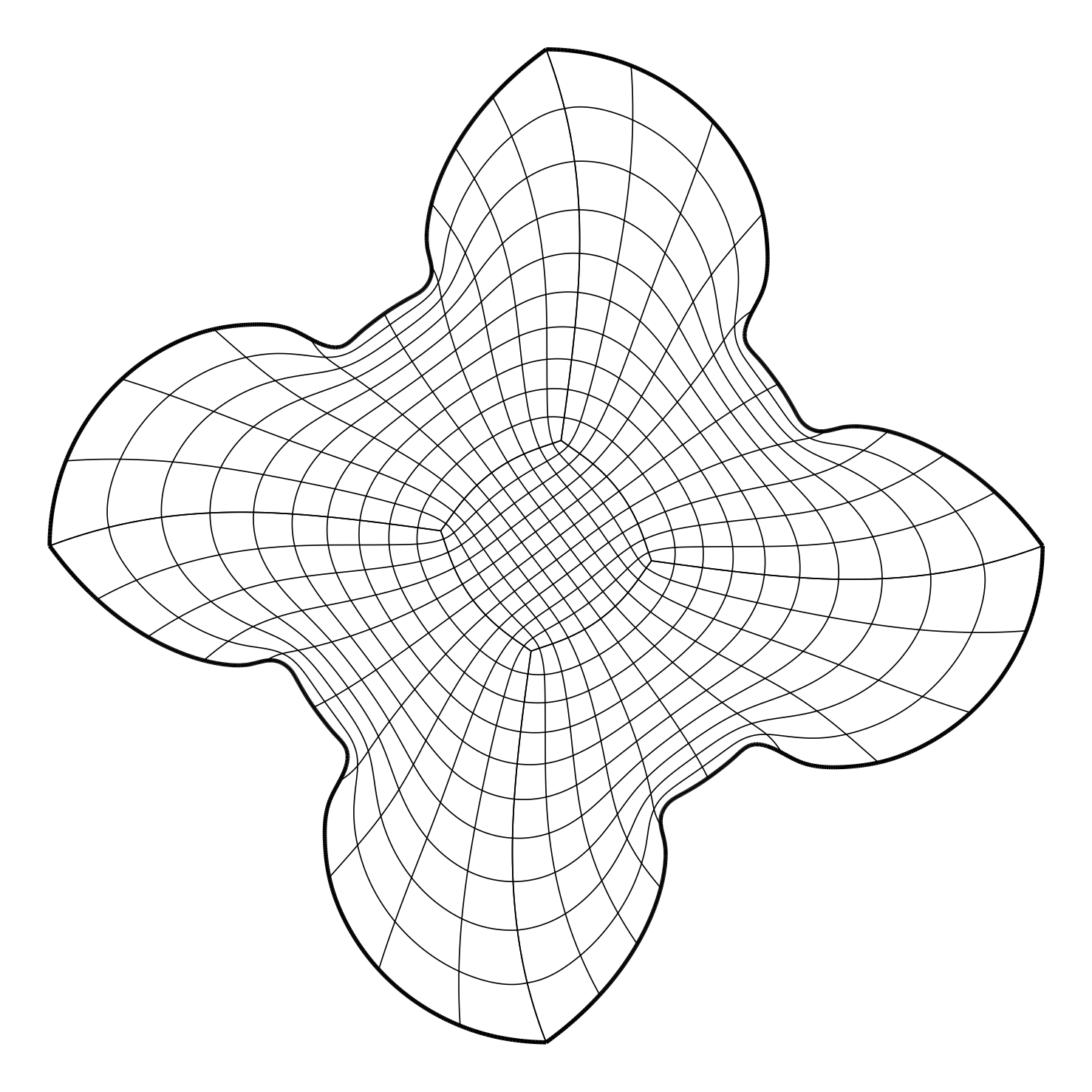}
\caption{The screw geometry after suppression of local transverse gradient jumps.}
\label{fig:male_screw_interface_remove_k0}
\end{figure}
\noindent We are considering the screw geometry from Figure \ref{fig:male_screw_benchmark_geom} and perform interface removal with normalisation. Figure \ref{fig:male_screw_interface_remove_k0} shows the resulting reparameterisation. As a measure of the degree of interface removal, we utilise the following value
\begin{align}
\label{eq:interface_removal_norm}
    L_{\Gamma}^2(\bx) = \sum \limits_{\gamma_{jk} \in \Gamma^I} \int \limits_{\br(\gamma_{jk})} \left\| [\![ \widehat{\partial}_{\bmu^{\perp}} \bx(\br)]\!] \right \|^2 \mathrm{d} \Gamma,
\end{align}
where $\widehat{\partial}_{\bmu^{\perp}}$ denotes the normalised directional derivative transverse to $\br(\gamma_{jk})$. With
$$\frac{L_{\Gamma}(\bx_h^{\bs})}{L_{\Gamma}(\bx_h)} \approx 0.0998$$ 
the technique is highly effective. \\ 
As stated in Theorem \ref{thrm:nondegeneracy_div_form_equations}, methods based on~\eqref{eq:domain_optimisation_coupled_equation} do not exclude singularities in $\bs: \hOmr \rightarrow \hOm^{\bs}$. While singularities are in practice avoided by discrete approximations, for merely essentially bounded $D^{\bs} \in \text{SPD}^{2 \times 2}(\hOmr)$, we may expect 
$$\inf \limits_{\hOmr} \det \partial_{\br} \bs \rightarrow 0 \quad \text{or} \quad \sup \limits_{\hOmr} \det \partial_{\br} \bs \rightarrow \infty$$
in a refinement study. For patchwise continuous $D^{\bs} \in \text{SPD}^{2 \times 2}(\hOmr)$, singularities (and unbounded gradients), if present, are located in the $\br(\mathbf{v}^i)$ with $\mathbf{v}^i \in \Gamma^v$. The creation of singularities can be avoided by employing the stabilisation from~\eqref{eq:patch_vertex_stabilisation}. We are considering the rectangular parametric domain $\hOm = \hOmr$ along with an irregular multipatch covering depicted in Figure \ref{fig:interface_removal_regularised_ref}.
\noindent We perform a refinement study of (normalised) interface removal with and without regularisation, initially assigning a uniform cubic knotvector with three internal knots to each side $L_i \in \Gamma^B$ and facet $\gamma_{ij} \in \Gamma^I$. Each refinement $h \rightarrow h/2$ halves the knotvector's knotspans. We are monitoring the value of $\nu_{\Gamma}( \, \cdot \,) := L_{\Gamma}(\, \cdot \,) / L_{\Gamma}(\br)$ for the original $\bs: \hOmr \rightarrow \hOmr$ and its regularised counterpart $\bs^{\text{reg}}: \hOmr \rightarrow \hOmr$, as well as the values $\min \det J_{\mathbf{v}}(\, \cdot \,)$ and $\max \det J_{\mathbf{v}}(\, \cdot \,)$ which are the minimum and maximum values of $\det \partial_{\br}( \, \cdot \, )$ over all patch vertices $\br(\mathbf{v}^i)$. As $\det \partial_{\br}(\, \cdot \,)$ is not single-valued in the $\mathbf{v}^i \in \Gamma^v$, we define this value as the minimum / maximum of taking the limit on each adjacent patch. Table \ref{tab:refinement_study_regularised_int_removal_mu_inf} contains the reference values in the absence of regularisation while Table \ref{tab:refinement_study_regularised_int_removal_mu_9} contains the corresponding values for the regularisation $D^{\bs} \rightarrow \overline{D}^{\kappa}(D^{\bs})$ with $\kappa = 9$. Furthermore, Figure \ref{fig:interface_removal_regularised} shows the controlmap with and without regularisation after the last refinement level. Table \ref{tab:refinement_study_regularised_int_removal_mu_inf} clearly demonstrates that $\min \det J_{\mathbf{v}}(\, \cdot \,)$ and $\max \det J_{\mathbf{v}}(\, \cdot \,)$ shrink / grow unboundedly in the absence of regularisation while Table \ref{tab:refinement_study_regularised_int_removal_mu_9} demonstrates that regularisation prevents further shrinkage / growth under refinement. 

\begin{table}[h!]
\renewcommand{\arraystretch}{1.5}
\centering
\begin{tabular}{c|c|c|c|c|c}
  & h & h/2 & h/4 & h/8 & h/16 \\
 \hline
$\nu_{\Gamma}(\bs)$ & 0.273 & 0.197 & 0.162 & 0.140 & 0.0991 \\
\hline
$\min \det J_{\mathbf{v}}(\bs)$ & 0.0679 & 0.0466 & 0.0357 & 0.0294 & 0.0184 \\
\hline
$\max \det J_{\mathbf{v}}(\bs)$ & 7.46 & 11.9 & 15.6 & 18.9 & 30.1
\end{tabular}
\caption{Reference values of performing interface removal on the quadrangulation from Figure \ref{fig:interface_removal_regularised_ref} in the absence of regularisation.}
\label{tab:refinement_study_regularised_int_removal_mu_inf}
\end{table}

\begin{table}[h!]
\renewcommand{\arraystretch}{1.5}
\centering
\begin{tabular}{c|c|c|c|c|c}
 $\kappa=9$ & h & h/2 & h/4 & h/8 & h/16 \\
 \hline
$\nu_{\Gamma}(\bs^{\text{reg}})$ & 0.426 & 0.413 & 0.407 & 0.406 & 0.406 \\
\hline
$\min \det J_{\mathbf{v}}(\bs^{\text{reg}})$ & 0.223 & 0.262 & 0.251 & 0.249 & 0.248 \\
\hline
$\max \det J_{\mathbf{v}}(\bs^{\text{reg}})$ & 4.44 & 3.93 & 3.84 & 3.86 & 3.84
\end{tabular}
\caption{The outcomes of performing regularised interface removal with $\kappa = 9$.}
\label{tab:refinement_study_regularised_int_removal_mu_9}
\end{table}

\noindent As expected, the regularisation also prevents monotone decrease of $\nu_{\Gamma}(\bs^{\text{reg}})$, eventually settling for a value of $\sim 0.4 L_{\Gamma}(\br)$ relative to the reference value of Figure \ref{fig:interface_removal_regularised_ref}. Meanwhile, the corresponding value shrinks unboundedly in the absence of regularisation. \\
For larger values of $\kappa$, we expect the discretisation to settle for a lower value of $\nu_{\Gamma}(\bs^{\text{reg}})$ at the expense of reducing / increasing the values of $\min \det J_{\mathbf{v}}(\bs^{\text{reg}})$ and $\max \det J_{\mathbf{v}}(\bs^{\text{reg}})$. Table \ref{tab:refinement_study_regularised_int_removal_mu_18} contains the outcomes for regularisation with $\kappa = 18$.

\begin{table}[h!]
\renewcommand{\arraystretch}{1.5}
\centering
\begin{tabular}{c|c|c|c|c|c}
 $\kappa=18$ & h & h/2 & h/4 & h/8 & h/16 \\
 \hline
$\nu_{\Gamma}(\bs^{\text{reg}})$ & 0.315 & 0.314 & 0.309 & 0.303 & 0.299 \\
\hline
$\min \det J_{\mathbf{v}}(\bs^{\text{reg}})$ & 0.0892 & 0.143 & 0.179 & 0.174 & 0.167 \\
\hline
$\max \det J_{\mathbf{v}}(\bs^{\text{reg}})$ & 6.36 & 6.98 & 6.33 & 6.16 & 6.08
\end{tabular}
\caption{The outcomes of performing regularised interface removal with $\kappa = 18$.}
\label{tab:refinement_study_regularised_int_removal_mu_18}
\end{table}
\noindent Indeed, the table confirms this expectation, settling for a value of $\sim 0.3 L_{\Gamma}(\br)$ while roughly doubling the shrinkage / growth of $\det \partial_{\br}(\bs^{\text{reg}})$ compared to $\kappa = 9$. \\
We conclude that the proposed regularisation is an effective means to tune the degree of interface removal at the expense of cell size homogeneity reduction. In practice, an appropriate choice of the decay rate $\kappa > 0$ is furthermore relatively insensitive to the average distance between the $\mathbf{v}^i \in \Gamma^v$, thanks to the scaling by $d_i^{\min}$ in~\eqref{eq:exponential_decay_regularisation}.

\subsection{Cell size homogenisation}
\label{subsubsect:cell_size_hom}
\begin{figure}[h!]
\centering
\begin{subfigure}[b]{0.8 \textwidth}
    \includegraphics[width=0.45 \textwidth]{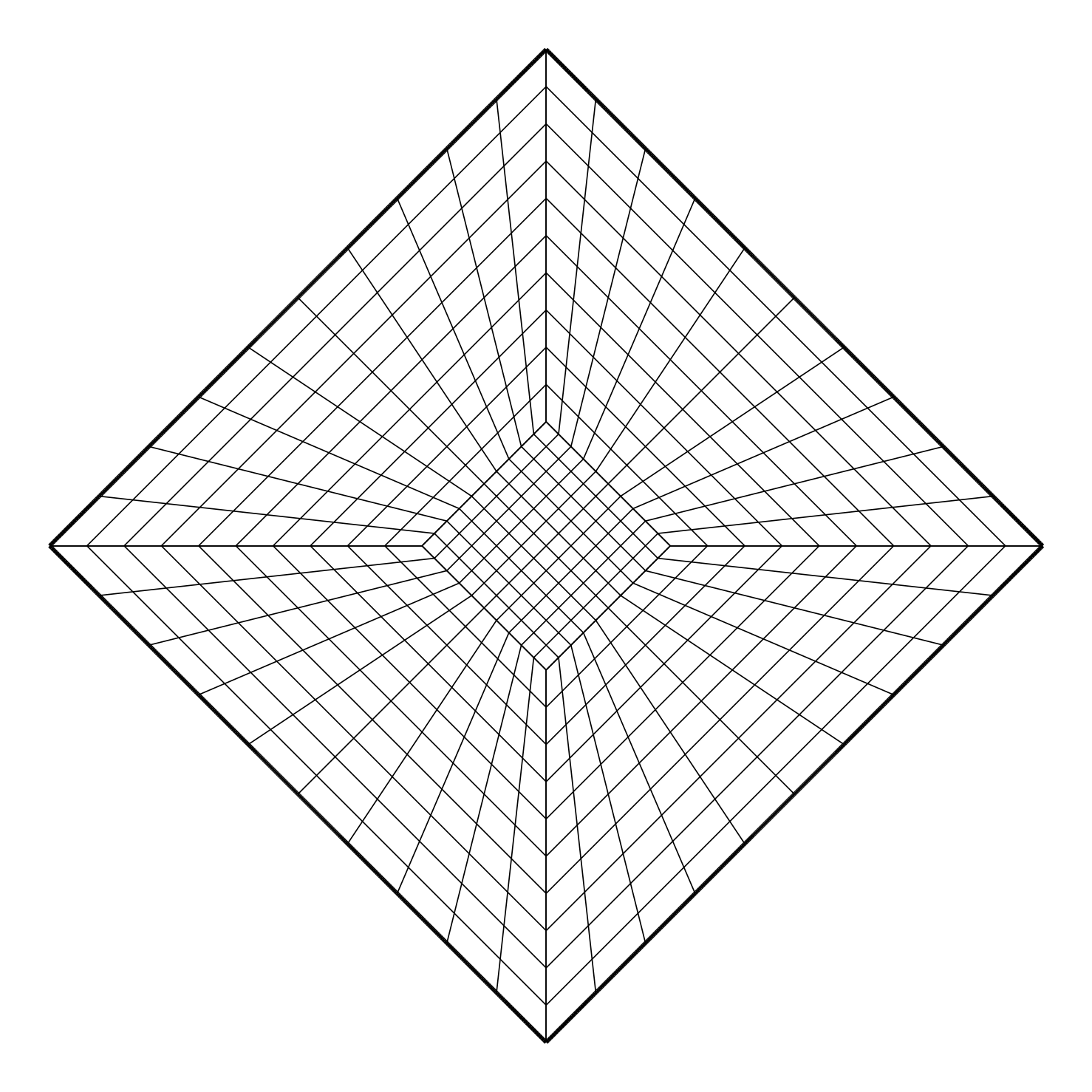} \hfill
    \includegraphics[width=0.45 \textwidth]{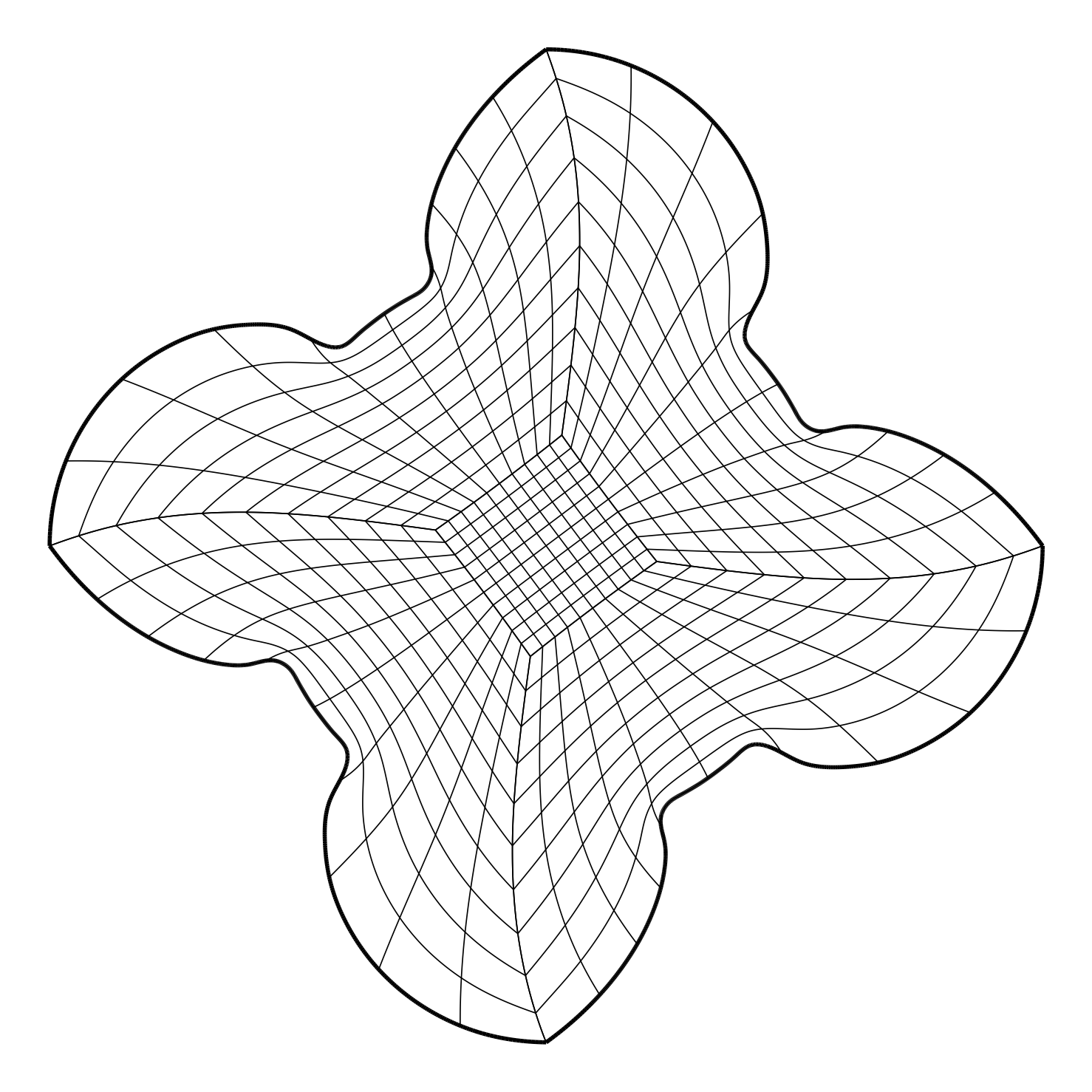}
    \caption{The reference geometry and controlmap.}
    \label{fig:male_screw_cell_k0}
\end{subfigure}
\\
\begin{subfigure}[b]{0.8 \textwidth}
    \includegraphics[width=0.45 \textwidth]{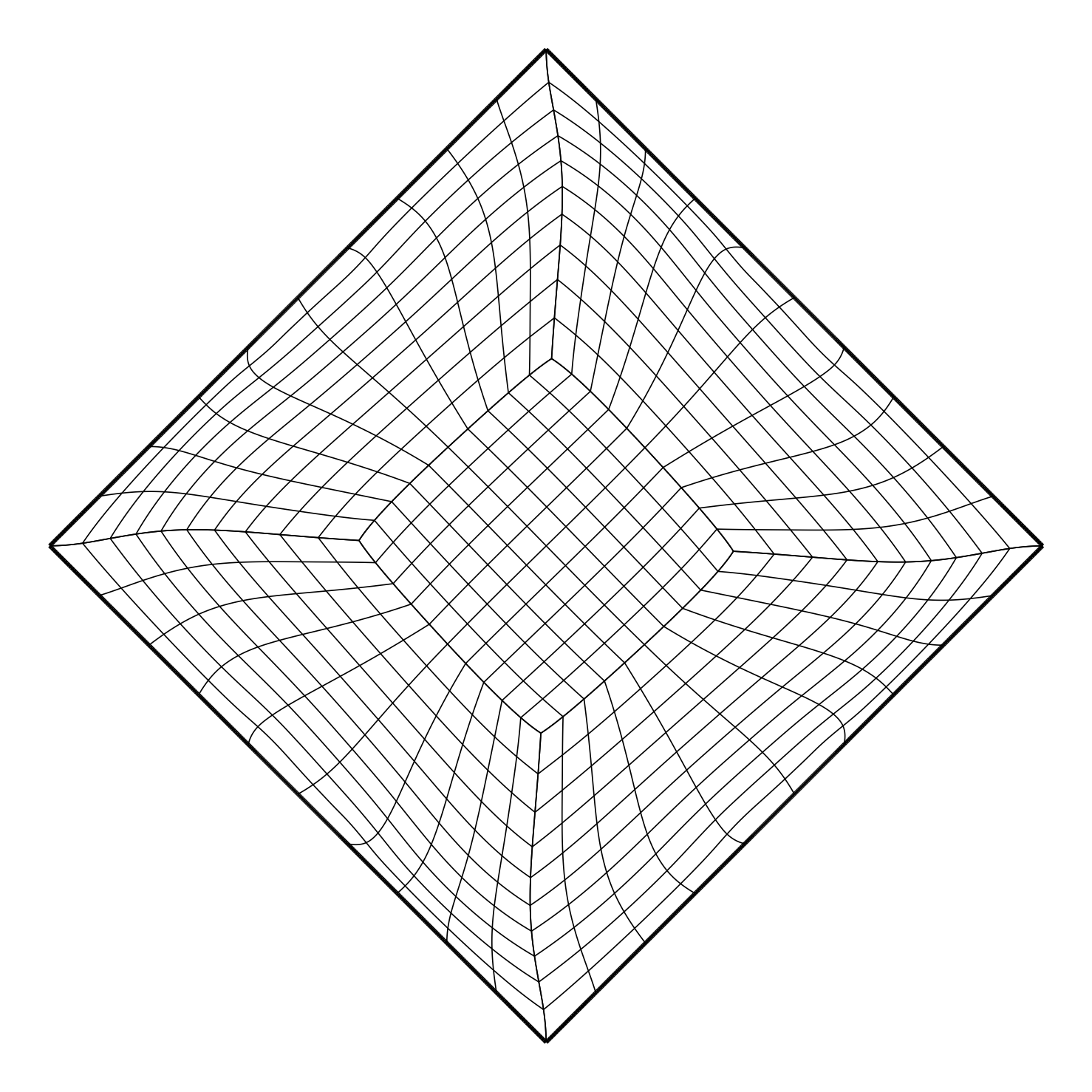} \hfill
    \includegraphics[width=0.45 \textwidth]{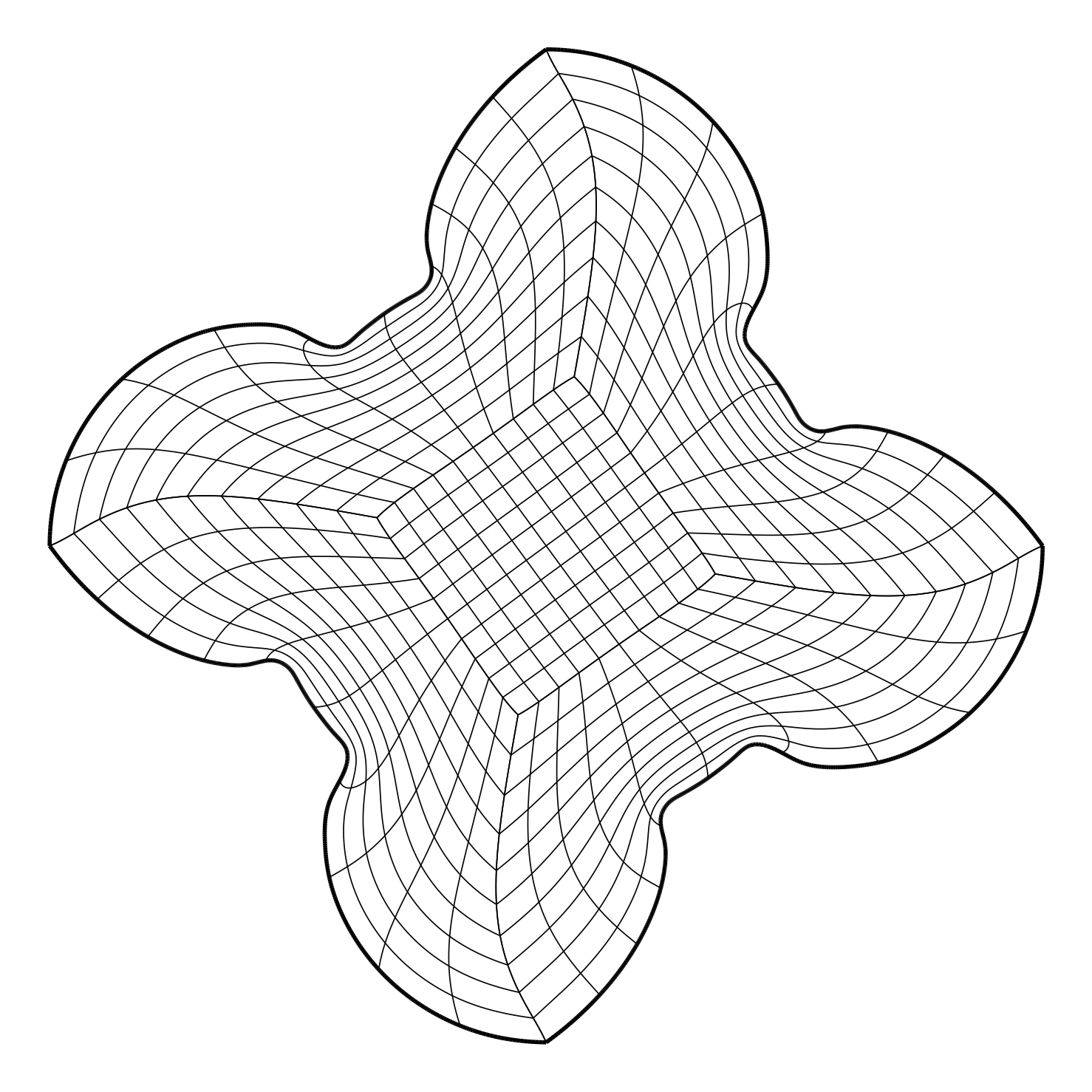}
    \caption{The screw geometry along with a plot of the corresponding controlmap $\bs: \hOm \rightarrow \hOm$ after reparameterisation with $k=2$.}
    \label{fig:male_screw_cell_k1}
\end{subfigure}
\caption{The screw geometry before and after reparameterisation under~\eqref{eq:cell_size_hom_diffusivity} with $k=2$.}
\label{fig:male_screw_cell_hom}
\end{figure}

A popular measure for the parameterisation's cell size homogeneity is the \textit{Area} functional
\begin{align}
    L_{\text{Area}}(\bx) = \int \limits_{\hOm} \left(\det J(\bx) \right)^2 \mathrm{d} \bxi,
\end{align}
which measures the variance of $\det J(\bx)$ over $\hOm$, with smaller values indicating better homogeneity. In the multipatch setting, it is more natural to measure the homogeneity on each individual patch and summing over all patches
\begin{align}
\label{eq:Area_multipatch}
    L_{\text{Area}}(\bx) := \sum_{i=1}^{N_p} \int \limits_{\Om^{\square}} \left( \det \partial_{\bmu} (\bx(\bxi) \circ \bm^i) (\bmu) \right)^2 \mathrm{d} \bmu,
\end{align}
wherein $\bmu = (\mu_1, \mu_2)^T$ denotes the free coordinate functions in $\Om^{\square}$. Direct minimisation of~\eqref{eq:Area_multipatch} over the map's inner controlpoints leads to a nonconvex problem which is furthermore prone to yielding degenerate maps. \\
In the context of the coupled system~\eqref{eq:domain_optimisation_coupled_equation}, there are two main ways to achieve homogenisation without having to resort to nonconvex optimisation:
\begin{enumerate}
    \item Designing a diffusivity $D^{\bs}(\br, \bx^{\bs})$ that contracts / expands the cell sizes of $\bs: \hOm^{\br} \rightarrow \hOm^{\br}$ wherever the cell sizes of $\bx^{\bs}: \hOm^{\br} \rightarrow \Om$ are large / small, while taking $D^{\bx}(\bs, \bx^{\bs}) = \mathcal{I}^{2 \times 2}$.
    \item Picking $D^{\bs} = D^{\bs}(\br)$ (i.e., $D^{\bs}$ has no dependency on $\bx^{\bs}$) while designing a diffusivity $D^{\bx} = D^{\bx}(\bx^{\bs})$ that encourages cell size homogenisation.
\end{enumerate}
As for method 1., we notice that for $D^{\bx} = \mathcal{I}^{2 \times 2}$, the solution of the inverse Laplace problem is merely a property of the shapes $\hOmr$ and $\Om$ as well as the diffeomorphic boundary correspondence $\mathbf{F}^{\br \rightarrow \bx}: \partial \hOmr \rightarrow \partial \Om$. As such, a controlmap $\bs: \hOmr \rightarrow \hOmr$ that is the identity on $\partial \hOmr$ computes the composition $\bx^{\bs}(\br) = \bx^{\br} \circ \bs(\br)$. Therefore, we may require $\bs(\br)$ to contract cells in $\br(\hOm_i)$ wherever $\det \partial_{\bmu} \bx^{\bs}(\br)$ is large and vice versa. We may also choose to penalise based on $\partial_{\bxi} \bx^{\bs}(\br)$ or $\partial_{\br} \bx^{\bs}(\br)$ to reduce the Area functional in a different coordinate system if desired. \\
To cast this problem into the form of~\eqref{eq:domain_optimization_coupled}, we contract the cells of $\bs(\br)$ by penalising the value of $\operatorname{tr}(G^{\br \rightarrow \bs})$, where $G^{\br \rightarrow \bs}$ denotes the metric between the coordinate systems induced by $\br: \hOm \rightarrow \hOm^{\br}$ and $\bs: \hOm^{\br} \rightarrow \hOm^{\br}$. This is accomplished by introducing $D^{\bs}(\br, \bx^{\bs}) = \sigma(\bx^{\bs}) \mathcal{I}^{2 \times 2}$, where $\sigma(\bx^{\bs})$ assumes large values in regions where contraction is desired and vice-versa. For instance
$$\sigma^k(\bx^{\bs})(\br) = (\det \partial_{\bmu} \bx^{\bs}(\br))^k \quad \text{on} \quad \br(\hOm_i),$$
where larger values of $k > 0$ lead to a more drastic homogenisation. As such, we are solving the coupled system~\eqref{eq:domain_optimization_coupled} with $D^{\bx} = \mathcal{I}^{2 \times 2}$ and 
\begin{align}
\label{eq:cell_size_hom_diffusivity}
    D^{\bs}(\bx^{\bs}) = \sigma^k(\bx^{\bs}) \mathcal{I}^{2 \times 2}.
\end{align}
The contraction under $\sigma^k(\bx)$ has a similar effect as operating on $(\det \partial_{\br} \bs)^2$ directly while being inherently less prone to yielding degenerate discrete maps. \\
While a formal proof is lacking, it is plausible to assume that the coupled system is well-posed under this choice of $D^{\bs}$ since for any bijective $\bs: \hOm^{\br} \rightarrow \hOm^{\br}$ (and given a diffeomorphic boundary correspondence $\mathbf{F}^{\br \rightarrow \bx}$), the coupled system~\eqref{eq:domain_optimization_coupled} approximates a UNDG map $\bx^{\bs}: \hOm^{\br} \rightarrow \Om$ that satisfies $\sigma^k(\bx^{\bs}) > 0$ (a.e.) such that $D^{\bs}(\bx^{\bs}) \in \text{SPD}^{2 \times 2}(\hOm^{\br})$. However, a root-finding algorithm may diverge in case the Newton increment accidentally causes $\bx^{\bs}$ to leave the set of UNDG maps. In practice, this is avoided by the barrier property of~\eqref{eq:weak_form_operator_controlmap} and the scheme converges reliably using Newton's method with line search for a wide range of choices $k > 0$ when the scheme is initialised with the tuple $(\bx^{\br}, \br)$ (i.e., the reference solution and reference controlmap). \\
\noindent We are again considering the screw geometry depicted with the bilinearly covered parametric domain before reparameterisation in Figure \ref{fig:male_screw_cell_k0}. Here, the reference solution corresponds to $\br(\bxi) = \bxi$.
\noindent Figure \ref{fig:male_screw_cell_hom} shows the geometry along with the associated controlmap $\bs(\bxi)$ after reparameterisation with $k=2$. Denoting the reparameterised maps by $\bx_h^k$, we define
\begin{align}
\label{eq:definition_nuArea_nudetJ}
    \nu_{\text{Area}}^k := \frac{L_{\text{Area}}(\bx_h^k)}{L_{\text{Area}}(\bx_h^0)} \quad \text{and} \quad \nu_{\det J}^k := \frac{\sup \partial_{\bmu} \bx_h^k}{\inf \partial_{\bmu} \bx_h^k},
\end{align}
where the latter is approximated by sampling over the abscissae of a dense Gauss-Legendre quadrature scheme. 
\begin{table}
\renewcommand{\arraystretch}{1.5}
\centering
\begin{tabular}{c|c|c|c|c}
 k & 0 & 1 & 2 & 3 \\
 \hline
$\nu_{\text{Area}}^k$ & 1 & 0.749 & 0.680 & 0.649 \\
\hline
$\nu_{\det J}^k$ & 23.3 & 6.67 & 4.22 & 5.69
\end{tabular}
\caption{Table showing the ratios between evaluating~\eqref{eq:Area_multipatch} in $\bx^k_h$ and the reference evaluation in $\bx_h^0$ for various values of $k$ as well as the ratio of the maximum and minimum values of $\det \partial_{\bmu} \bx_h^k$ sampled over a dense quadrature scheme.}
\label{tab:poisson_area_ratio}
\end{table}

\noindent Table \ref{tab:poisson_area_ratio} contains the values of $\nu_{\text{Area}}^k$ and $\nu_{\det J}^k$ for various values of $k \in [0, 3]$. The table clearly demonstrates that the methodology has the desired effect, reaching saturation for larger values of $k$. Furthermore, all parameterisations with $k > 0$ significantly reduce the anisotropy of $\det \partial_{\bmu} \bx_h^k$. We mention that the diffusivity is merely essentially bounded since $\sigma^{k}(\bx^{\bs})$ is generally patchwise discontinuous. This may lead to singularities or unbounded growth (for an example of a scalar elliptic problem in which a diffusivity that is a scaling times the identity creates a singularity, see \cite{morin2002convergence, da2014mathematical}). However, in practice singularities are avoided and the problem requires no stabilisation. A possible explanation is that cell size homogenisation counteracts the tendency to generate singularities for discrete approximations. \\

\begin{figure}[h!]
\centering
    \includegraphics[width=0.3 \textwidth, align=c]{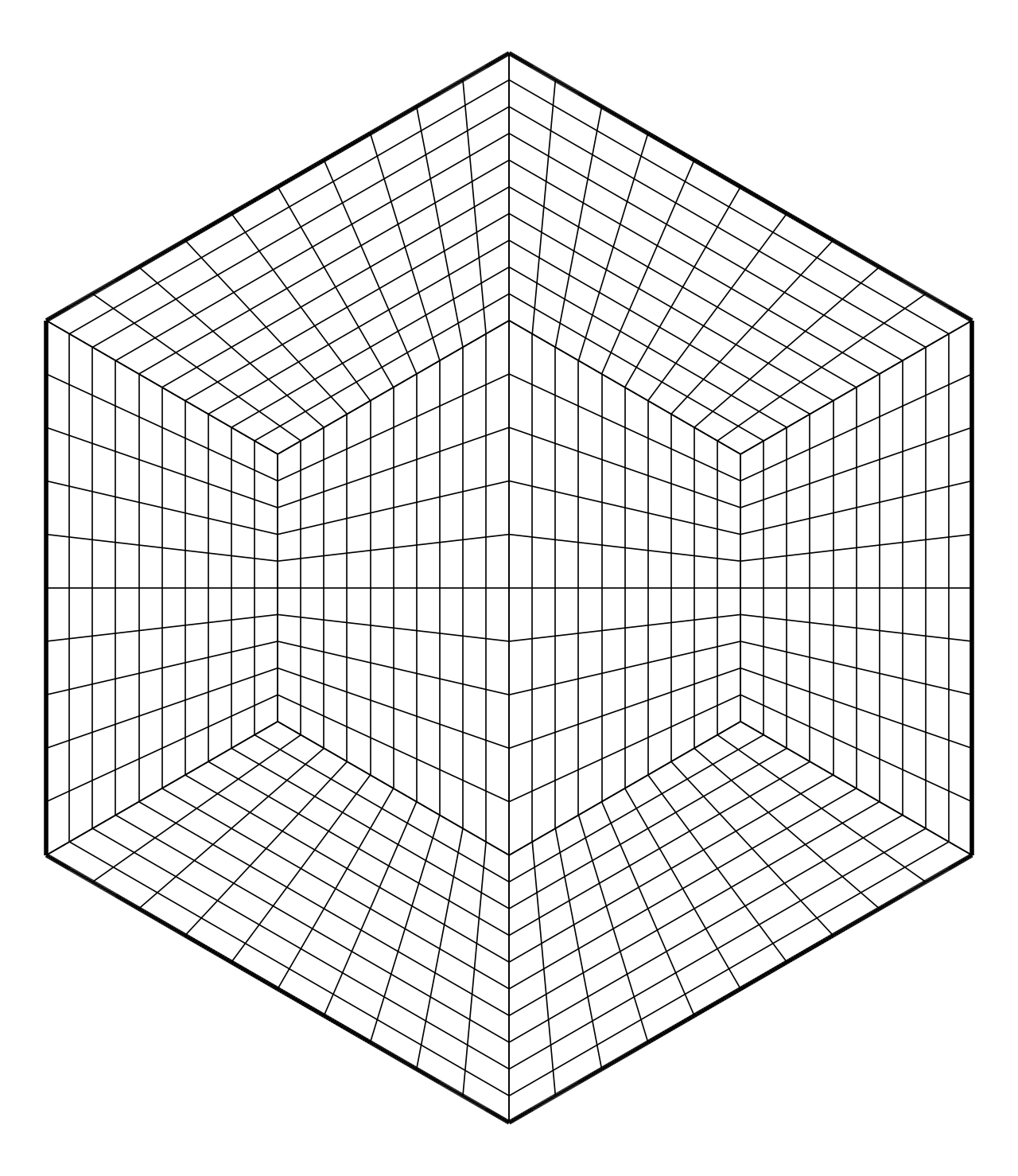}
    \includegraphics[width=0.3 \textwidth, align=c]{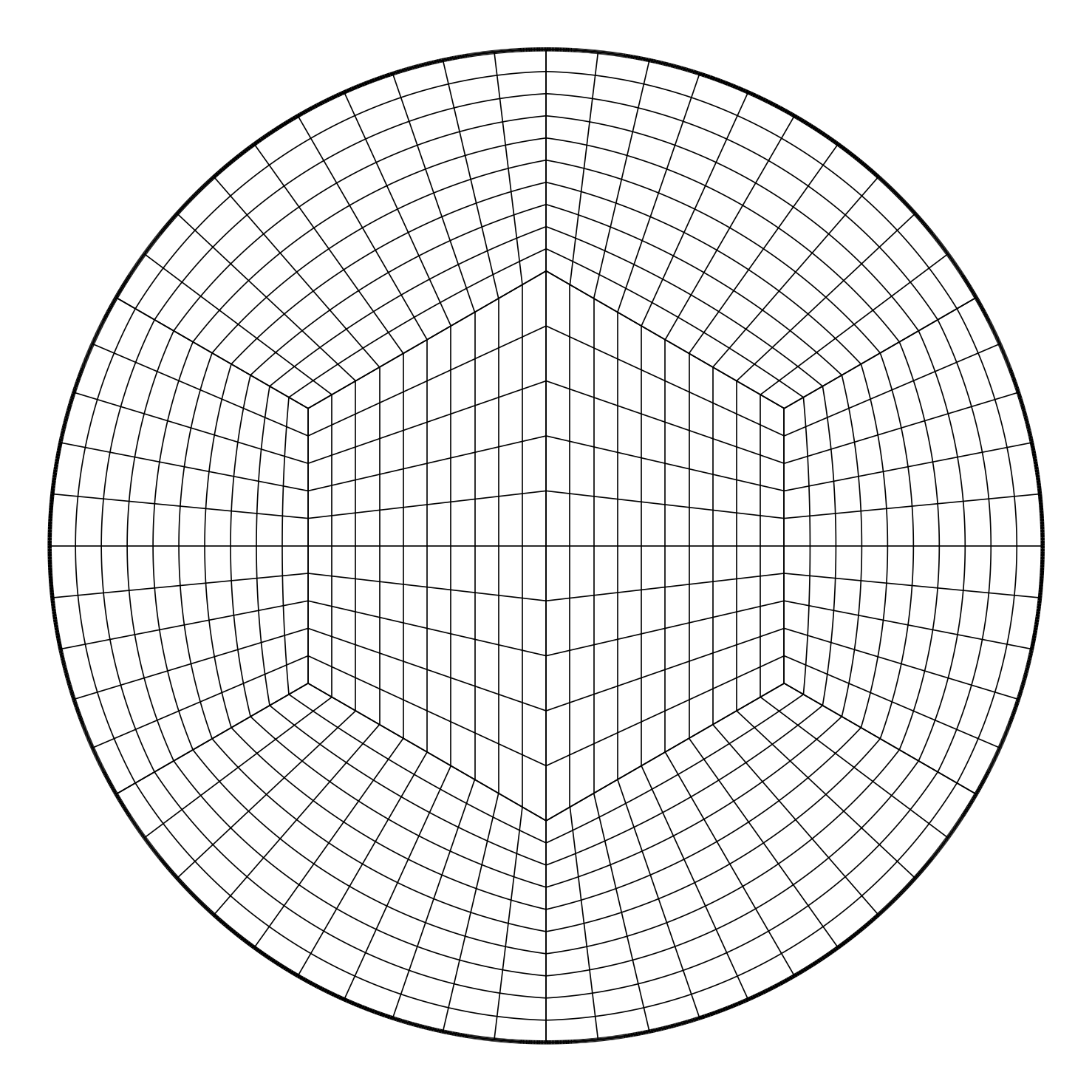}
    \includegraphics[width=0.3 \textwidth, align=c]{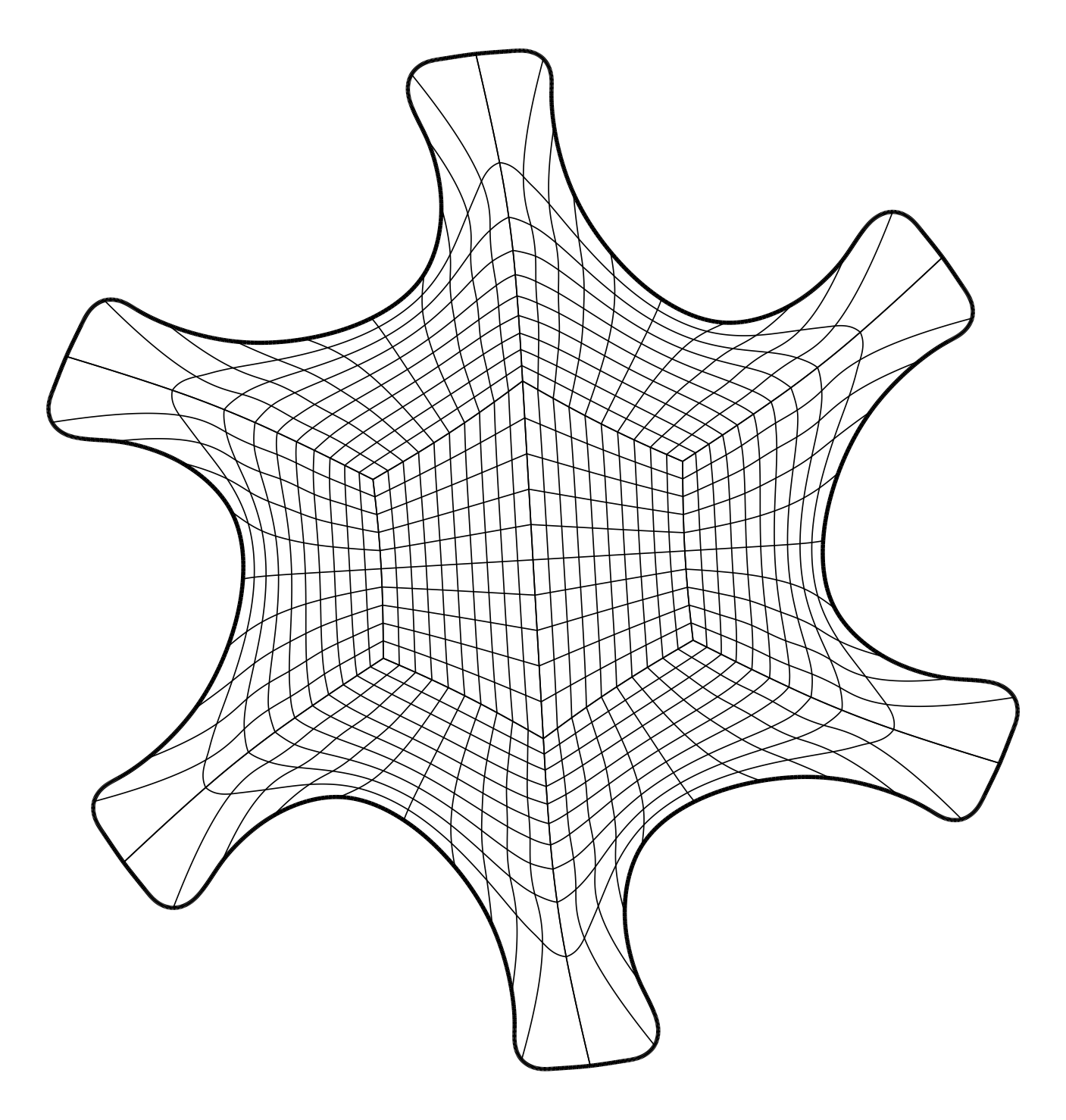}
\caption{Depiction of $\hOm$, the parameterisation of $\hOmr$ under $\widehat{\br}: \hOm \rightarrow \hOmr$ and the reference parameterisation of $\Om$ under $\bx^{\widehat{\br}}(\bxi): \hOm \rightarrow \hOmr$.}
\label{fig:female_screw_reference_domain}
\end{figure}

\begin{figure}[h!]
\centering
\begin{subfigure}[b]{0.8 \textwidth}
    \includegraphics[width=0.45 \textwidth]{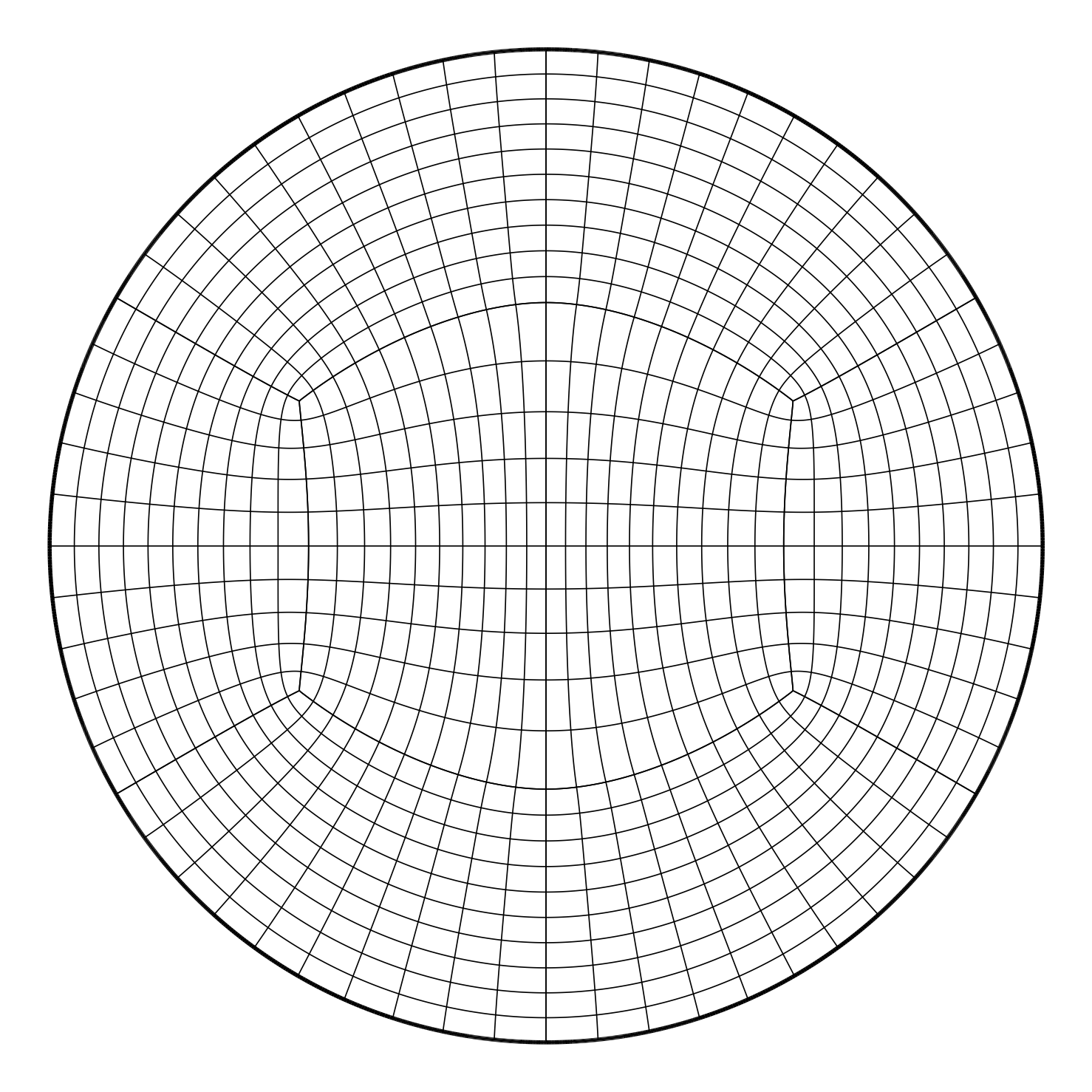} \hfill
    \includegraphics[width=0.45 \textwidth]{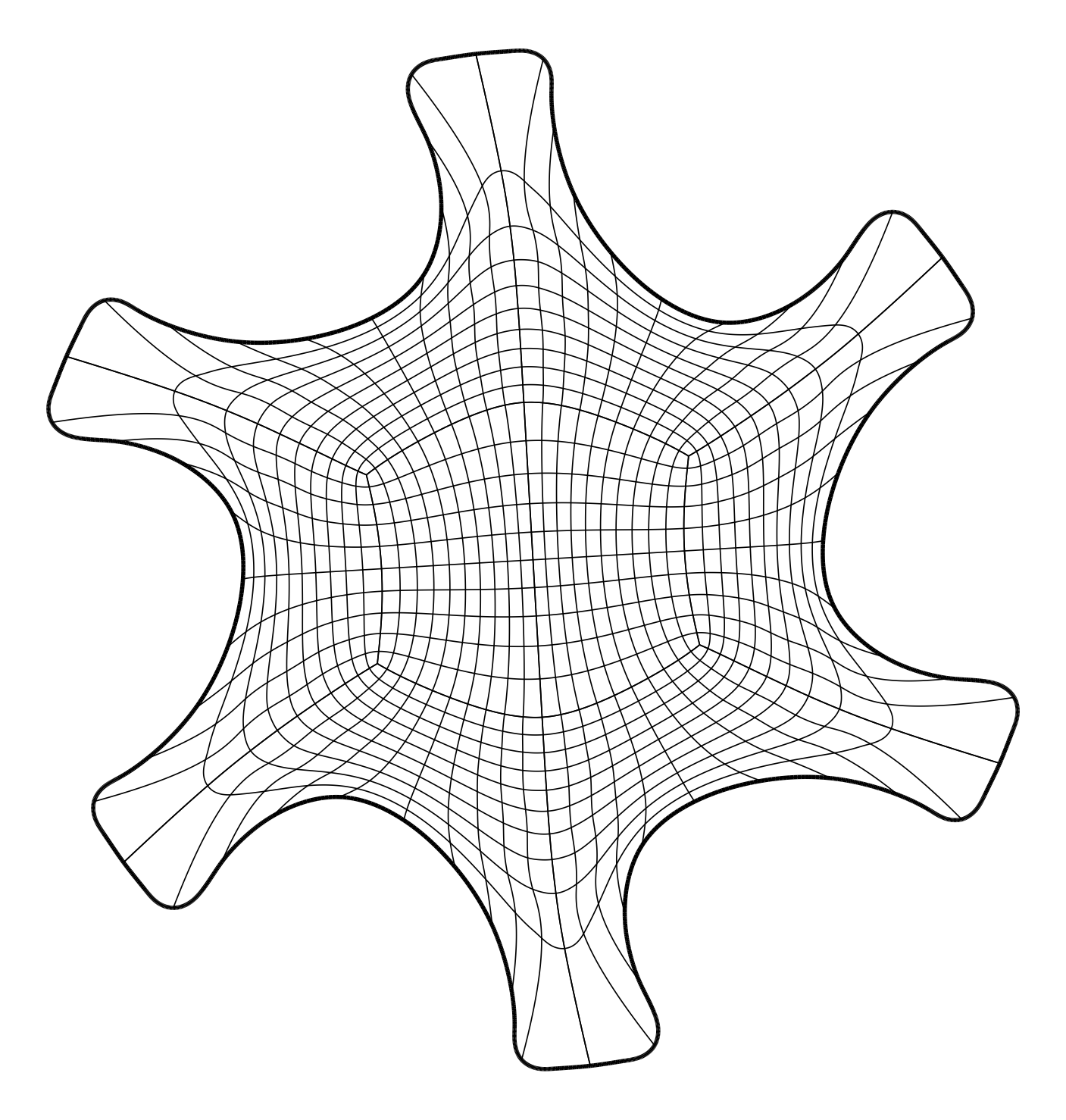}
    \caption{$\bs(\hOm^{\br})$ and the female screw geometry after reparameterisation using~\eqref{eq:cell_size_hom_diffusivity} with $k=0$.}
    \label{fig:female_screw_cell_k0}
\end{subfigure}
\\
\begin{subfigure}[b]{0.8 \textwidth}
    \includegraphics[width=0.45 \textwidth]{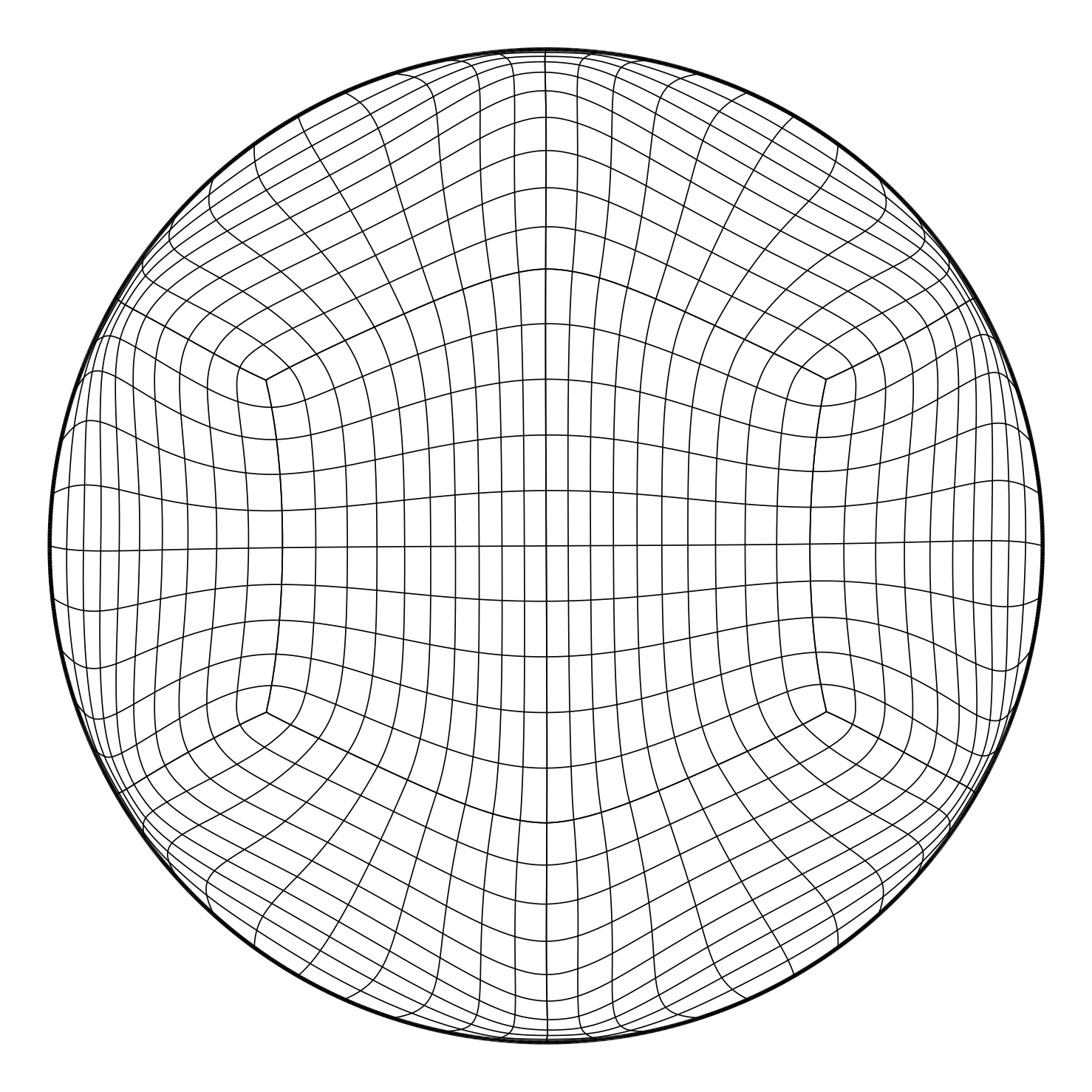} \hfill
    \includegraphics[width=0.45 \textwidth]{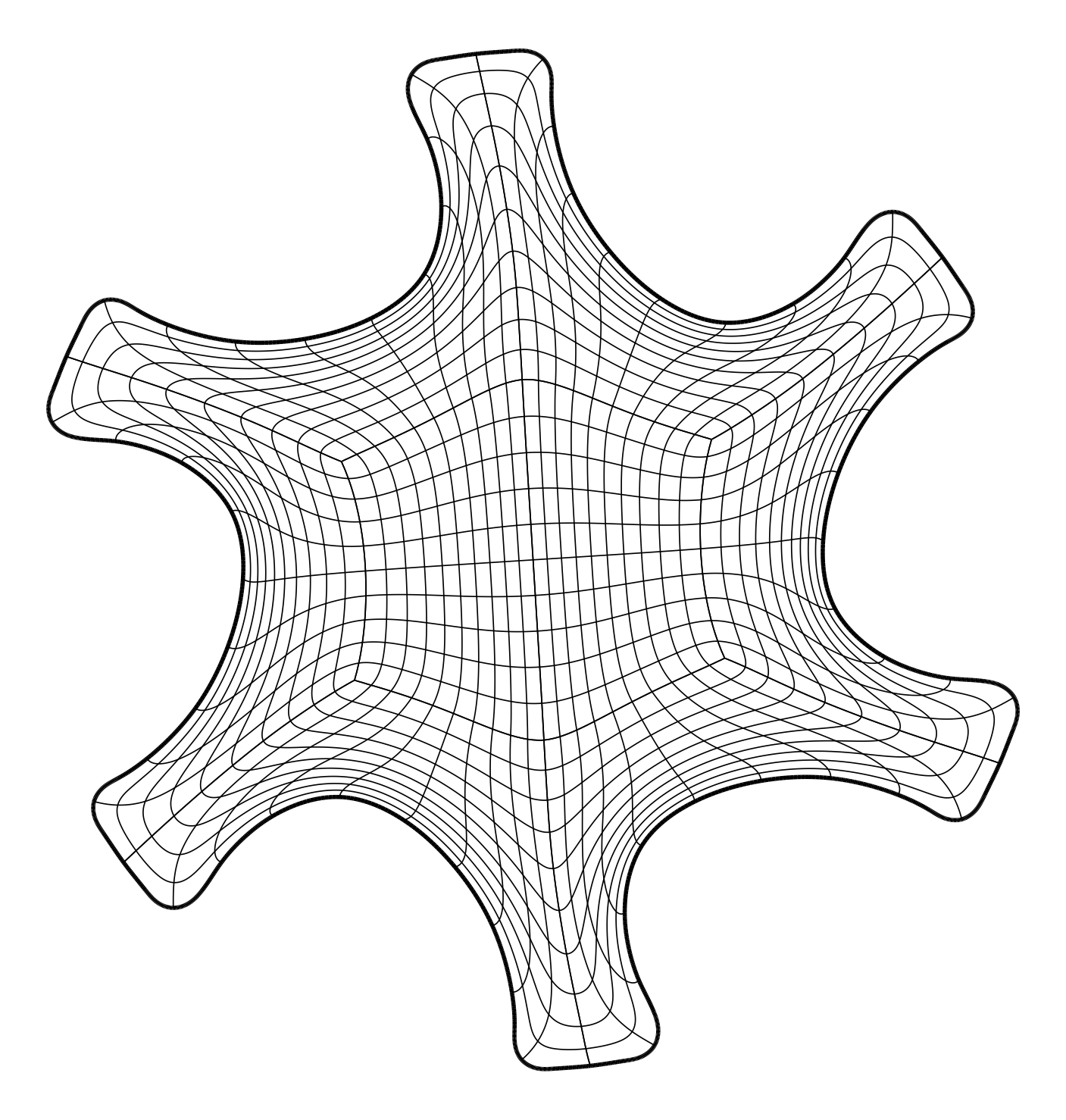}
    \caption{$\bs(\hOm^{\br})$ and the female screw geometry after reparameterisation using~\eqref{eq:cell_size_hom_diffusivity} with $k=4$.}
    \label{fig:female_screw_cell_k2}
\end{subfigure}
\caption{The female screw geometry after reparameterisation under~\eqref{eq:cell_size_hom_diffusivity} with $k=0$ and $k=4$.}
\label{fig:female_screw_cell_hom_local_trace_penalty}
\end{figure}

\noindent As a second example, we are considering the geometry depicted in Figure \ref{fig:female_screw_reference_domain} (right) along with the parametric domain $\hOm$ given by a regular six-sided polygon. Unlike the geometry from Figure \ref{fig:male_screw_cell_k0}, this geometry has no corners. As such, we take $\hOmr$ to be the unit disc where the boundary correspondence $\mathbf{F}^{\bxi \rightarrow \br}: \partial \hOm \rightarrow \partial \hOmr$ is chosen such that the induced correspondence $\mathbf{F}^{\br \rightarrow \bx}: \partial \hOmr \rightarrow \partial \Om$ is diffeomorphic. The interior of $\hOmr$ is parameterised by applying the bilinearly-blended Coons' patch approach to each patch of $\hOm$ individually. The reference parameterisation $\bx^{\widehat{\br}}(\bxi): \hOm \rightarrow \Om$ and the associated reference controlmap $\widehat{\br}: \hOm \rightarrow \hOmr$ are depicted in Figure \ref{fig:female_screw_reference_domain}.
\noindent In this example, we are combining the techniques from this section with those from Section \ref{subsect:interface_removal}. We employ the reference parameterisation $\widehat{\br}: \hOm \rightarrow \hOmr$ in combination with the diffusivity from~\eqref{eq:laplace_local_normalised_trace_penalty} to map $\hOmr$ onto itself. This results in the new reference controlmap $\br: \hOm \rightarrow \hOmr$ that removes the patch interfaces. The new parameterisation of $\hOmr$ and the associated new reference map $\bx^{\br}: \hOm \rightarrow \hOmr$ are depicted in Figure \ref{fig:female_screw_cell_k0}. Here, we do not stabilise using~\eqref{eq:patch_vertex_stabilisation} since the discrete approximation remains uniformly nondegenerate with acceptable behaviour in the vicinity of the patch vertices.
\begin{table}[h!]
\renewcommand{\arraystretch}{1.5}
\centering
\begin{tabular}{c|c|c|c|c|c}
 k & 0 & 1 & 2 & 3 & 4 \\
 \hline
$\nu_{\text{Area}}^k$ & 1.02 & 0.496 & 0.414 & 0.386 & 0.373  \\
\hline
$\nu_{\Gamma}^k$ & 0.117 & 0.142 & 0.158 & 0.171 & 0.180
\end{tabular}
\caption{Table showing the ratios between evaluating~\eqref{eq:Area_multipatch} in $\bx^k_h$ and the reference evaluation in $\bx_h^0$ for various values of $k$.}
\label{tab:female_poisson_area_interface_ratio}
\end{table}

Table \ref{tab:female_poisson_area_interface_ratio} contains the values of
$$\nu^k_{\text{Area}} := \frac{L_{\text{Area}}(\bx^k_h)}{L_{\text{Area}}(\bx^{\widehat{\br}}_h)} \quad \text{and} \quad \nu^k_{\Gamma} := \frac{L_{\Gamma}(\bx_h^k)}{L_{\Gamma}(\bx^{\widehat{\br}}_h)}, \quad \text{with } L_{\Gamma}(\, \cdot \,) \text{ as in~\eqref{eq:interface_removal_norm}},$$
for various $k \in [0, 4]$. Figure \ref{fig:female_screw_cell_hom_local_trace_penalty} shows the parameterisation of the geometry for $k = 0$ and $k=4$. The table clearly demonstrates a monotonous reduction of $\nu^k_{\text{Area}}$ (eventually reaching saturation for larger values of $k$) at the expense of a slight increase in the value of $\nu^k_{\Gamma}$. \\

For more precise control over the expansion / contraction of cells, it can be helpful to decompose the diffusivity into the scaling $\sigma^k(\bx^{\bs})$ times the sum of two symmetric rank one tensors, i.e.,
\begin{align}
\label{eq:poisson_trace_penalty_general}
    D(\bx^{\bs}) = \frac{2}{a_i + 1} \sigma^k(\bx^{\bs}) \left( a_i \hat{\mathbf{v}}^{i, 1} \otimes \hat{\mathbf{v}}^{i, 1} + \hat{\mathbf{v}}^{i, 2} \otimes \hat{\mathbf{v}}^{i, 2} \right), \quad \text{ on } \br(\hat{\Omega}^i), \quad \text{where} \quad a_i > 0.
\end{align}
Here, $\hat{\mathbf{v}}^{i, 1}$ and $\hat{\mathbf{v}}^{i, 2}$ have length one and are not parallel. Note that
$$\operatorname{tr}\left(\frac{2}{a_i + 1}\left( a_i \hat{\mathbf{v}}^{i, 1} \otimes \hat{\mathbf{v}}^{i, 1} + \hat{\mathbf{v}}^{i, 2} \otimes \hat{\mathbf{v}}^{i, 2} \right) \right) = 2, \quad \text{as before}.$$
If the $\hat{\mathbf{v}}^{i, j}$ are patchwise discontinuous, we note that the diffusivity may require stabilisation. Taking $a_i$ large will force $\bs: \hOmr \rightarrow \hOmr$ to predominantly slide in the direction of $\mathbf{v}^{i, 1}$ on $\br(\hat{\Omega}^i)$.
\begin{figure}[h!]
\centering
\begin{subfigure}[b]{0.95 \textwidth}
    \includegraphics[width=0.3 \textwidth]{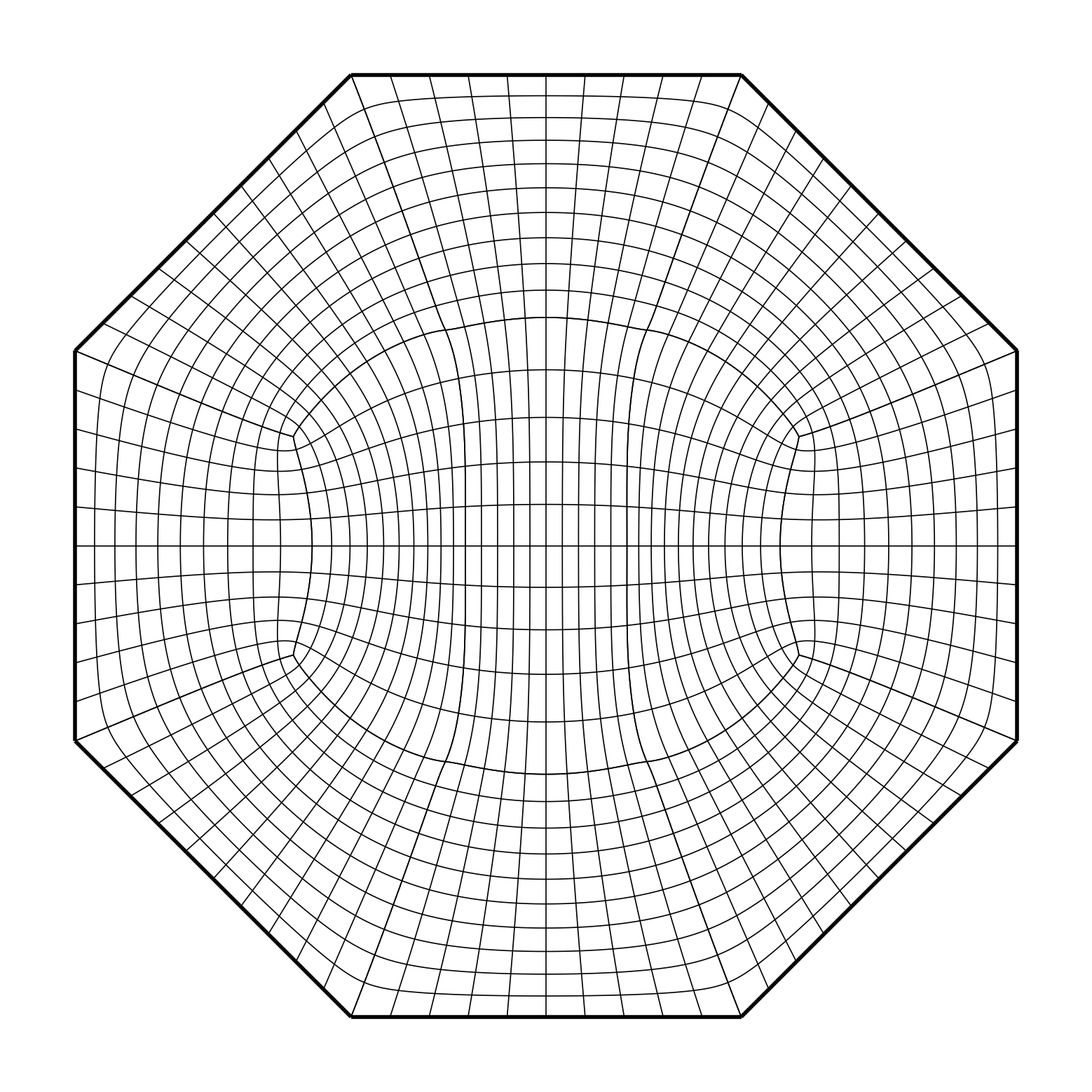} \hfill
    \includegraphics[width=0.3 \textwidth]{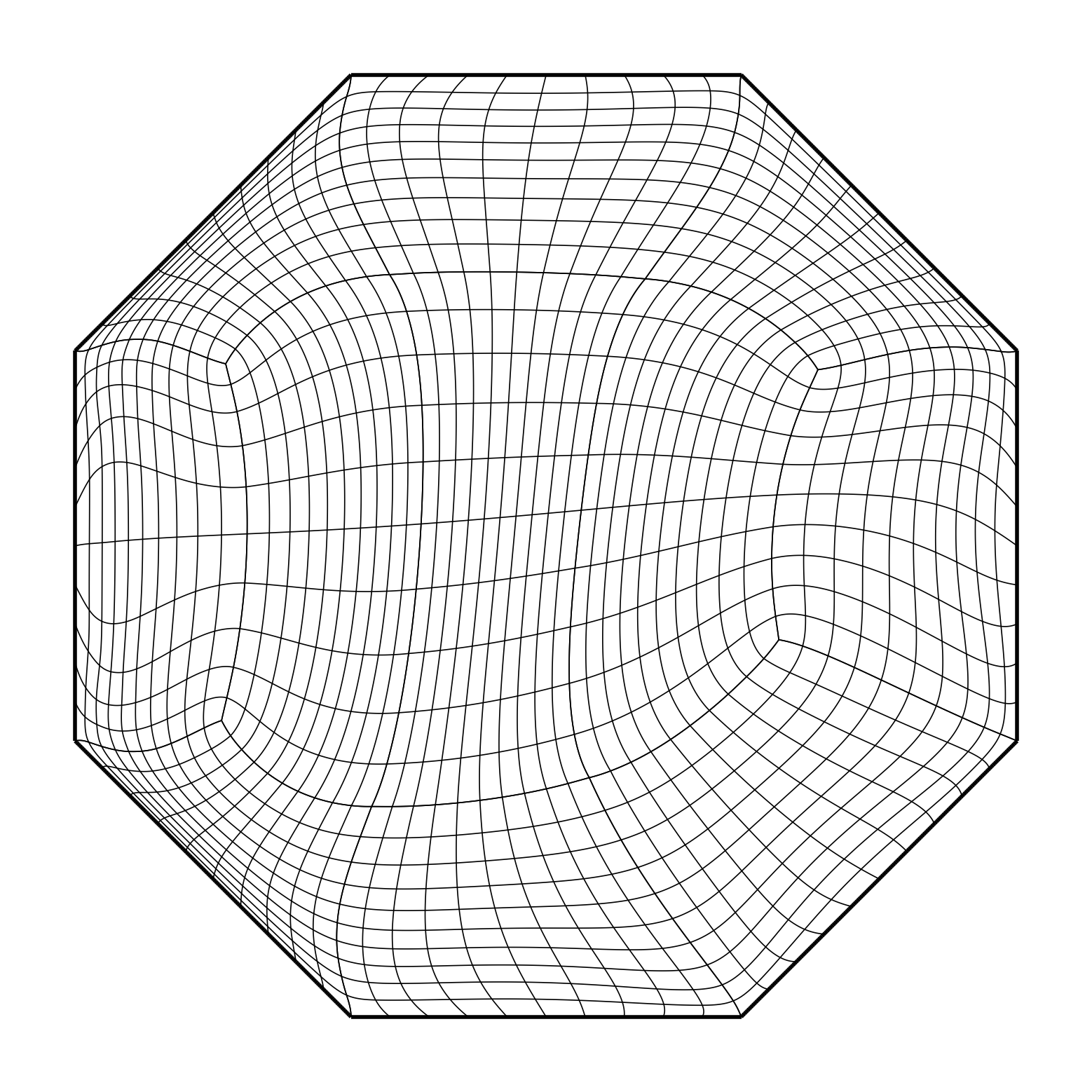}
    \hfill
    \includegraphics[width=0.3 \textwidth]{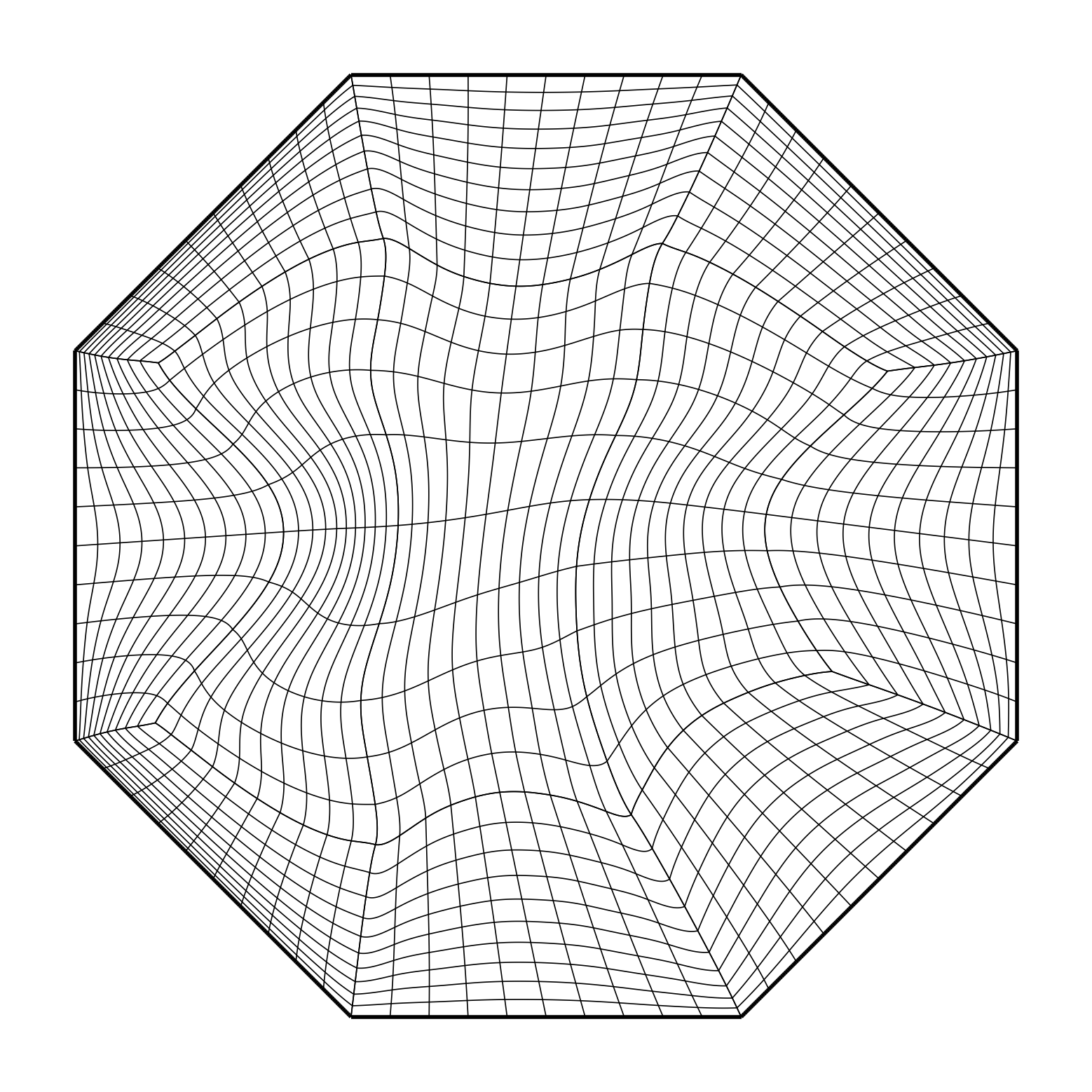}
    \caption{Visualisation of the reference controlmap (left) and the controlmaps $\bs: \hOmr \rightarrow \hOmr$ using $D(\br, \bx^{\bs}) = \sigma(\bx^{\bs})^k \mathcal{I}^{2 \times 2}$ and~\eqref{eq:poisson_trace_penalty_general}, both with $k = 3.5$. In the latter, $a_i$ is taken large close to $\partial \hOmr$ while $\hat{\mathbf{v}}^{i, 1}$ is the normalised component of $\partial_{\bmu} \br$ directed transversal to $\partial \hOmr$ for boundary patches.}
    \label{fig:throwing_star_cell_size_hom_param}
\end{subfigure}
\\
\begin{subfigure}[b]{0.95 \textwidth}
    \includegraphics[width=0.3 \textwidth]{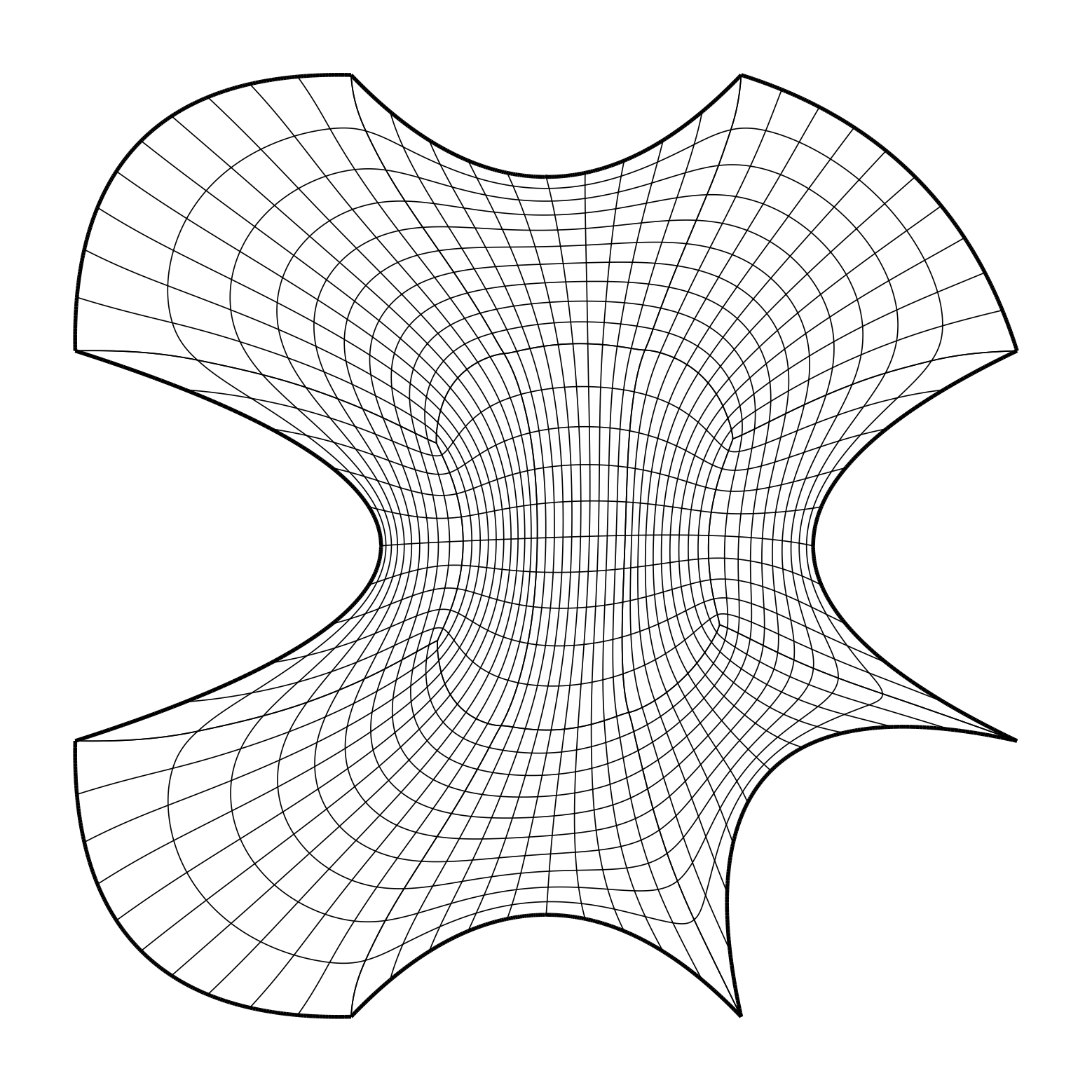} \hfill
    \includegraphics[width=0.3 \textwidth]{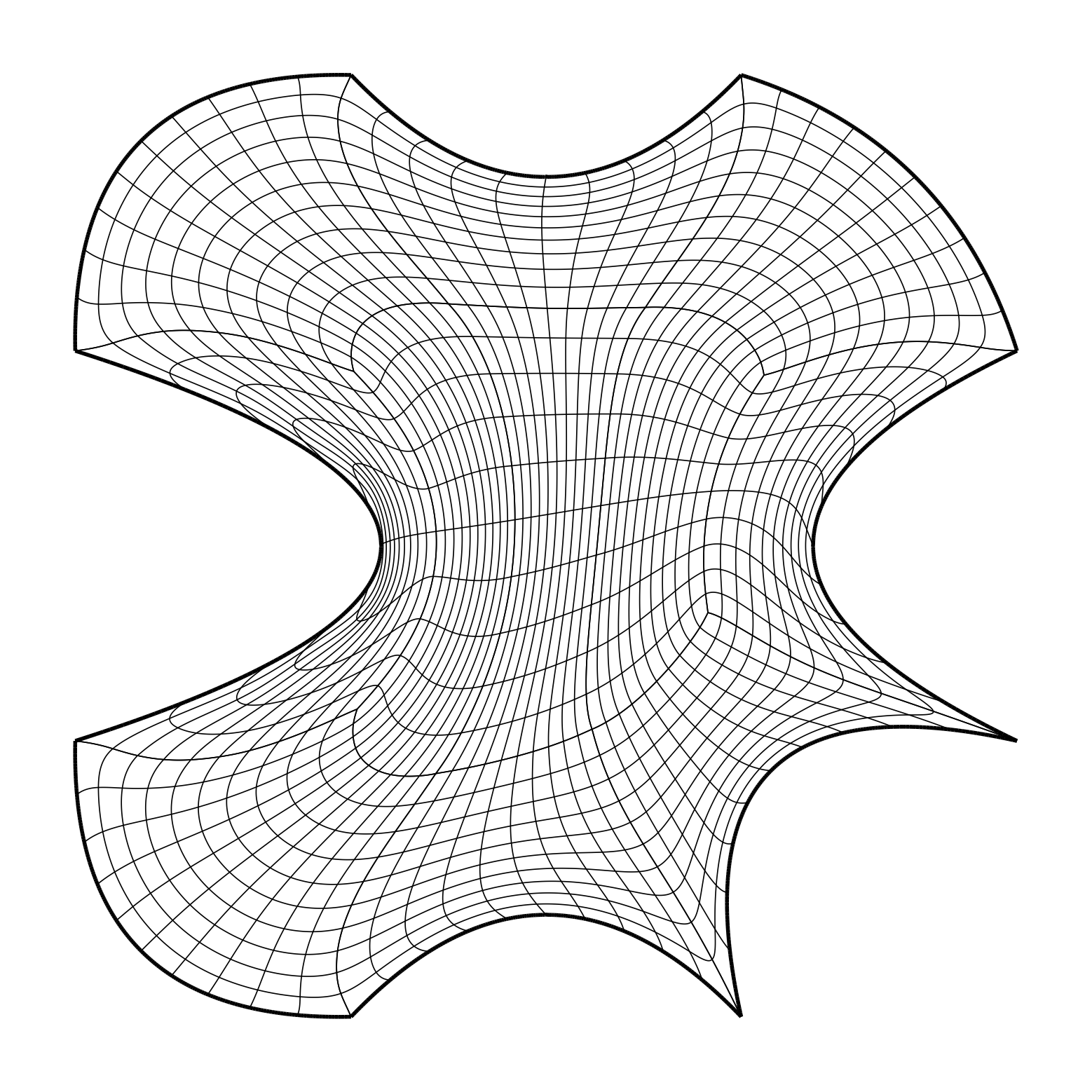} \hfill
    \includegraphics[width=0.3 \textwidth]{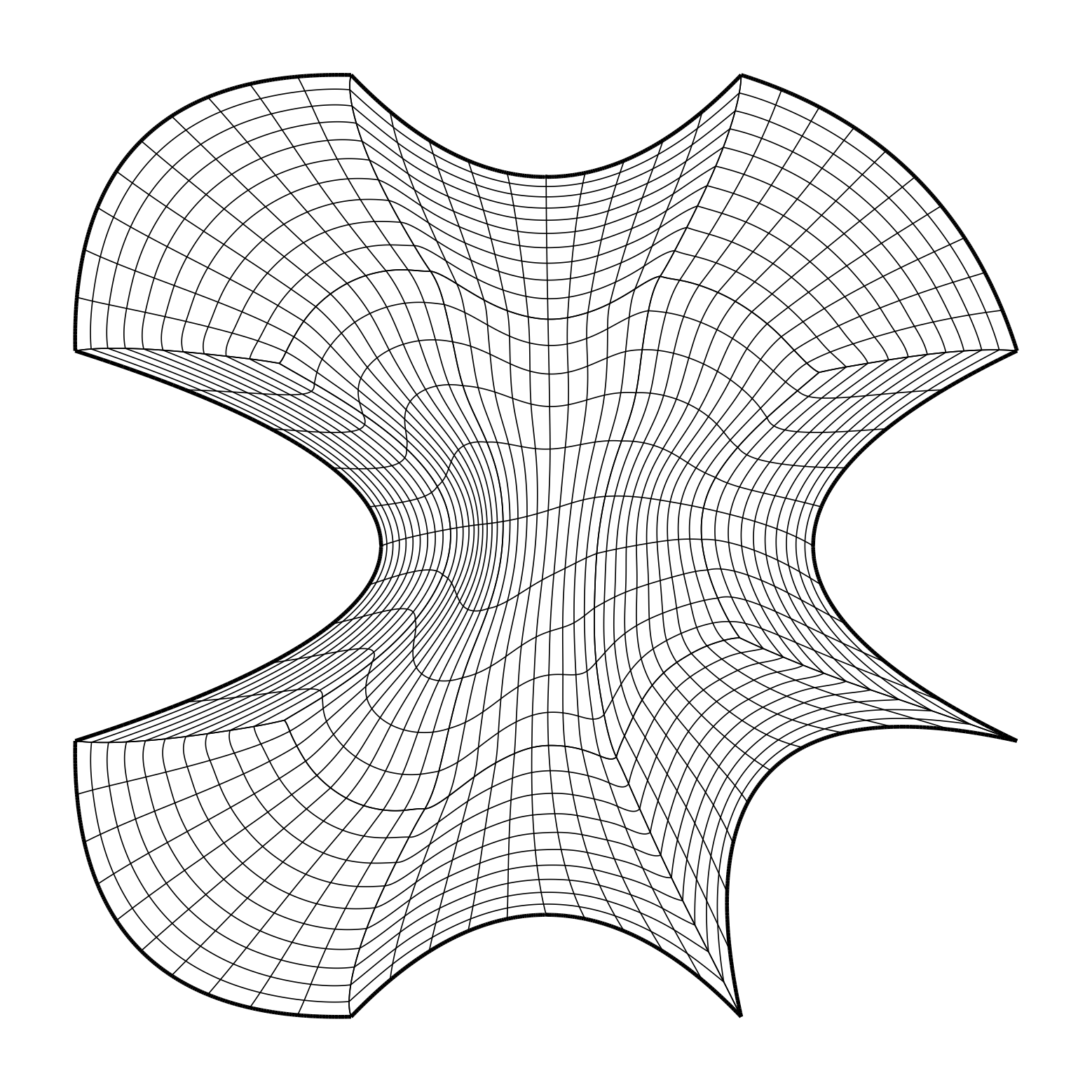}
    \caption{The associated parameterisations. The parameterisations under Figure \ref{fig:throwing_star_cell_size_hom_param} (center) and Figure \ref{fig:throwing_star_cell_size_hom_param} (right) both show an outstanding homogenisation. The former exhibits a large degree of cell skewness close to the boundary, especially for the center left patch. Meanwhile, the latter avoids the boundary skewness but increases the isoline angles by the patch interfaces.}
    \label{fig:throwing_star_cell_size_hom_geom}
\end{subfigure}
\caption{Reparameterisation of a irregular geometry using $D(\br, \bx^{\bs}) = \sigma(\bx^{\bs})^k \mathcal{I}^{2 \times 2}$ (center) and~\eqref{eq:poisson_trace_penalty_general} (right). The latter avoids cell skewness close to the boundary. }
\label{fig:throwing_star_cell_size_hom}
\end{figure}
For large values of $k$, cell size homogenisation can lead to a considerable degree of cell skewness close to the boundary. As an example, Figure \ref{fig:throwing_star_cell_size_hom} (center) shows the homogenisation of a geometry with reference parameterisation depicted in Figure \ref{fig:throwing_star_cell_size_hom} (left) using~\eqref{eq:poisson_trace_penalty_general} under $D^{\bs}(\bx^{\bs}) = \sigma^k(\bx^{\bs}) \mathcal{I}^{2 \times 2}$ with $k=3.5$. This effect can be avoided by setting $\hat{\mathbf{v}}^{i, j} = \widehat{\partial}_{\bmu_j} \br$, where $\hat{\mathbf{v}}^{i, 1} = \widehat{\partial}_{\bmu^{\perp}} \br$ for boundary patches, while taking the $a_i$ large close to the boundary of $\hOmr$. The result for $k = 3.5$ is depicted in Figure \ref{fig:throwing_star_cell_size_hom} (right). With $\nu^k_{\text{Area}} = 0.515$ for the former and $\nu^k_{\text{Area}} = 0.518$ for the latter, the homogenisation is only marginally less effective under~\eqref{eq:poisson_trace_penalty_general}. However, the latter completely avoids cell skewness close to the boundary while sacrificing some regularity across the patch interfaces. \\

\begin{figure}[h!]
\centering
    \includegraphics[width=0.4 \textwidth, align=c]{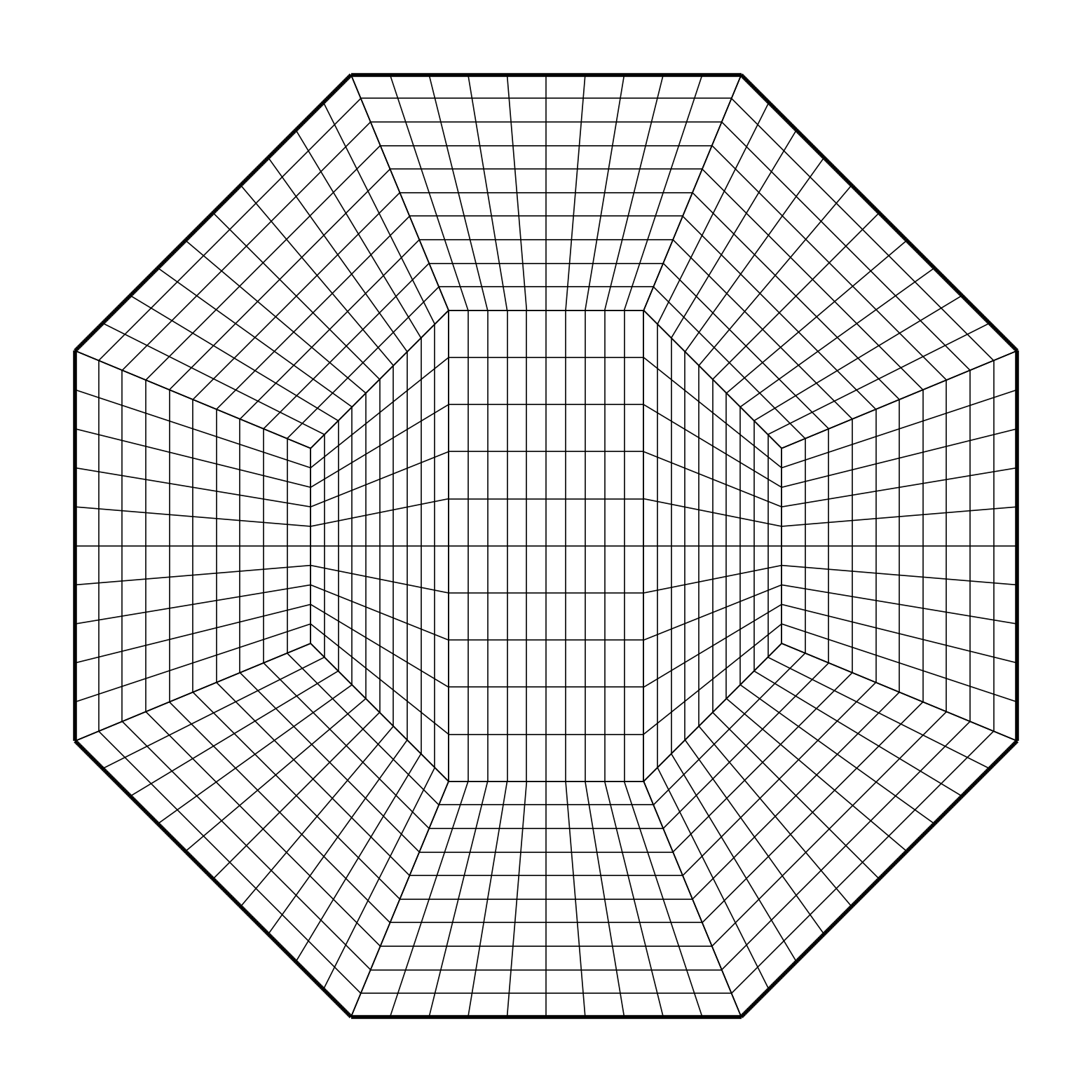}
    \includegraphics[width=0.4 \textwidth, align=c]{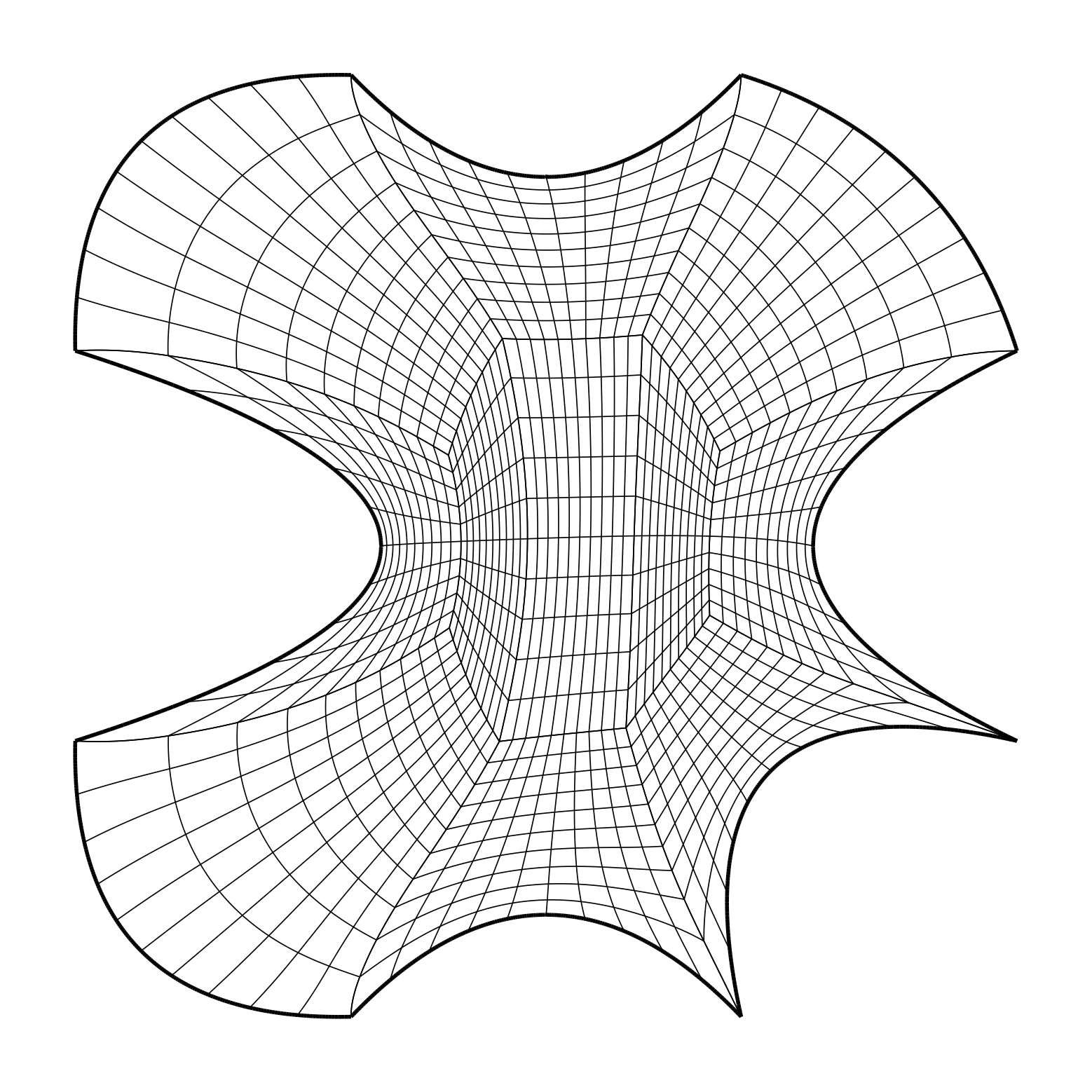}
\caption{Figure depicting the parameterisations of the tuple $(\hOm, \Om)$ for $D^{\bx} = D^{\bs} = \mathcal{I}^{2 \times 2}$.}
\label{fig:throwing_star_cell_size_hom_vanilla}
\end{figure}

\begin{figure}[h!]
\centering
    \includegraphics[width=0.3 \textwidth, align=c]{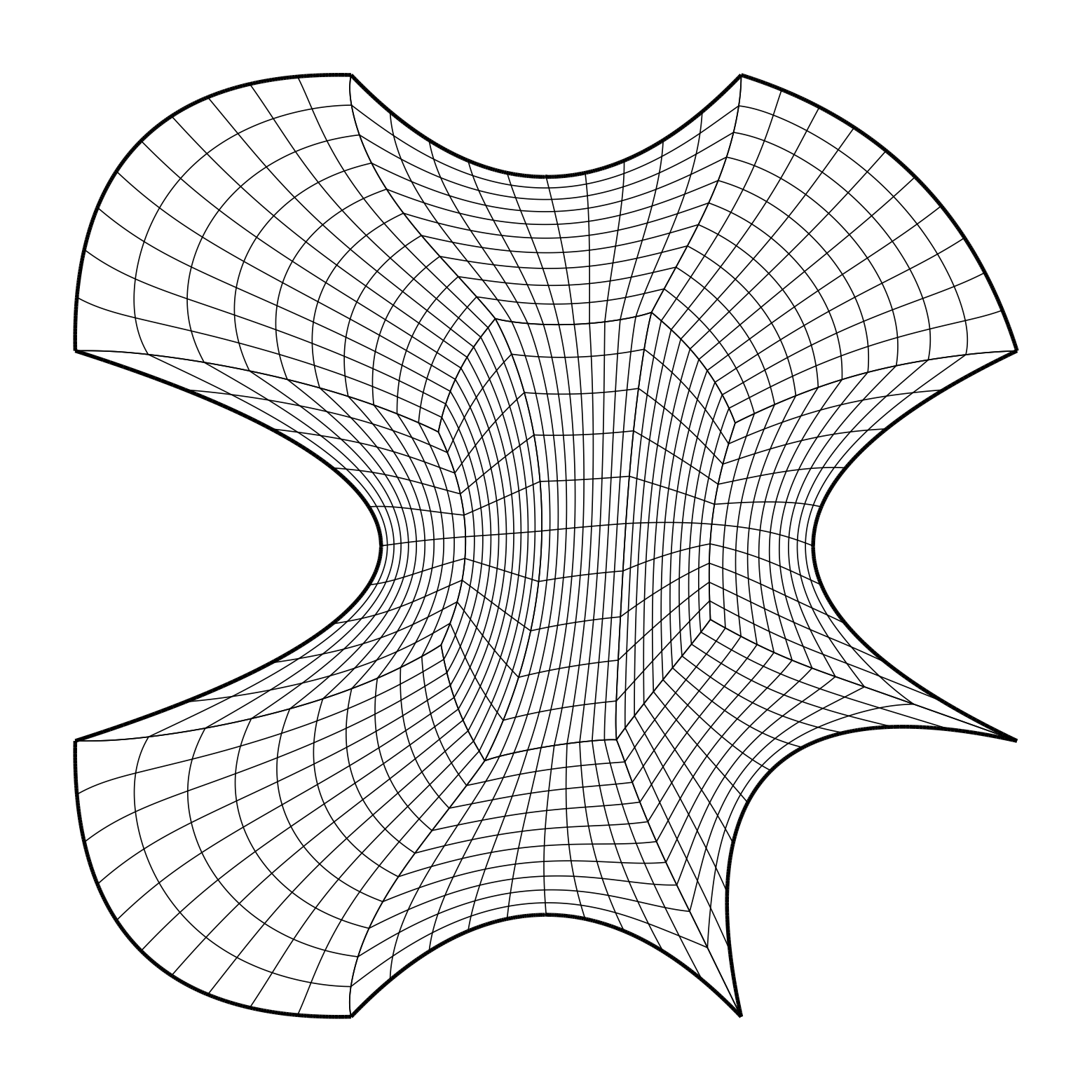}
    \includegraphics[width=0.3 \textwidth, align=c]{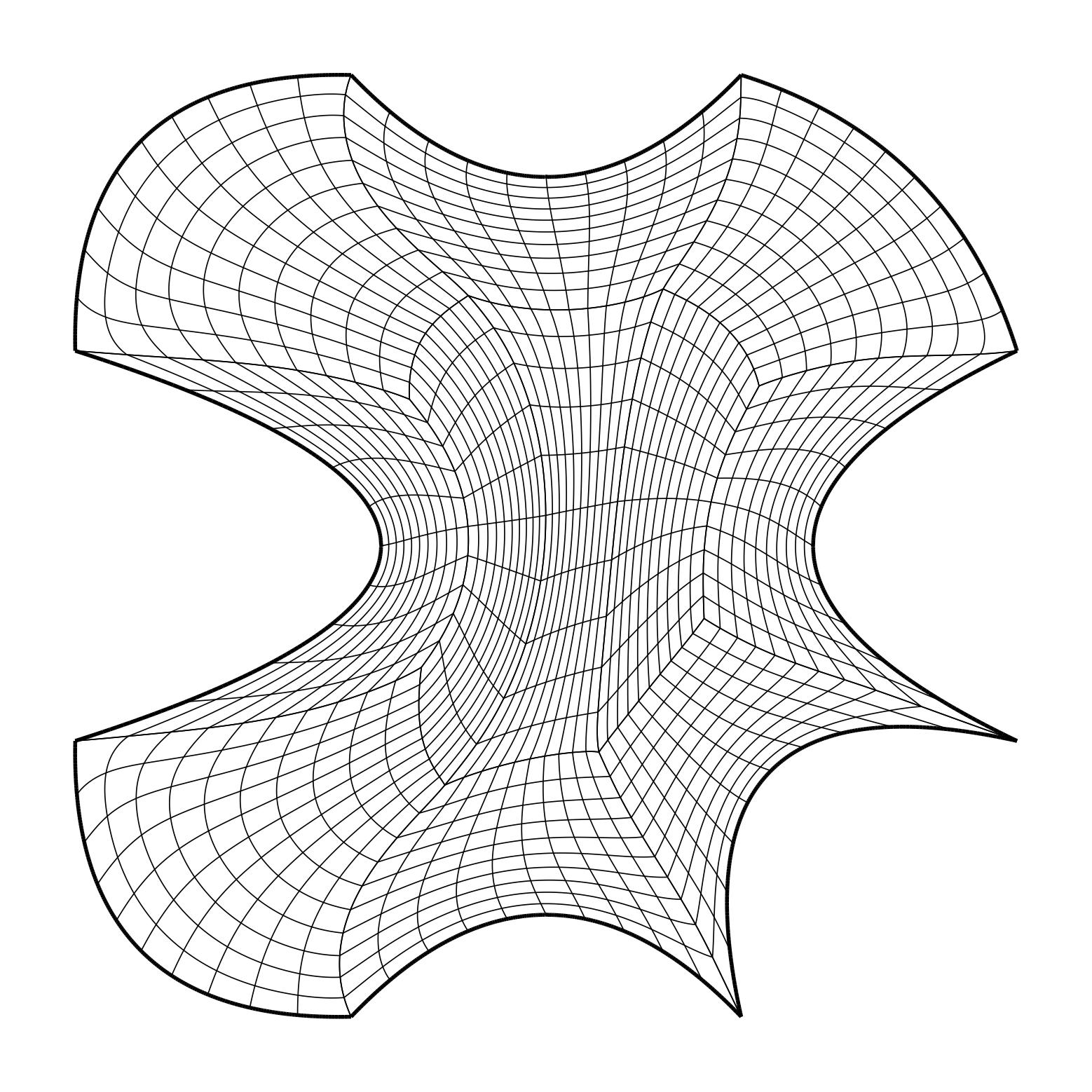}
    \includegraphics[width=0.3 \textwidth, align=c]{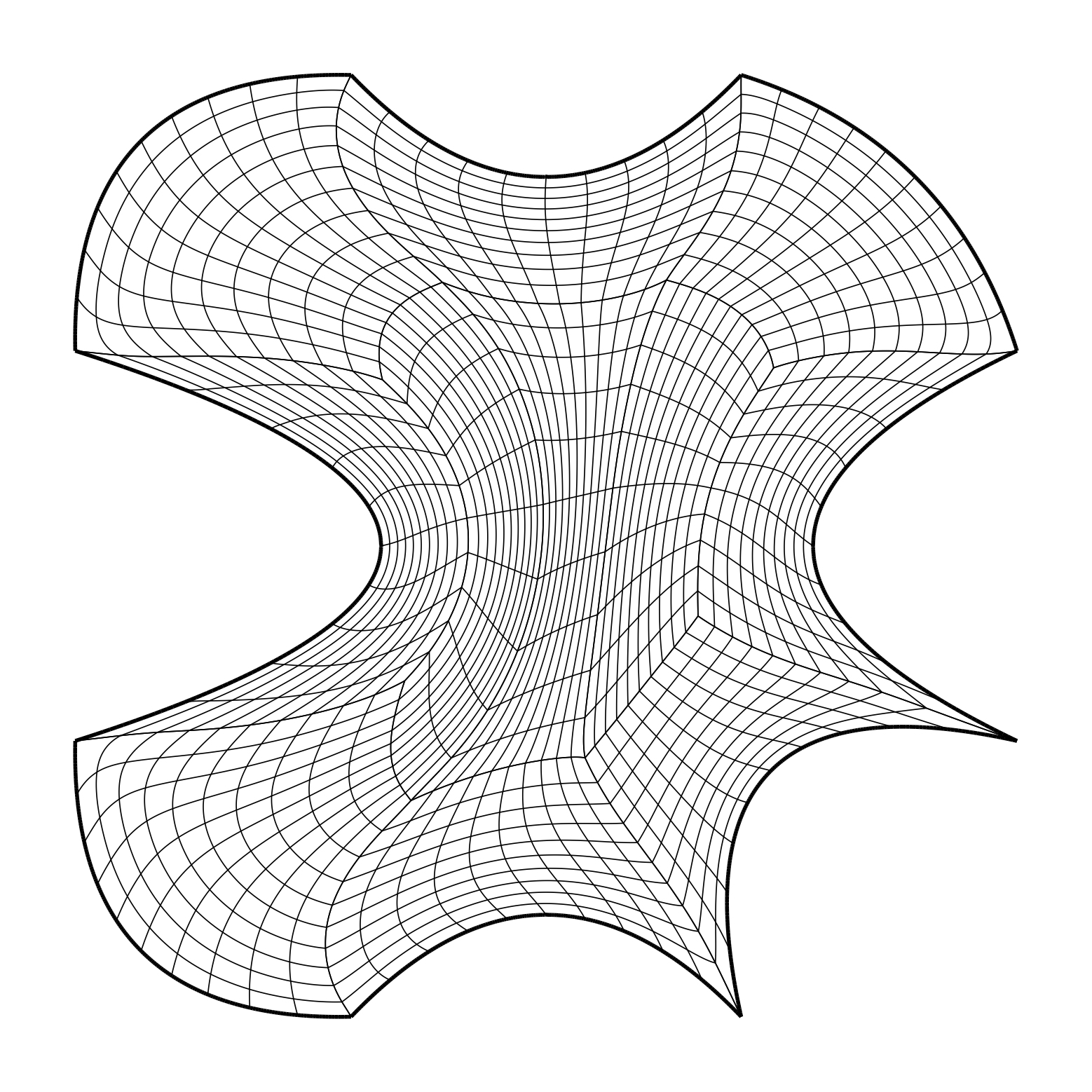}
\caption{Figure depicting the reparameterisations of the reference parameterisation from Figure \ref{fig:throwing_star_cell_size_hom_vanilla} under $\omega^k(\bx^{\bs})$ for $k=1$, $k=4$ and $k=8$.}
\label{fig:throwing_star_cell_size_hom_direct}
\end{figure}

\noindent In contrast to method 1., method 2. requires $D^{\bs}$ to be a function of $\br$ only and, for convenience, we assume $D^{\bs} = \mathcal{I}^{2 \times 2}$ such that $\bs(\br) = \br$. The cell size homogenisation is now encouraged through a proper choice of $D^{\bx}(\bs, \bx^{\bs})$. Similar to method 1., we take $D^{\bx}(\bs, \bx^{\bs}) = \omega(\bx) \mathcal{I}^{2 \times 2}$, for some $\omega(\bx) > 0$. Method 2. has the advantage of decoupling the system from~\eqref{eq:domain_optimisation_coupled_equation}. This reduces the problem size of the iterative root-finding algorithm, which now computes only $\bx^{\bs}$ instead of the tuple $(\bx^{\bs}, \bs)$, with $\bs(\bxi) = \br(\bxi)$. \\
This choice of $D^{\bx}(\bs, \bx^{\bs})$ encourages the contraction of $\br(\bx)$ isolines in $\hOmr$ wherever $\omega(\bx)$ is large and vice-versa. Exchanging the dependency $\br(\bx) \rightarrow \bx(\br)$, the isolines will now be contracted in regions where $\omega(\bx)$ is small. Inspired by method 1., we define the family monitor functions
$$
\omega(\bx^{\bs})^k(\br) := \left(\det \partial_{\bmu} \bx^{\bs}(\br) \right)^{-k} \quad \text{on} \quad \br(\hOm_i).
$$
As such, we are solving the decoupled system with $D^{\bx}(\bx^{\bs}) = \omega^{k}(\bx^{\bs}) \mathcal{I}^{2 \times 2}$. Clearly, for a root-finding algorithm to converge, the value of $\det \partial_{\bmu} \bx^{\bs}(\br)$ has to stay positive. 
As before, the barrier property of~\eqref{eq:weak_form_operator_controlmap} prevents intermediate iterates from leaving the set of UNDG maps and the scheme converges reliably for a wide range of $k > 0$ when initialised with the solution of one of the NDF formulations. We are considering the geometry with reference parameterisation from Figure \ref{fig:throwing_star_cell_size_hom_vanilla}. The same figure shows the bilinearly covered parametric domain.
We are monitoring the values of $\nu^{k}_{\text{Area}}$ and $\nu^{k}_{\det J}$ (cf.~\eqref{eq:definition_nuArea_nudetJ}) for $k \in \{0, \ldots, 8\}$. Table \ref{tab:poisson_area_ratio_throwing_star} contains the associated values, while Figure \ref{fig:throwing_star_cell_size_hom_direct} depicts the homogenised parameterisations for three different values of $k$.
\begin{table}
\renewcommand{\arraystretch}{1.5}
\centering
\begin{tabular}{c|c|c|c|c|c|c|c|c|c}
 k & 0 & 1 & 2 & 3 & 4 & 5 & 6 & 7 & 8 \\
 \hline
$\nu_{\text{Area}}^k$ & 1 & 0.716 & 0.623 & 0.583 & 0.558 & 0.542 & 0.532 & 0.524 & 0.518 \\
\hline
$\nu_{\det J}^k$ & 54.7 & 23.1 & 15.6 & 12.2 & 10.2 & 8.87 & 7.88 & 6.98 & 6.05
\end{tabular}
\caption{Table showing the values of $\nu_{\text{Area}}^k$ and $\nu_{\det J}^k$ as defined in~\eqref{eq:definition_nuArea_nudetJ} after reparameterising the reference parameterisation from Figure \ref{fig:throwing_star_cell_size_hom_vanilla} under $\omega^k(\bx^{\bs})$ for various values of $k \geq 0$.}
\label{tab:poisson_area_ratio_throwing_star}
\end{table}
The table clearly demonstrates that the methodology is highly effective homogenising the cell sizes in the local $\bmu$ coordinate systems. We also observe a significant reduction in the anisotropy of $\det \partial_{\bmu} \bx_h^{k}$, which is reduced from the initial $\nu_{\det J}^k = 54.7$ to $\nu_{\det J}^k = 6.05$ for $k=8$.

\subsection{Grid Adaptation}
\label{subsect:grid_adaptation}
\begin{figure}[h!]
\centering
\includegraphics[height=7.5cm]{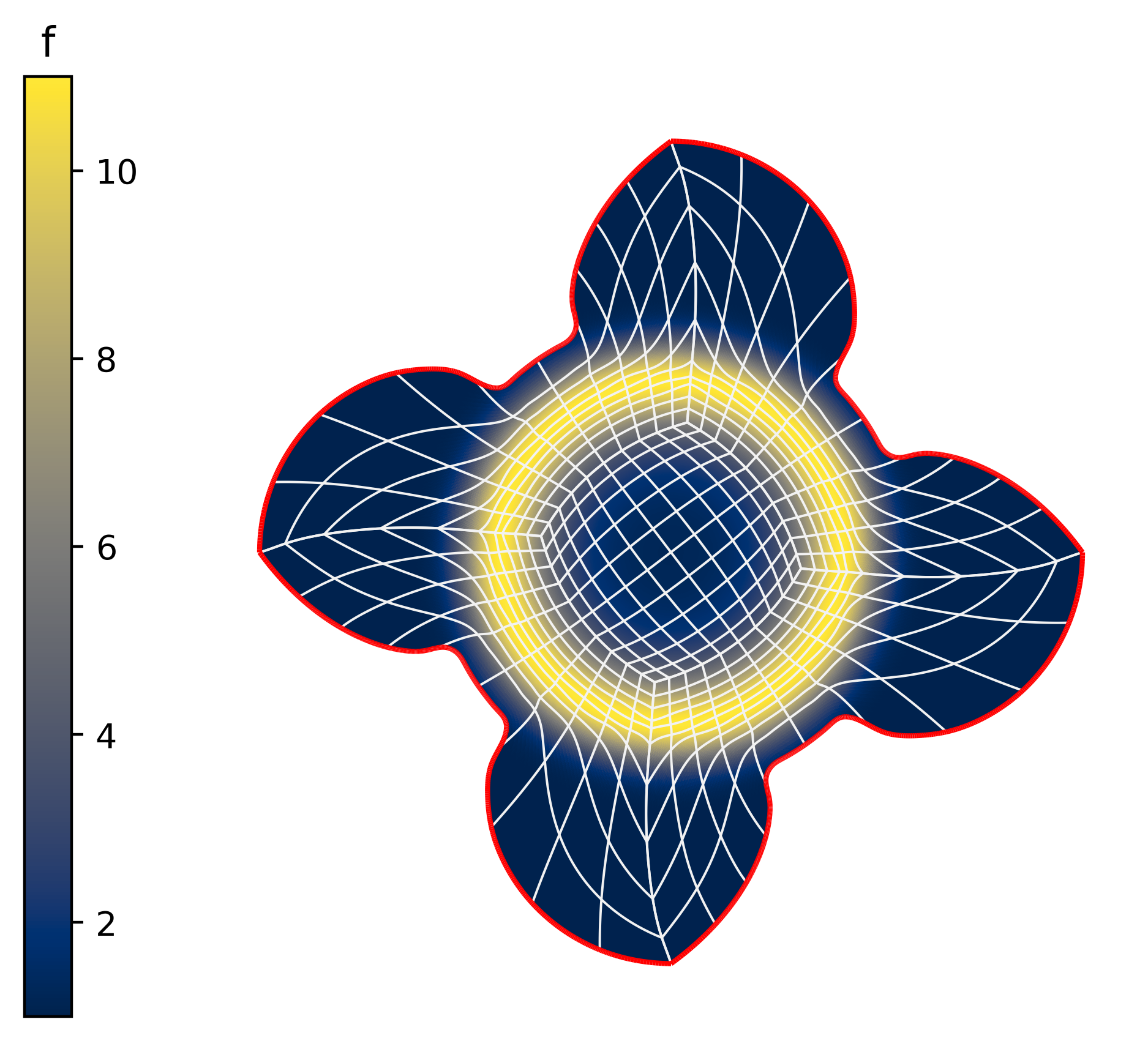} \hspace*{0.1cm}
\includegraphics[height=6.8cm]{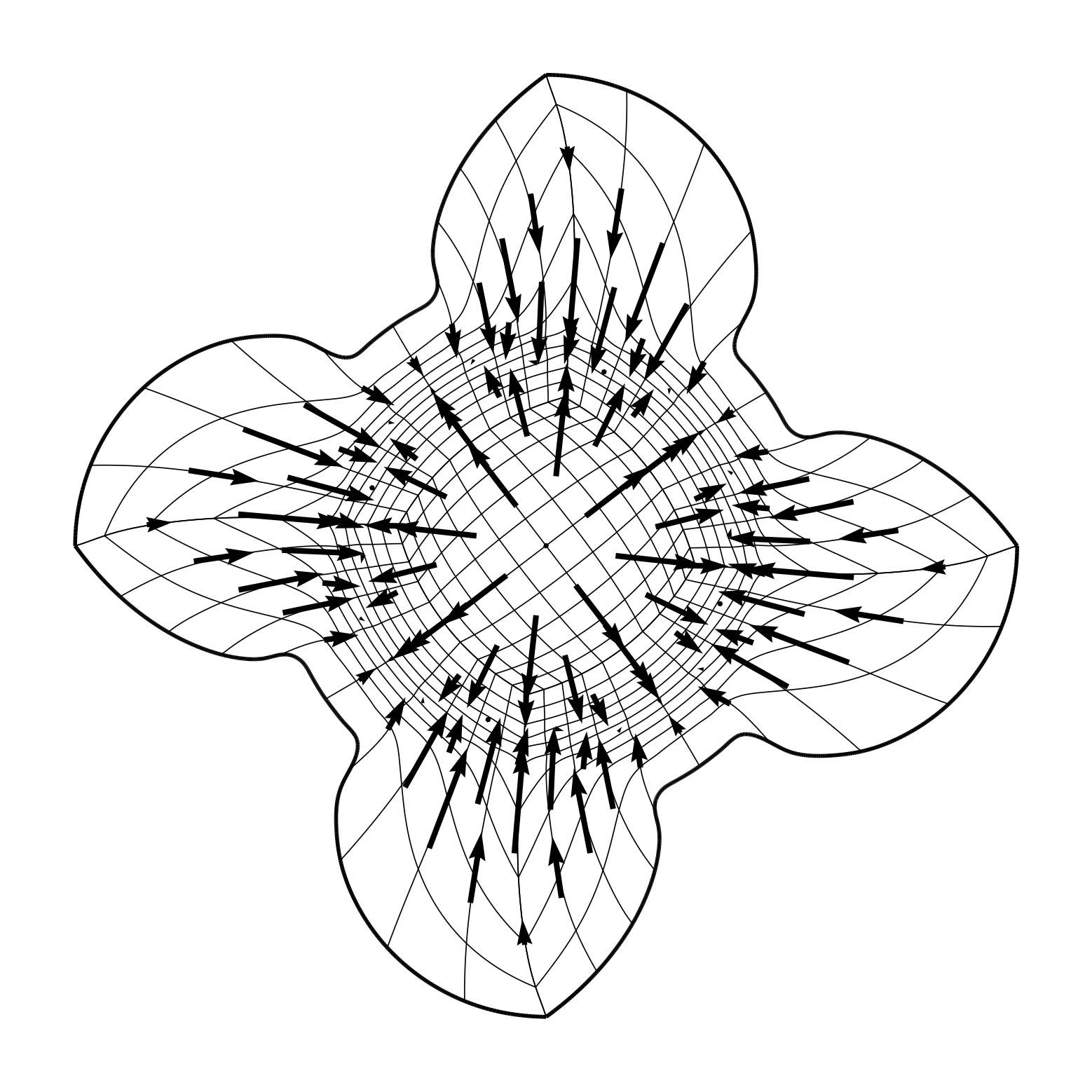}
\caption{Figure showing the cell size contraction in the vicinity of large function values of a ring-shaped function $f \in C^{\infty}(\Om)$. The arrows in the right figure show the movement of select points with respect to the reference parameterisation from Figure \ref{fig:male_screw_benchmark_geom}.}
    \label{fig:male_screw_adaptivity_Id}
\end{figure}

\begin{figure}[H]
\centering
\includegraphics[height=6.8cm]{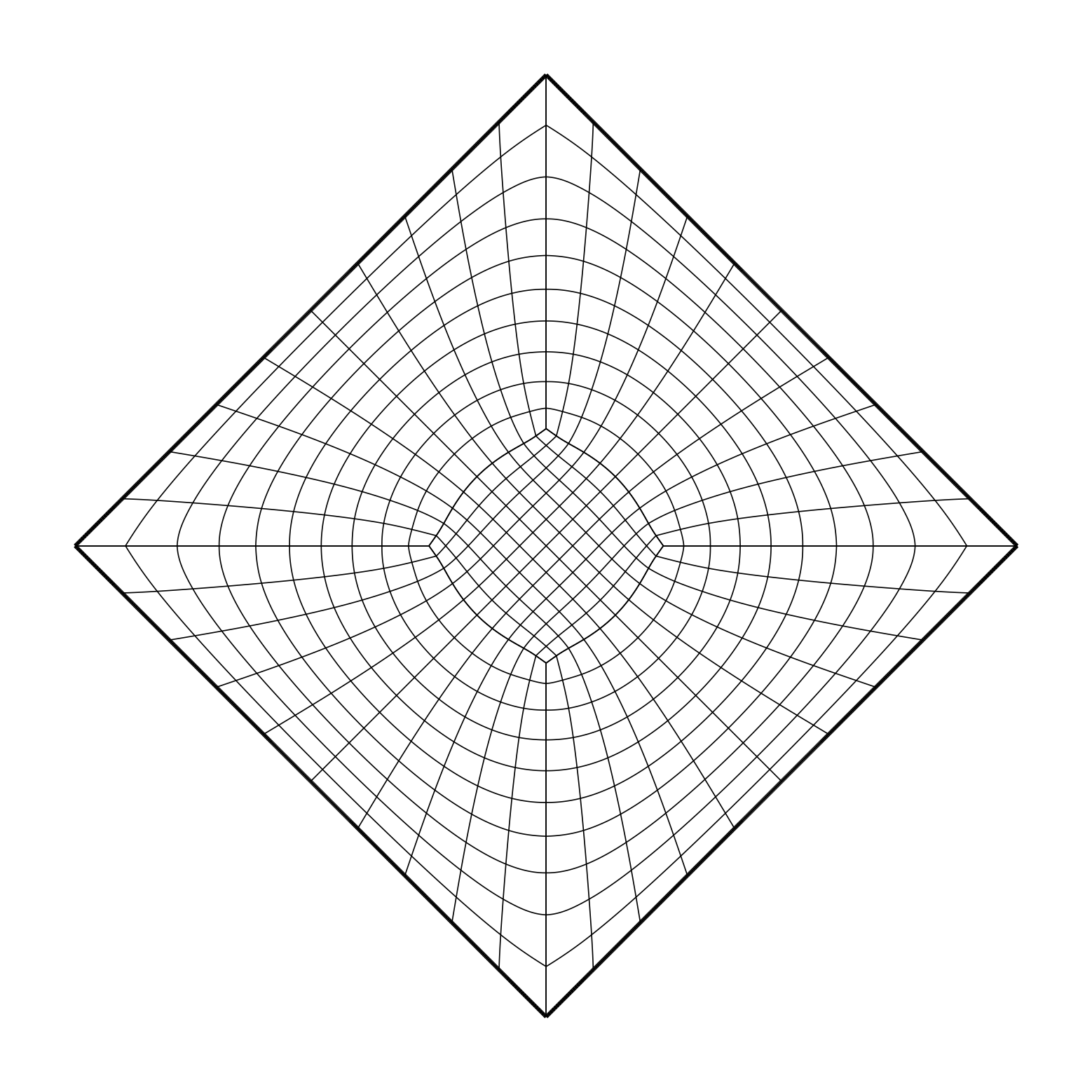} \hspace*{0.1cm}
\includegraphics[height=6.8cm]{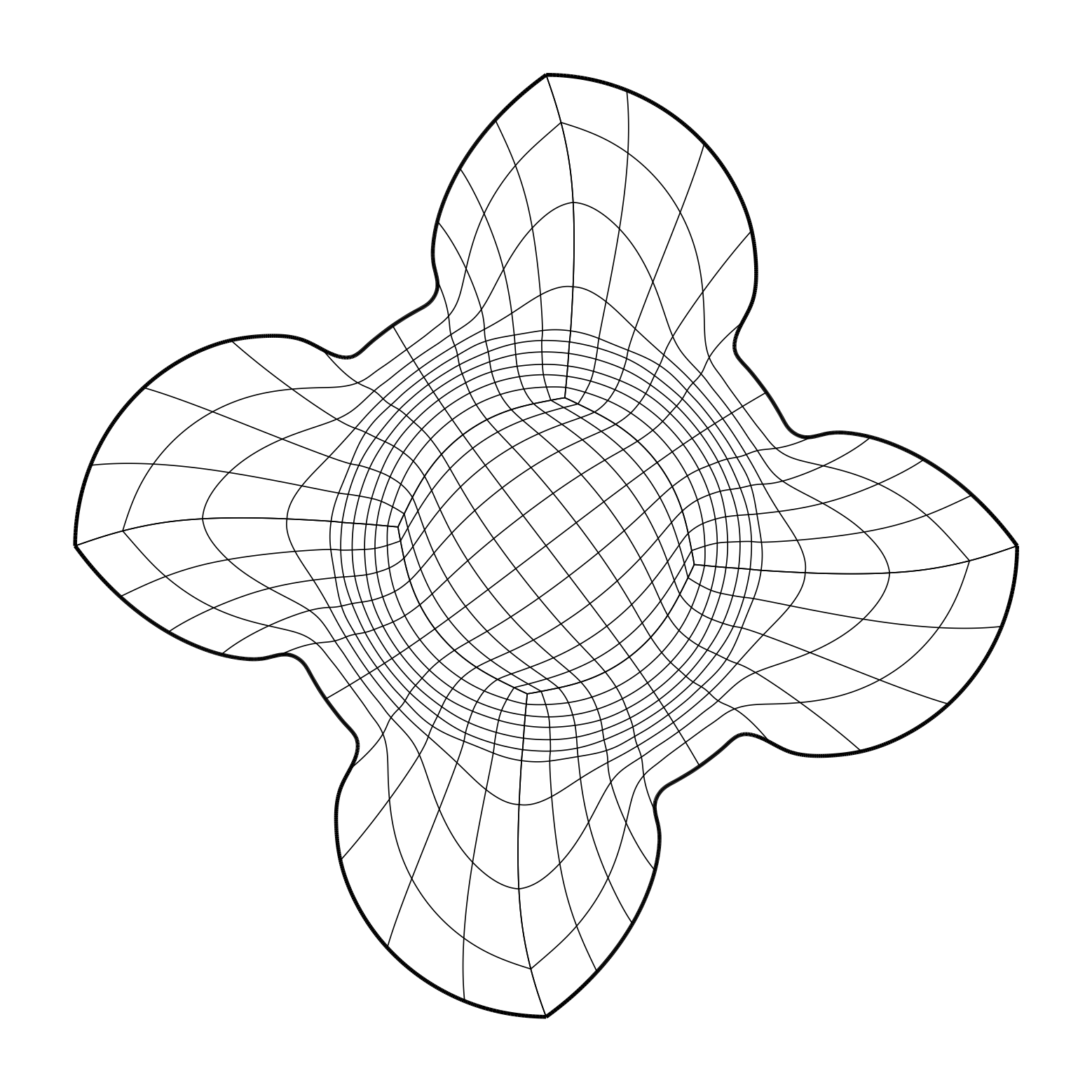}
\caption{Figure showing the parameterisation of $\hOmr = \hOm$ and the result of performing cell size contraction under the coordinate transformation.}
    \label{fig:male_screw_adaptivity_Id_smooth}
\end{figure}

\begin{figure}[H]
\centering
\includegraphics[height=6.8cm]{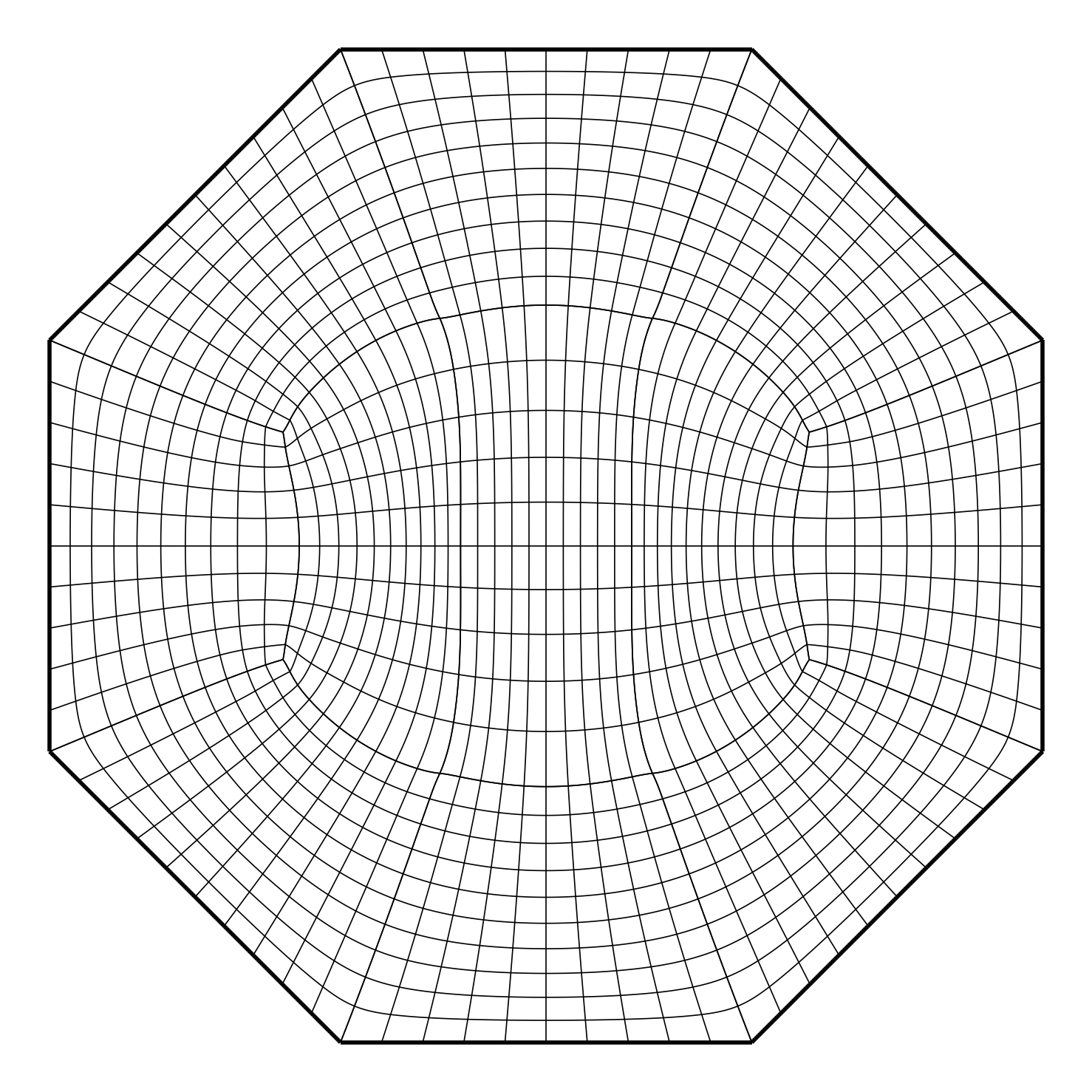} \hspace*{0.1cm}
\includegraphics[height=6.8cm]{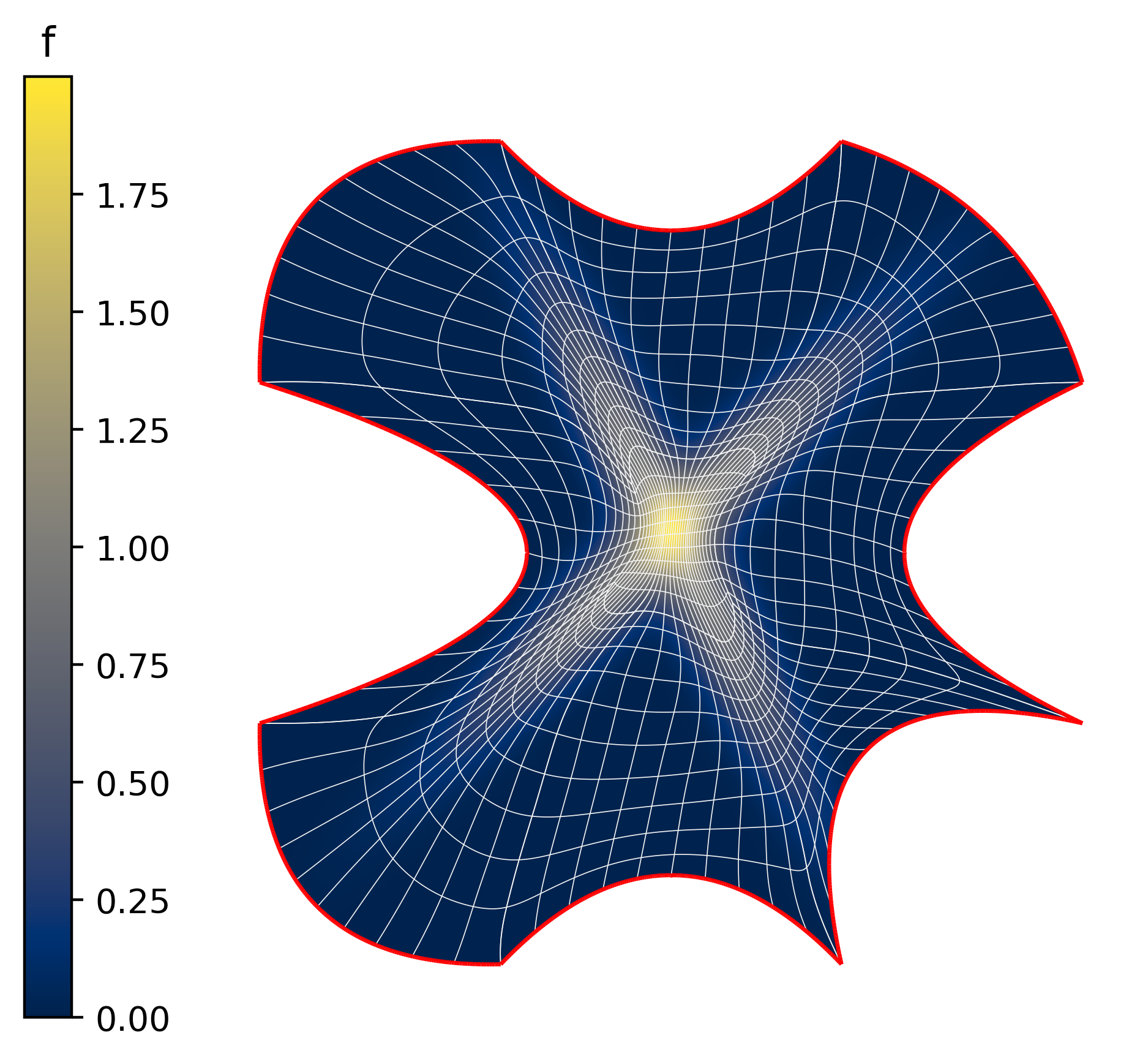}
\caption{A further example of grid adaptation in the vicinity of large function values with interface removal.}
    \label{fig:throwing_star_adaptivity_Id_smooth}
\end{figure}

In various applications it can be desirable to contract the map's isolines in regions where a large value of a function or its gradient is assumed. Given a function $f: \Om \rightarrow \mathbb{R}^{+}$, with $f \in C^{\infty}(\Om)$, the clustering of isolines can be achieved by designing diffusivities that contract the map's cell sizes in regions where $f: \Om \rightarrow \mathbb{R}^{+}$ is large (and vice-versa). The most basic choice is $D^{\bs}(\br, \bx) = \mathcal{I}^{2 \times 2}$ and $D^{\bx}(\bs, \bx) = \sigma(\bx) \mathcal{I}^{2 \times 2}$. As in Section \ref{subsubsect:cell_size_hom}, this choice decouples~\eqref{eq:domain_optimization_coupled} and the first equation can be regarded as the Euler-Lagrange equation of
\begin{align}
    \min \limits_{\br: \Om \rightarrow \hOmr} \int \limits_{\Om} \sigma(\bx) \operatorname{tr} \left( G^{\bx \rightarrow \br} \right) \mathrm{d} \bx \quad \text{s.t.} \quad \br = \left(\mathbf{F}^{\br \rightarrow \bx} \right)^{-1} \text{ on } \partial \Om.
\end{align}
As before, upon exchanging the dependencies $\br(\bx) \rightarrow \bx(\br)$, the isolines will be contracted in regions where $\boldsymbol{\sigma}(\bx)$ is small. To contract cells in the vicinity of large function values of $f: \Omega \rightarrow \mathbb{R}^+$, we design a suitable monitor function $\sigma(\cdot)$. A possible choice is given by \cite[Chapter~9]{liseikin1999grid}:
$$\sigma(\bx) = \frac{1}{\nu_1 f(\bx)^k + \nu_2}, \quad \text{or} \quad \sigma(\bx) = \frac{1}{\nu_1 \| \nabla f(\bx) \|^k + \nu_2} \quad \text{for gradient penalisation}.$$
Here, $\nu_2 > 0$ avoids division by zero in case $f \rightarrow 0$ and the parameters $\nu_1 > 0$ and $k > 0$ tune the degree of penalisation. A numerical scheme is best initialised with the nondegenerate reference solution (i.e., the solution for $\sigma(\bx) = 1$). \\
We are considering the screw geometry with reference parameterisation from Figure \ref{fig:male_screw_benchmark}. Here, we take $\br(\bxi) = \bxi$ and $\hOmr = \hOm$. We would like to contract cells based on the function value of a ring-shaped function $f \in C^{\infty}(\Om)$ using $\nu_1 = 1$, $\nu_2 = 0.01$ and $k = 1$. Figure \ref{fig:male_screw_adaptivity_Id} depicts the result of reparameterising under $D^{\bx}(\bx)$ along with an arrow plot showing the movement of a select number of points with respect to the reference map.
\noindent The figure shows a strong contraction of cells in the vicinity of large function values, clearly demonstrating that the methodology has the desired effect. The cell contraction can be increased by increasing the value of $\nu_1$ or $k$. However, strong penalisation can have unpredictable effects on the cells, in particular close to the patch vertices. \\
In practice, this can be avoided by performing patch interface removal using the techniques of Section \ref{subsect:interface_removal}. We are considering the same example as before while now performing interface removal with stabilisation in $\hOmr$. Figure \ref{fig:male_screw_adaptivity_Id_smooth} shows the images of locally drawn isolines in $\hOmr$ along with the result of performing cell contraction using the same parameters.
Compared to Figure \ref{fig:male_screw_adaptivity_Id}, isolines crossing the patch interfaces exhibit less erratic behaviour and align well with the ring-shaped function. \\
To demonstrate that the proposed technique (in combination with interface removal) is effective when applied to geometries with fewer symmetries than in Figure \ref{fig:male_screw_adaptivity_Id}, we refer to Figure \ref{fig:throwing_star_adaptivity_Id_smooth}.

\subsection{Boundary Orthogonality}
\label{subsect:boundary_orthogonality}
\begin{figure}[h!]
\centering
\begin{subfigure}[b]{0.8 \textwidth}
    \includegraphics[width=0.45 \textwidth]{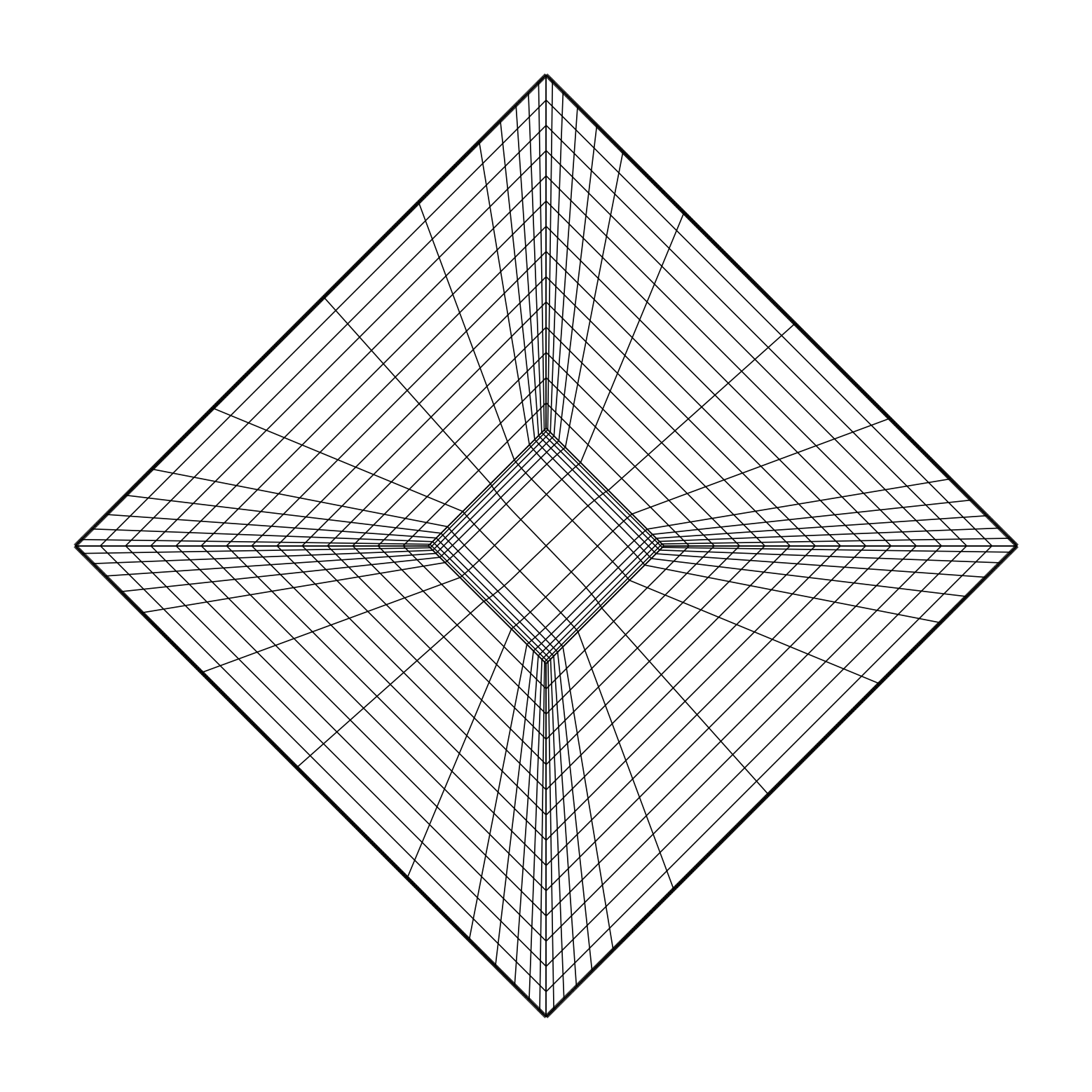} \hfill
    \includegraphics[width=0.45 \textwidth]{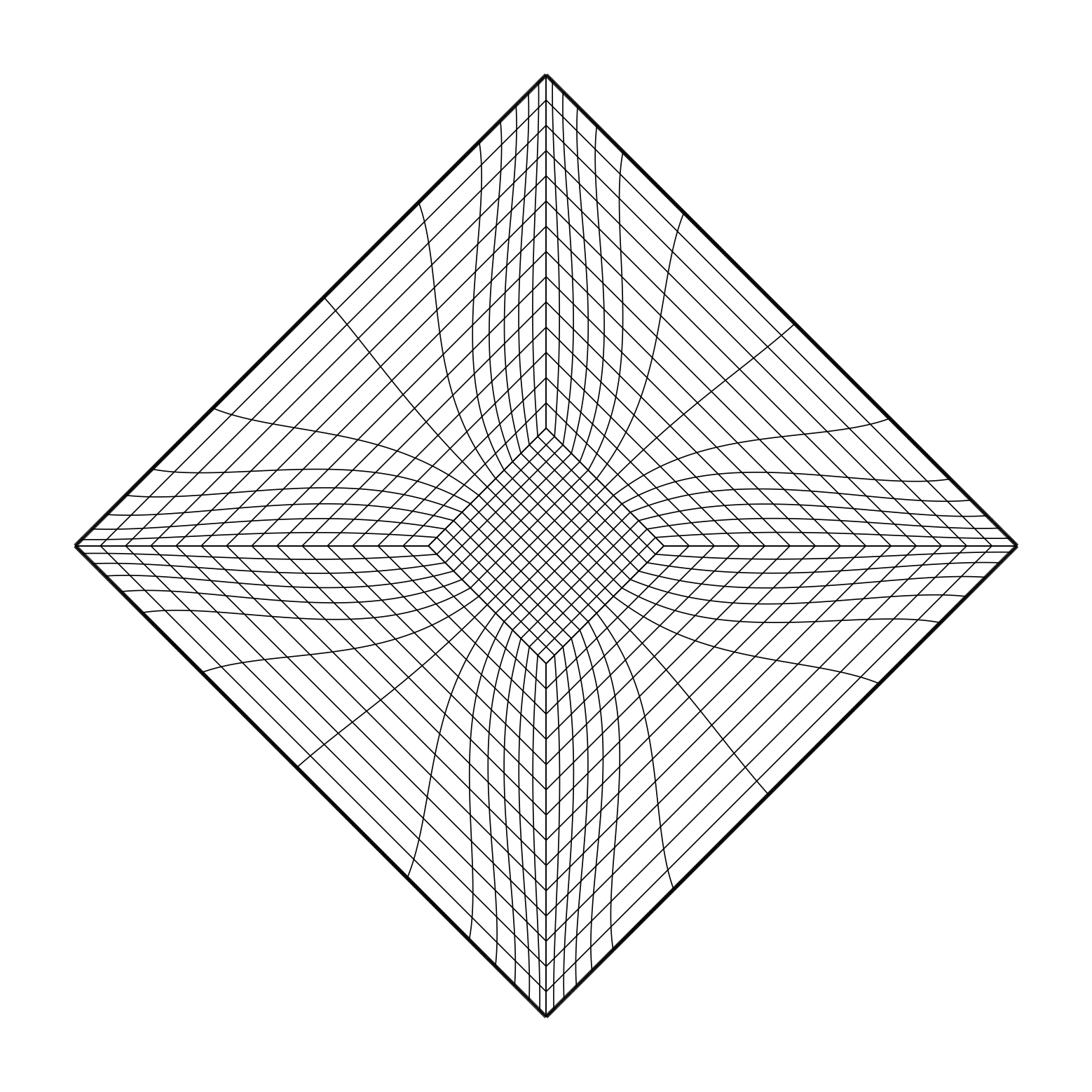}
    \caption{Reparameterisation of $\hOm$ under $\bs: \hOm \rightarrow \hOm$ with the choices $\bmu_1^i(\bmu) = q_i(\mu_1)$ (left) and $\bmu_1^i(\bmu) = t_i(\bmu)$ (right).}
    \label{fig:male_screw_boundary_orth_param}
\end{subfigure}
\\
\begin{subfigure}[b]{0.8 \textwidth}
    \includegraphics[width=0.45 \textwidth]{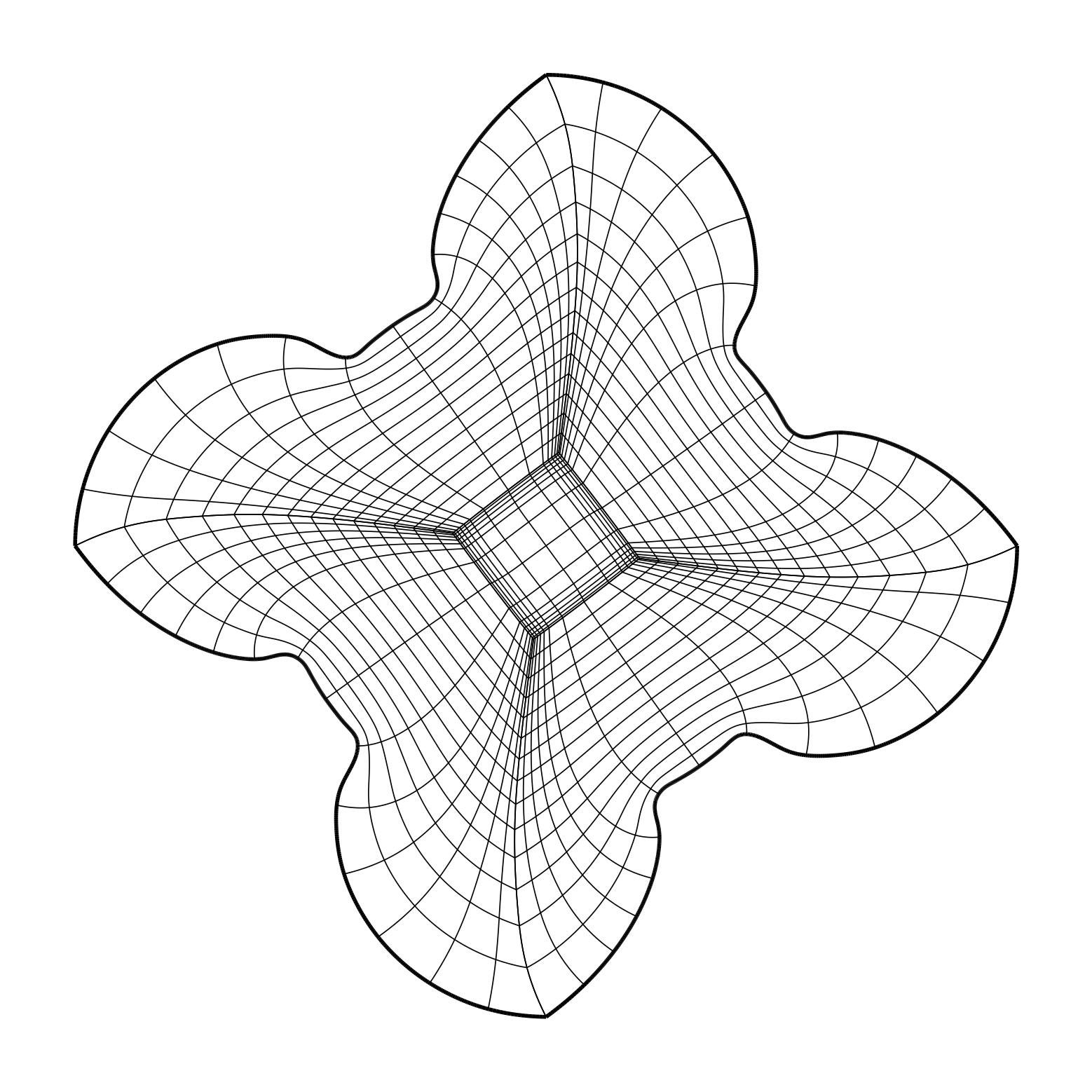} \hfill
    \includegraphics[width=0.45 \textwidth]{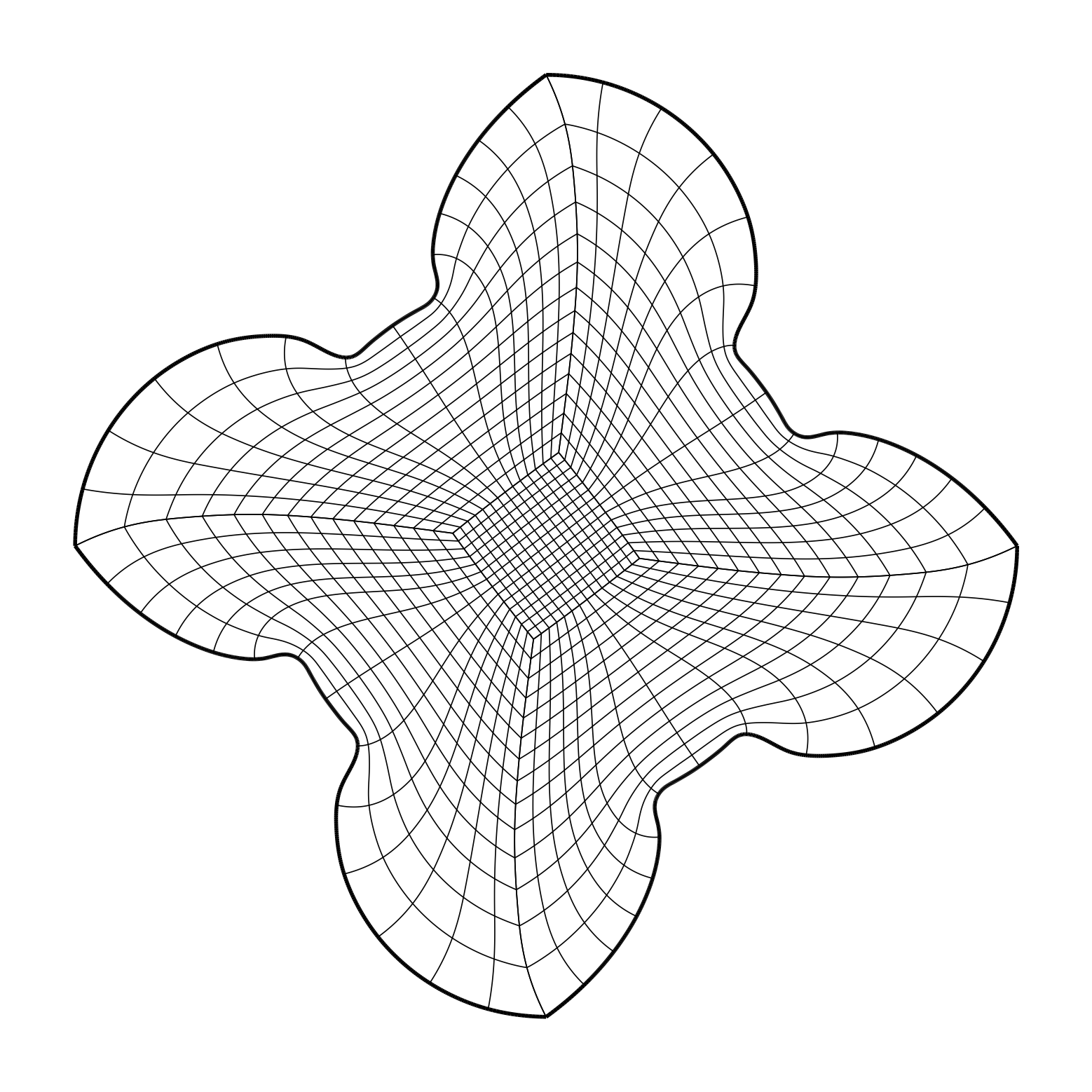}
    \caption{The parameterisations of $\Om$ under the controlmaps from \textbf{(a)}.}
    \label{fig:male_screw_boundary_orth_geom}
\end{subfigure}
\caption{}
\label{fig:male_screw_orth}
\end{figure}

Many applications favour a parameterisation in which locally drawn $\bmu$-isolines intersect the boundary $\partial \Omega$ at a right angle. Unfortunately, it is not possible to simultaneously impose Dirichlet and Neumann data on the inverted elliptic equations. As such, boundary orthogonality has to be enforced through an appropriate coordinate transformation in the parametric domain $\hOmr$, wherein we assume $\hOmr = \hOm$ for convenience. Furthermore, we assume that each $\overline{\hOm_i}$ coincides on $\partial \hOm$ with exactly one of the edges $\overline{L}_k \in \Gamma^B$. Here, we give a multipatch generalisation of the singlepatch method proposed in \cite{hinz2020goal}. Given the reference solution $\bx_h: \hOm \rightarrow \Om$ over the original bilinearly covered parametric domain $\hOm$, we introduce the maps 
$$\bx_h^i(\bmu) := \bx_h \circ \bm^i, \quad \text{for boundary patches } \hOm_i \in \mathcal{Q}^B$$
and we denote $\bx_h^i(\Om^{\square}) := \Omega_i \subset \Omega$. Without loss of generality, we may assume that each $\hOm_i$ is oriented such that $\mu_1$ and $\mu_2$ correspond to the directions tangential and transversal to $\partial \Omega$, respectively. Denoting the eastern, western, southern and northern segments of $\partial \Om^\square$ by $\gamma_e, \gamma_w, \gamma_s$ and $\gamma_n$, respectively, we may furthermore assume the orientation is such that $\gamma_n$ is mapped onto $C_k \subset \partial \Omega$ under $\bx_h^i$. We denote the associated sides of $\Om_i$ under the map $\bx_h^i$ by $\Gamma_e^i, \Gamma_w^i, \Gamma_s^i$ and $\Gamma_n^i = C_k \subset \partial \Omega$. We are seeking a function $f_i: \Omega_i \rightarrow \mathbb{R}$ that satisfies homogeneous Neumann boundary conditions on $\Gamma_n^i$. For this, we solve
$$\Delta f_i = 0 \text{ in } \Om_i, \quad \text{s.t.} \quad f_i(\bx(\bmu)) = \mu_1 \text{ on } \overline{\Gamma_e^i} \cup \overline{\Gamma_w^i} \cup \overline{\Gamma_s^i} \quad \text{and} \quad \frac{\partial f_i}{\partial \mathbf{n}} = 0 \text{ on } \Gamma_n^i,$$
where $\mathbf{n}$ denotes the unit outward normal vector on $\partial \Omega_i$. A pullback leads to
$$\Delta_{\bx_h^i} f_i = 0 \text{ in } \Om^\square, \quad \text{s.t.} \quad f_i(\bmu) = \mu_1 \text{ on } \overline{\gamma_e} \cup \overline{\gamma_w} \cup \overline{\gamma_s} \quad \text{and} \quad \frac{\partial f_i}{\partial \mathbf{n}} = 0 \text{ on } \gamma_n,$$
whose discretisation imposes the Neumann data through partial integration in the usual way. The restriction $q_i(\mu_1) := f_i \vert_{\gamma_n}$ will be a monotone function over $\mu_1 \in [0, 1]$ thanks to the imposed Dirichlet data and the maximum principle. Given a diffeomorphism $\bnu^i: \Om^{\square} \rightarrow \Om^{\square}$ that satisfies $\bnu^i \vert_{\gamma_n} = q_i$ and $\partial_{\mu_2} \bnu^i_1 = 0$ on $\gamma_n$, a map $\tilde{\bx}^i: \Om^{\square} \rightarrow \Omega_i$ that maps inversely harmonically into $\Omega^{\square}$ with the coordinate system induced by $\bnu^i$, will map local $\mu_1$ isolines onto isolines in $\Omega_i$ that intersect $\overline{C_k} = \partial \Omega \cap \overline{\Omega}_i$ at a right angle, thanks to the Neumann data we imposed on $\Gamma_n^i$. Two possible choices are given by
\begin{align}
\label{eq:boundary_orth_mu1_functions}
    \bnu^i_1(\bmu) = q_i(\mu_1) \quad \text{and} \quad \bnu^i_1(\bmu) = t_i(\bmu) \quad \text{with} \quad t_i(\bmu) := (1 + 2 \mu_2) \, (1 - \mu_2)^2 \, \mu_1 + (3 - 2 \mu_2) \, \mu_2^2 \, q_i(\mu_1),
\end{align}
while $\bnu^i_2(\bmu) = \mu_2$. Here, the former maps straight $\mu_1$ isolines onto straight $q_i(\mu_1)$ isolines in $\Om^{\square}$ while the latter maps the same isolines onto curves that start at $\bmu = (\mu_1, 0)$ and end in $\bmu = (q_i(\mu_1), 1)$ while intersecting $\gamma_n$ at a right angle. Note that the latter furthermore satisfies $\bnu^i(\bmu) = \bmu$ on $\gamma_s$. \\
With this choice of $\bnu^i$, a map $\tilde{\bx}: \hOm \rightarrow \Om$ that maps inversely harmonically into $\hOm$ with the coordinate system induced by the controlmap $\bs: \hOm \rightarrow \hOm$ that satisfies 
$$\bs(\bmu) = \bm^i \circ \bnu^i, \quad \text{or equivalently} \quad \bs(\bxi) = \bm^i \circ \bnu^i \circ (\bm^i)^{-1} \quad \text{on} \quad \hOm_i \in \mathcal{Q}^B,$$
will now map the images of local $\mu_1$-isolines under $\bm^i$ onto isolines in $\Om_i$ that intersect $\partial \Omega$ at a right angle. The controlmap $\bs: \hOm \rightarrow \hOm$ that leads to boundary orthogonality is hence known for boundary patches $\hOm_i \in \mathcal{Q}^B$. For the choice $\bnu^i_1(\bmu) = t_i(\bmu)$, the controlmap $\bs: \hOm \rightarrow \hOm$ can be taken as the identity on patches $\hOm_k \notin \mathcal{Q}^B$. \\
For the choice $\bnu^i_1 = q_i(\mu_1)$, on the other hand, the partially-known controlmap induces a reparameterisation of the interior facets $\gamma_{ij} \in \Gamma^I$ with $\hOm_i \in \mathcal{Q}^B$ or $\hOm_j \in \mathcal{Q}^B$. As such, the original bilinear parameterisations of patches $\hOm_k \notin \mathcal{Q}^B$ will no longer be conforming to the images of boundary patches under $\bs(\bxi)$. In this case, we may require the controlmap to satisfy $\bs(\bxi) \vert_{\gamma_{ij}} = \bxi$ for facets $\gamma_{ij} \in \Gamma^I$ that are not associated with one of the boundary patches. With that, the map $\bs: \hOm \rightarrow \hOm$ is now known on the boundary patches and on all interior facets. The interior of the remaining patches $\hOm_k \notin \mathcal{Q}^B$ can now be parameterised from the curves $\bs \vert_{\gamma_{ij}}, \, \gamma_{ij} \in \Gamma^I$ one-by-one using, for instance, the bilinearly-blended Coons' patch approach. The result is a controlmap $\bs: \hOm \rightarrow \hOm$ that leads to boundary orthogonality and is furthermore conforming across all interior facets. \\
The controlmap can be projected onto $\mathcal{V}_h^2$, where we note that for $\bnu^i_1 = t_i(\bmu)$, $\mathcal{V}_h$ has to be patchwise bicubic in order for the projection to be exact. \\
We are again considering the screw geometry whose reference controlmap and parameterisation are depicted in Figure \ref{fig:male_screw_cell_k0}. As a measure for the degree of boundary orthogonalisation, we utilise
\begin{align}
    L_{\perp}^2(\bx) := \sum \limits_{\hOm_i \in \mathcal{Q}^B} \, \int \limits_{\partial \hOm \cap \overline{\hOm}_i} \left(\widehat{\partial}_{\bmu^{\perp}} \bx \cdot \widehat{\partial}_{\bmu^\parallel} \bx \right)^2 \mathrm{d} \Gamma,
\end{align}
wherein $\widehat{\partial}_{\mu^{\parallel}}( \, \cdot \, )$ and $\widehat{\partial}_{\mu^{\perp}}( \, \cdot \, )$ denote the normalised tangential and transverse derivatives with respect to $\bmu$ on $\hOm_i \in \mathcal{Q}^B$, respectively. \\
Figure \ref{fig:male_screw_boundary_orth_geom} shows the reparameterisations of $\hOm$ under $\bs: \hOm \rightarrow \hOm$ for the choices $\bnu^i_1(\bmu) = q_i(\mu_1)$ and $\bnu^i_1(\bmu) = t_i(\bmu)$, respectively, while Figure \ref{fig:male_screw_boundary_orth_geom} shows the associated parameterisations of $\Om$ after recomputation under $\bs: \hOm \rightarrow \hOm$.
With
$$\frac{L_{\perp}(\bx^{\bs})}{L_{\perp}(\bx^{\bxi})} = 0.172 \quad \text{for} \quad \bmu_1^i(\bmu) = q_i(\mu_1) \quad \text{and} \quad \frac{L_{\perp}(\bx^{\bs})}{L_{\perp}(\bx^{\bxi})} = 0.169 \quad \text{for} \quad \bmu_1^i(\bmu) = t_i(\bmu),$$
both choices are similarly effective, wherein the small discrepancy is explained by differing truncation errors. Figure \ref{fig:male_screw_boundary_orth_geom}, left, reveals that the choice $\bmu_1^i(\bmu) = q_i(\mu_1)$ leads to a strong clustering of cells close to the patch interfaces, which may not be desirable. This is not the case for \ref{fig:male_screw_boundary_orth_geom}, right. Meanwhile, the former maps straight isolines in $\Om^{\square}$ onto straight isolines on the $\hOm_i \in \mathcal{Q}^B$, a property that can be exploited when combining boundary orthogonalisation with the creation of a boundary layer (see Section \ref{subsect:boundary_layers}).

\subsection{Boundary layers}
\label{subsect:boundary_layers}
\begin{figure}[H]
\centering
\includegraphics[height=6.8cm]{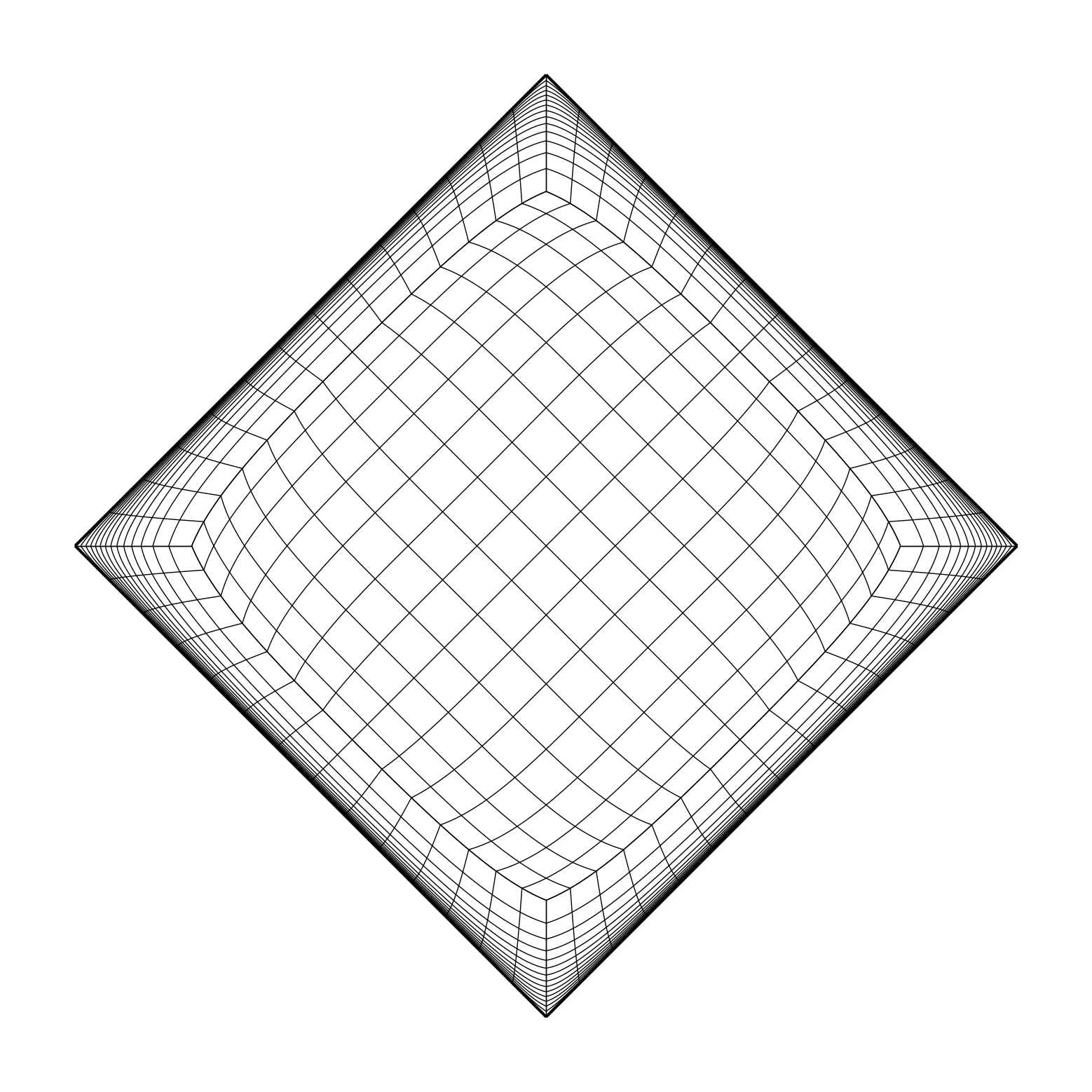} \hspace*{0.1cm}
\includegraphics[height=6.8cm]{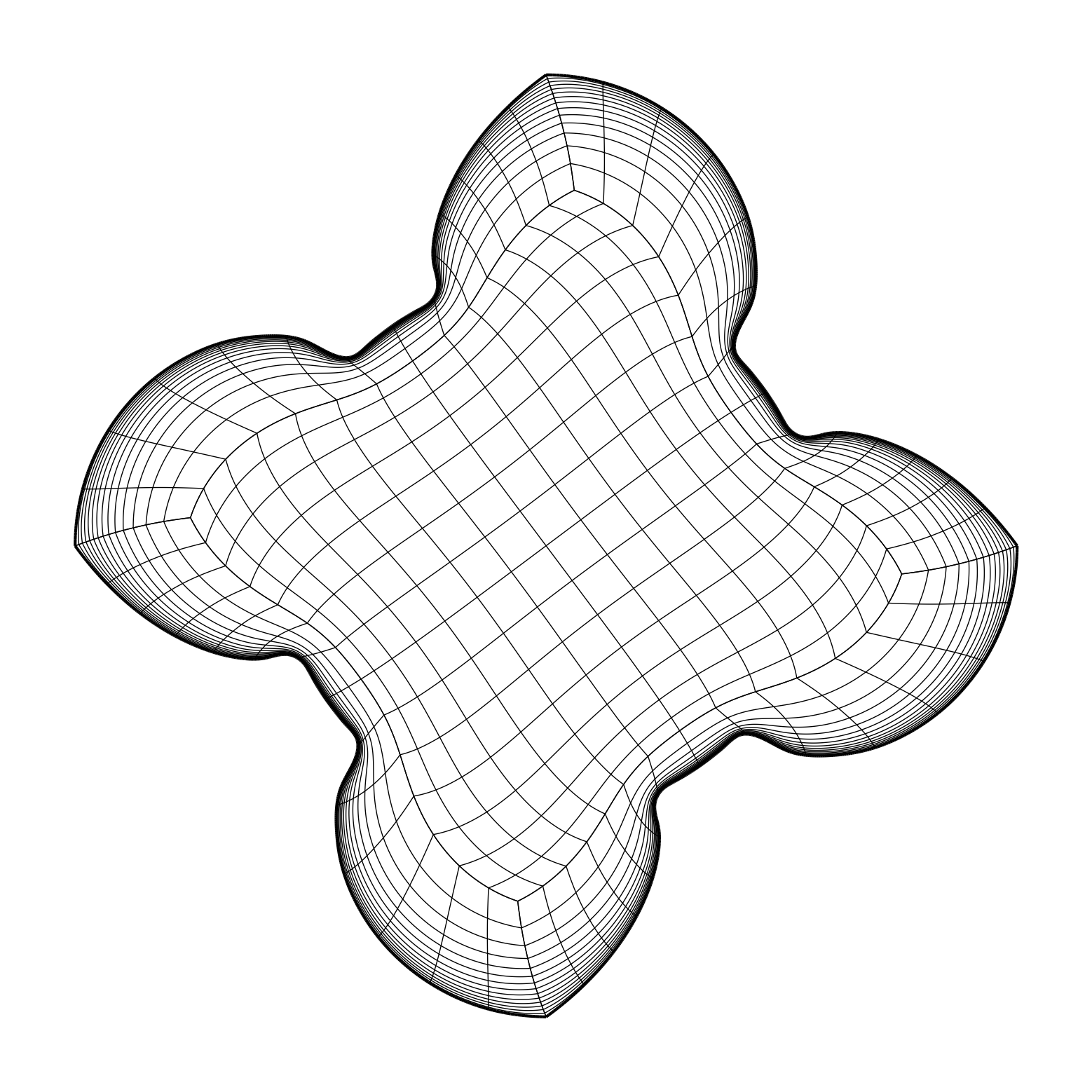}
\caption{Example showing the creation of a boundary layer using~\eqref{eq:boundary_layer_diffusivity} with $\mu = 30$, $k=2$ and $\nu = 0.005$.}
\label{fig:male_screw_boundary_layer_vanilla}
\end{figure}

\begin{figure}[H]
\centering
\includegraphics[height=6.8cm]{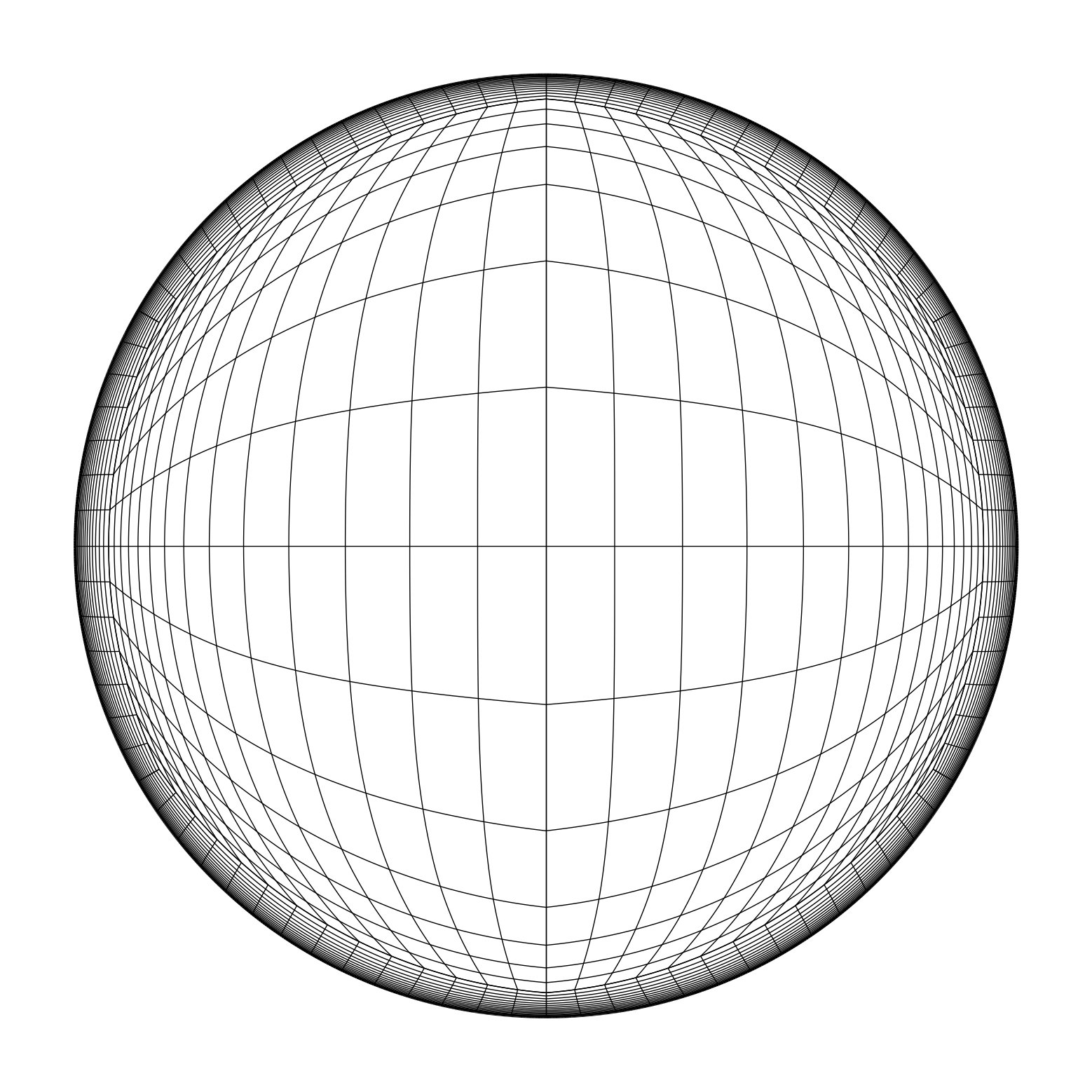} \hspace*{0.1cm}
\includegraphics[height=6.8cm]{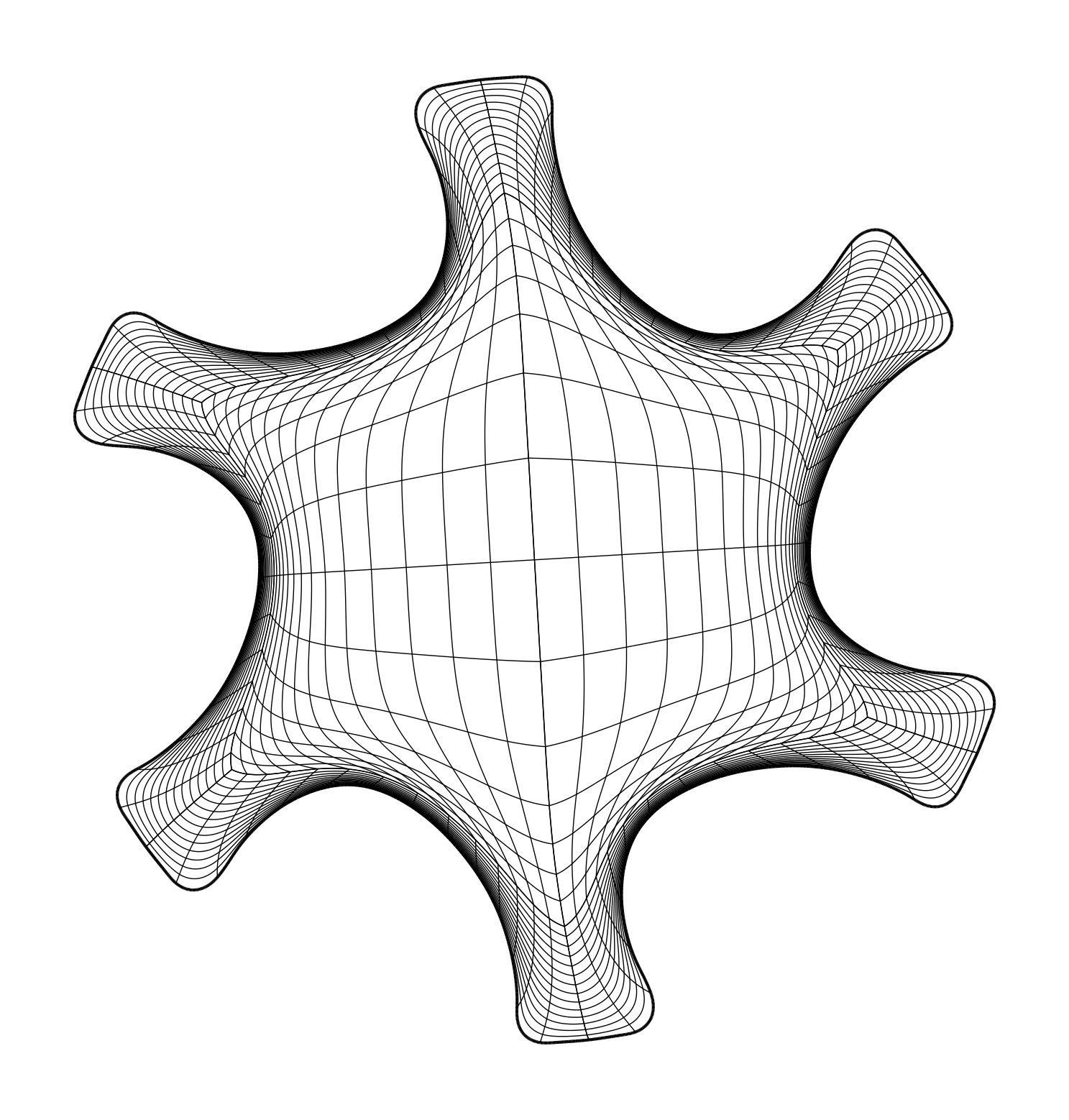}
\caption{Example showing the creation of a boundary layer using~\eqref{eq:boundary_layer_diffusivity} with $\mu = 30$, $k=2$ and $\nu = 0.005$.}
\label{fig:female_screw_boundary_layer_vanilla}
\end{figure}

Many applications in computational fluid dynamics deal with PDE-problems whose solutions are known to create a steep gradient in the vicinity of the boundary $\partial \Om$. To capture important features of the solution, such applications favour parameterisations with a dense clustering of cells in the vicinity of $\partial \Omega$. This can be achieved by introducing a controlmap $\bs: \hOm \rightarrow \hOm$ that clusters locally-drawn $\bmu$-isolines close to $\partial \hOm$, while potentially sacrificing some cell density in the interior. We assume that $\hOmr = \hOm$ is a regular $2n$-sided polygon with radius one, centered at $\bxi = (0, 0)^T$. As before, we assume that $\mu_1 = \text{const}$ refers to the local $\bmu$-direction transverse to $\partial \hOm$ for boundary patches. A convenient way to create a boundary layer is utilising a diffusivity of the form
\begin{align}
\label{eq:boundary_layer_diffusivity}
    D^{\bs}(\bxi) = (1 - \exp{(-\mu \|\bxi \|^2)}) \, \| \bxi \|^k \, \left( \hat{\bxi} \otimes \hat{\bxi} \right) + \nu \mathcal{I}^{2 \times 2},
\end{align}
for some $\mu \gg 1$, $\nu \ll 1$ and $k > 0$. Here, $\hat{\bxi} := \bxi / \| \bxi \|$ while the prefactor $(1 - \exp{(-\mu \|\bxi \|^2)})$ in~\eqref{eq:boundary_layer_diffusivity} avoids the singularity in the origin. Taking $\nu$ small urges the controlmap $\bs: \hOm \rightarrow \hOm$ to map points $\bxi \in \hOm$ in radially outward direction with a tight clustering close to $\partial \hOm$. If the locally-drawn $\mu_1$ isolines of boundary patches $\hOm_i \in \mathcal{Q}^B$ align with straight rays drawn from the origin in radially outward direction as in Figure \ref{fig:male_screw_cell_k0}, the diffusivity from~\eqref{eq:boundary_layer_diffusivity} will furthermore largely preserve the intersection angle of the $\mu_1$ isolines with the boundary $\partial \Omega$ (w.r.t. to the choice $D^{\bs} = \mathcal{I}^{2 \times 2}$). Since the controlmap creation is an a priori step, the methodology is compatible with all operators from Section \ref{sect:numerical_schemes}. \\
We are considering the screw geometry with the reference parameterisation from Figure \ref{fig:male_screw_cell_k0}. Figure \ref{fig:male_screw_boundary_layer_vanilla} shows the reparameterisation using~\eqref{eq:boundary_layer_diffusivity} with $\mu = 30$, $k = 2$ and $\nu = 0.005$.
The figure shows a strong clustering of transverse isolines close to the boundary. The clustering intensity can be increased by increasing the value of $k$. \\
As a second example, we are considering the female screw geometry with reference controlmap and parameterisation shown in Figure \ref{fig:female_screw_reference_domain}. While $\hOmr \neq \hOm$ is now given by the unit disc, the boundary patches are oriented such that the $\mu_1$ transverse isolines align with straight rays that intersect the origin, as required. As such, we replace $\bxi \rightarrow \br$ in~\eqref{eq:boundary_layer_diffusivity} and expect the diffusivity to largely preserve the intersection angle of transverse isolines with $\partial \hOmr$ and $\partial \Om$. Figure \ref{fig:female_screw_boundary_layer_vanilla} shows the result of boundary layer creation under~\eqref{eq:boundary_layer_diffusivity} using the same parameters as in Figure \ref{fig:male_screw_boundary_layer_vanilla}.
In both examples, the layer density is intensified / reduced close to the inward-facing and protruded parts of $\partial \Omega$, respectively. This effect may be counteracted by introducing a position-dependent clustering parameter, i.e., $k \rightarrow k(\br)$. However, this is beyond the scope of this paper. \\
As a final example, we are combining boundary layer creation with boundary orthogonality. For this purpose, we assume that we are in the possession of a controlmap $\bs: \hOm \rightarrow \hOm$ that orthogonalises transverse $\mu_1$ isolines for boundary patches $\hOm_i \in \mathcal{Q}^B$. Assuming again that the orientation is chosen such that the transverse direction is given by $\mu_1 = \text{const}$ for the $\hOm_i \in \mathcal{Q}^B$, we restrict ourselves to the choice $\bnu_1^i(\bmu) = q_i(\mu_1)$ (cf.~\eqref{eq:boundary_orth_mu1_functions}) such that $\bs: \hOm \rightarrow \hOm$ maps straight $\mu_1$ isolines onto straight $\mu_1$ isolines for boundary patches, as in Figure \ref{fig:male_screw_boundary_orth_param} (left). For the purpose of boundary layer creation, we may compose the $\bnu^i(\bmu) = (q_i(\mu_1), \mu_2)^T$ with a function of the form $\blambda^i(\bmu) = (\mu_1, f_i(\mu_1, \mu_2))^T$, where $f_i(\mu_1, 0) = 0$, $f_i(\mu_1, 1) = 1$ and $\partial_{\mu_2} f_(\mu_1, \mu_2) > 0$. We note that the map $\blambda^i \circ \bnu^i: \Om^{\square} \rightarrow \Om^{\square}$ still satisfies $\partial_{\mu_2} (\blambda^i \circ \bnu^i )_1 = 0$ on $\gamma_n$ such that boundary orthogonality is preserved. The purpose of $f_i$ is to create a boundary layer near $\mu_2 = 1$, i.e., $\partial_{\mu_2} f_i(\mu_1, \mu_2)$ in $\mu_2 = 1$ should be small. A function that satisfies aforementioned requirements is given by
\begin{align}
    f_i(\mu_1, \mu_2) = \frac{1 - \exp{(-d_i(\mu_1) \mu_2)}}{1 - \exp{(-d_i(\mu_1))}}, \quad \text{ for some } \quad d_i(\mu_1) > 0,
\end{align}
where larger values of $d_i(\mu_1)$ create a (locally) stronger boundary layer. Note that we have
\begin{align}
    \left. \partial_{\mu_2} f_i(\mu_1, \mu_2) \right \vert_{\gamma_n} := f^{\prime}_{i, n}(\mu_1) = \frac{d_i(\mu_1) \exp{(-d_i(\mu_1))}}{1 - \exp{(-d_i(\mu_1))}}.
\end{align}
Thanks to boundary orthogonalisation, the local tangent bundle of $\bx^{\bs}$ is diagonal in the basis that is spanned by the unit tangent and normal vectors on $\partial \Omega$, i.e,
\begin{align}
\label{eq:boundary_layer_local_tangent_bundle_diagonal}
    \partial_{\bmu} \bx^{\bs} = a_i \mathbf{t} \otimes \mathbf{t} + b_i \mathbf{n} \otimes \mathbf{n}, \quad \text{on} \quad L_{k^i} := \overline{\hOm_i} \cap \partial \Omega \quad \text{for } \hOm_i \in \mathcal{Q}^B,
\end{align}
where $a_i$ follows from the boundary correspondence $\mathbf{F}: \partial \hOm \rightarrow \partial \Om$ while $b_i$ is a property of the parameterisation of the interior. We would like to create a new controlmap $\tilde{\bs}: \hOm \rightarrow \hOm$ whose associated map $\bx^{\tilde{\bs}}: \hOm \rightarrow \Om$ satisfies
\begin{align}
    \partial_{\bmu} \bx^{\tilde{\bs}} = a_i \mathbf{t} \otimes \mathbf{t} + k \mathbf{n} \otimes \mathbf{n} \quad \text{on} \quad L_{k^i} := \overline{\hOm_i} \cap \partial \Omega \quad \text{for } \hOm_i \in \mathcal{Q}^B,
\end{align}
where $k > 0$ is a user-specified parameter. Requiring $\tilde{\bs}: \hOm \rightarrow \hOm$ to be of the form
\begin{align}
\label{eq:orth_boundary_layer_controlmap_boundary_patches}
    \tilde{\bs}(\bxi) = \bm^i \circ \blambda^i \circ \bnu^i  \circ (\bm^i)^{-1} \quad \text{on } \hOm_i \in \mathcal{Q}^B,
\end{align}
it is clear that the recomputed map will satisfy
\begin{align}
    \left. \partial_{\bmu} \bx^{\tilde{\bs}}(\bmu) \right \vert_{\gamma_n} = a_i(\mu_1) \mathbf{t} \otimes \mathbf{t} + b_i(\mu_1) f_{i, n}^{\prime}(\mu_1) \mathbf{n} \otimes \mathbf{n} \quad \text{on} \quad L_{k^i} \subset \partial \hOm \quad \text{for } \hOm_i \in \mathcal{Q}^B.
\end{align}
Given that each $f_{i, n}^{\prime} = f_{i, n}^{\prime}(d_i)$, we may find the $d_i(\mu_1)$ by minimising the nonlinear cost function
\begin{align}
\label{eq:constant_transverse_deriv_opt}
    \sum \limits_{\{i \, \vert \, \hOm_i \in \mathcal{Q}^B \}} \, \, \int \limits_{L_{k^i}} (b_i f^{\prime}_{i, n} -k )^2 \mathrm{d} \Gamma \rightarrow \min \limits_{d},
\end{align}
where $b_i$ and $f^{\prime}_{i, n}$ are now taken as functions over $\bxi \in L_{k^i} \subset \partial \hOm$ while $d: \partial \hOm \rightarrow \mathbb{R}$ satisfies
\begin{align}
    d(\bxi) = d_i \circ (\mathbf{m}^i)^{-1} \quad \text{on } L_{k^i} \subset \partial \hOm.
\end{align}
A discretisation then constructs $d(\bxi)$ from the linear span of the $\phi_h \in \mathcal{V}_h$ that are nonvanishing on $\partial \hOm$ and finds the minimum over this subspace of $C^0(\partial \hOm)$ using Newton's method. Upon completion, the function $d(\bxi)$ is known on $\partial \hOm$ and the restriction to $L_{k^i} \subset \partial \hOm$ can be expressed as a function over $\gamma_n \subset \partial \Om^{\square}$ via the $(\mathbf{m}^i)^{-1}$. The canonical extension from $\gamma_n$ into $\Om^{\square}$ then defines $\blambda^i(\bmu) = (\mu_1, f_i(\mu_1, \mu_2))^T$ and we define $\tilde{\bs}: \hOm \rightarrow \hOm$ by~\eqref{eq:orth_boundary_layer_controlmap_boundary_patches} for boundary patches while requiring $\tilde{\bs}(\bxi) = \bs(\bxi)$ for interior patches. Since $d(\bxi) \in C^0(\partial \hOm)$, the controlmap defined in this way will be conforming across all interior interfaces. Finally, to relieve the computational burden, the controlmap $\tilde{\bs}: \hOm \rightarrow \hOm$ is expressed in $\mathcal{V}_h^2$ via an $L^2(\hOm)$ projection. \\

\begin{figure}[H]
\centering
\includegraphics[height=6.8cm]{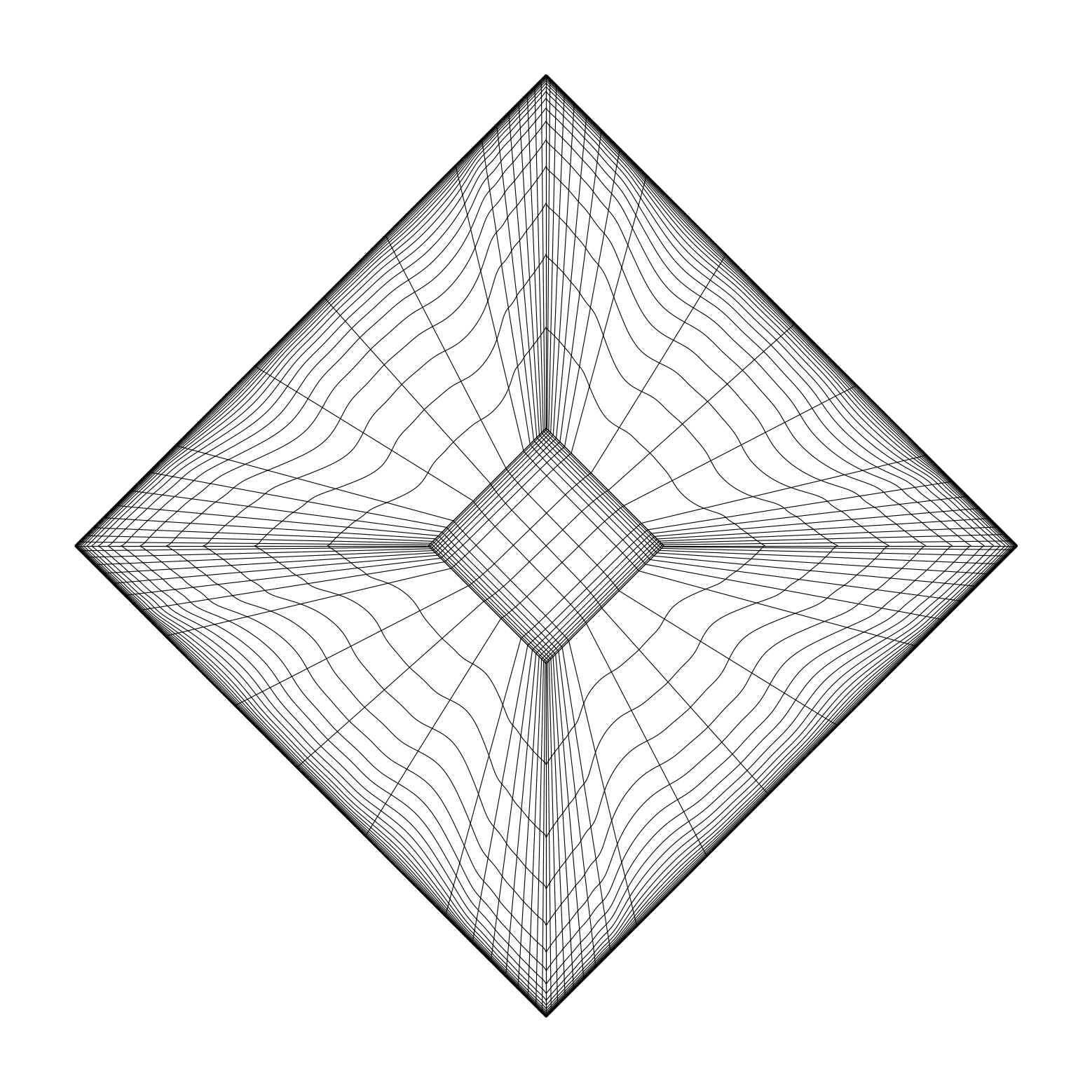} \hspace*{0.1cm}
\includegraphics[height=6.8cm]{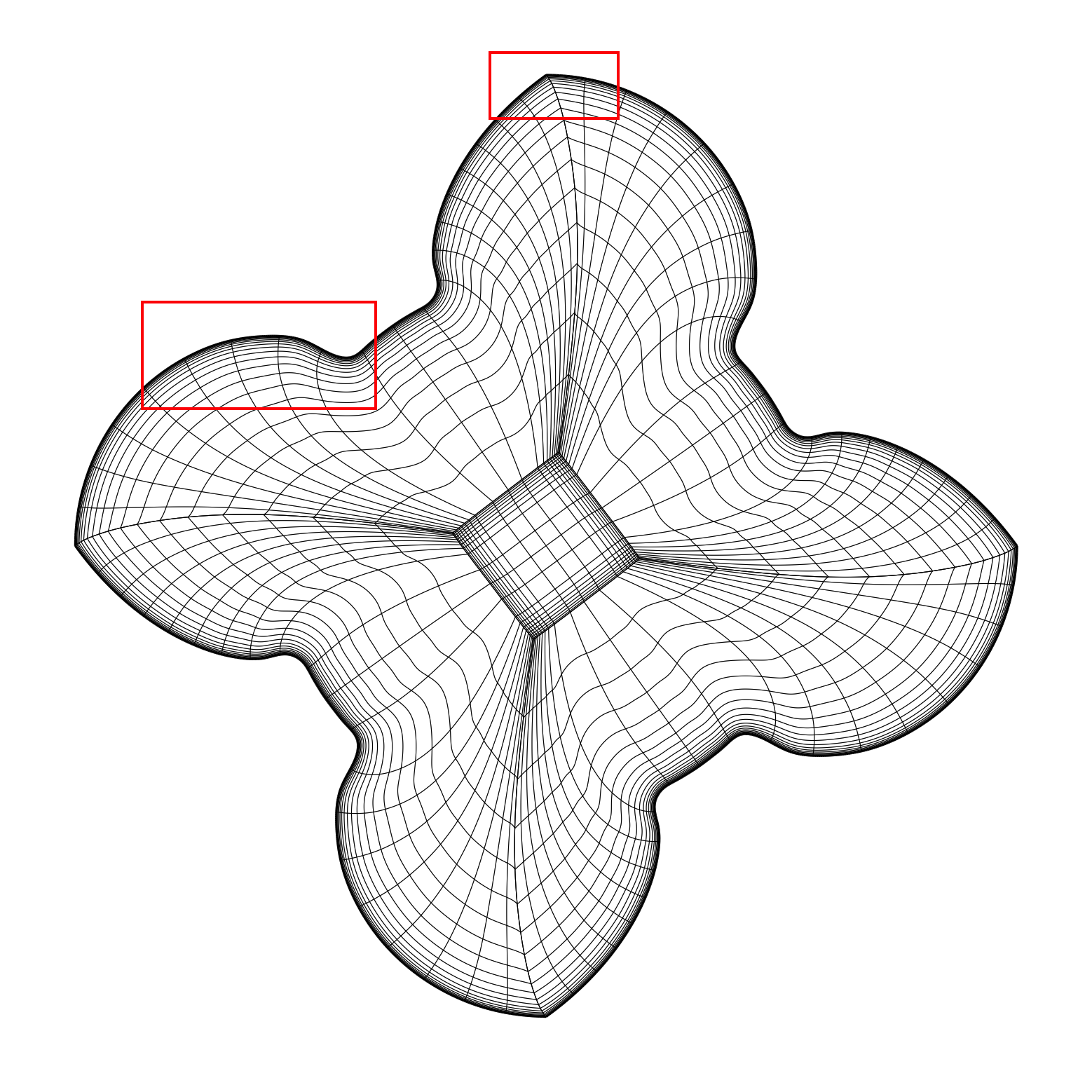}
\caption{Controlmap and the associated parameterisation of the male screw geometry orthogonalised by the boundaries along with a boundary layer that has been created using a sequence of algebraic operations on the boundary patches. Figure \ref{fig:male_screw_boundary_layer_algebraic_zoom} shows a zoom-in on the segments that have been highlighted in the right figure.}
\label{fig:male_screw_boundary_layer_algebraic}
\end{figure}

\begin{figure}[H]
\centering
\includegraphics[height=4cm]{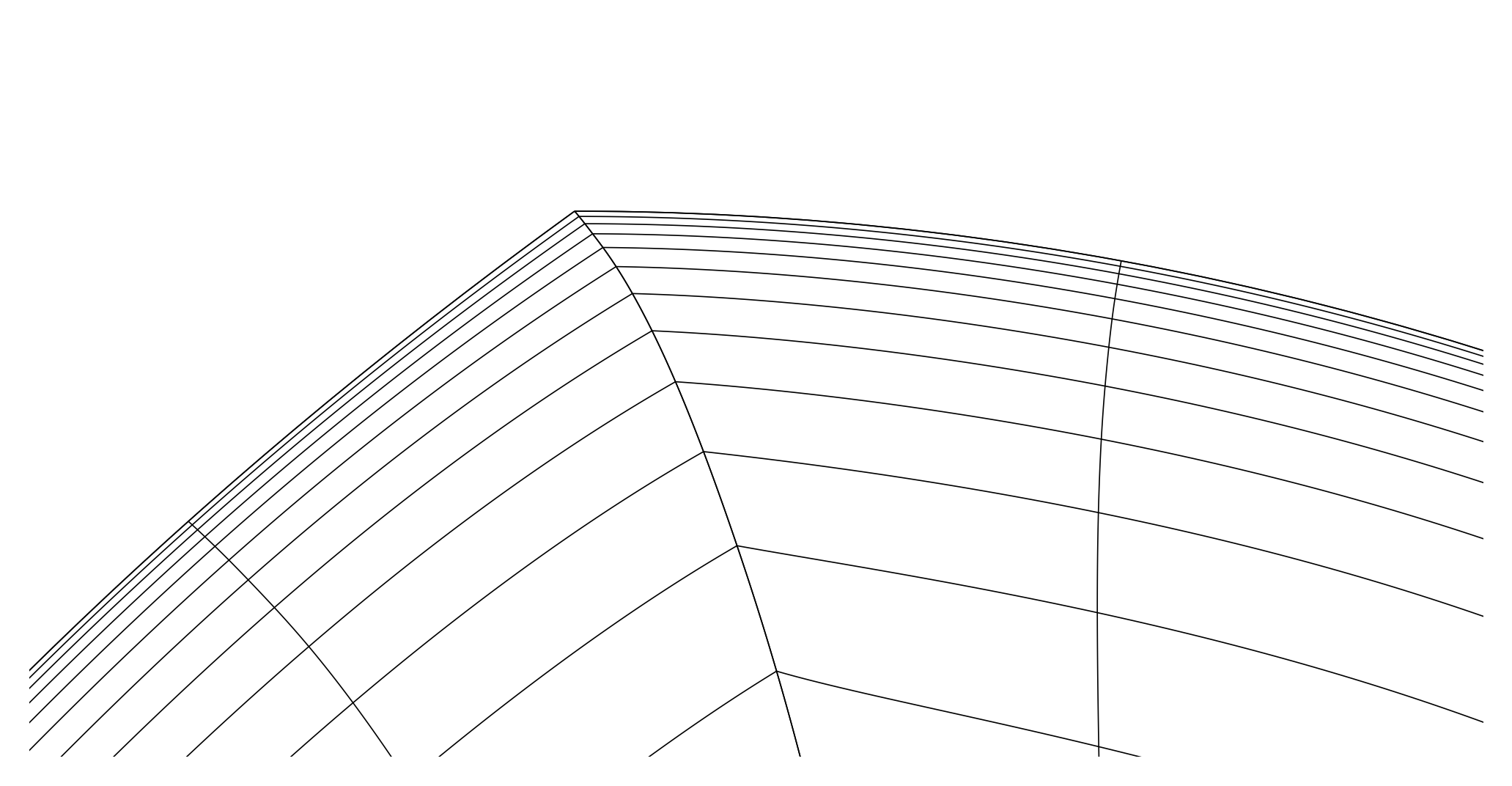} \hspace*{0.1cm}
\includegraphics[height=3cm]{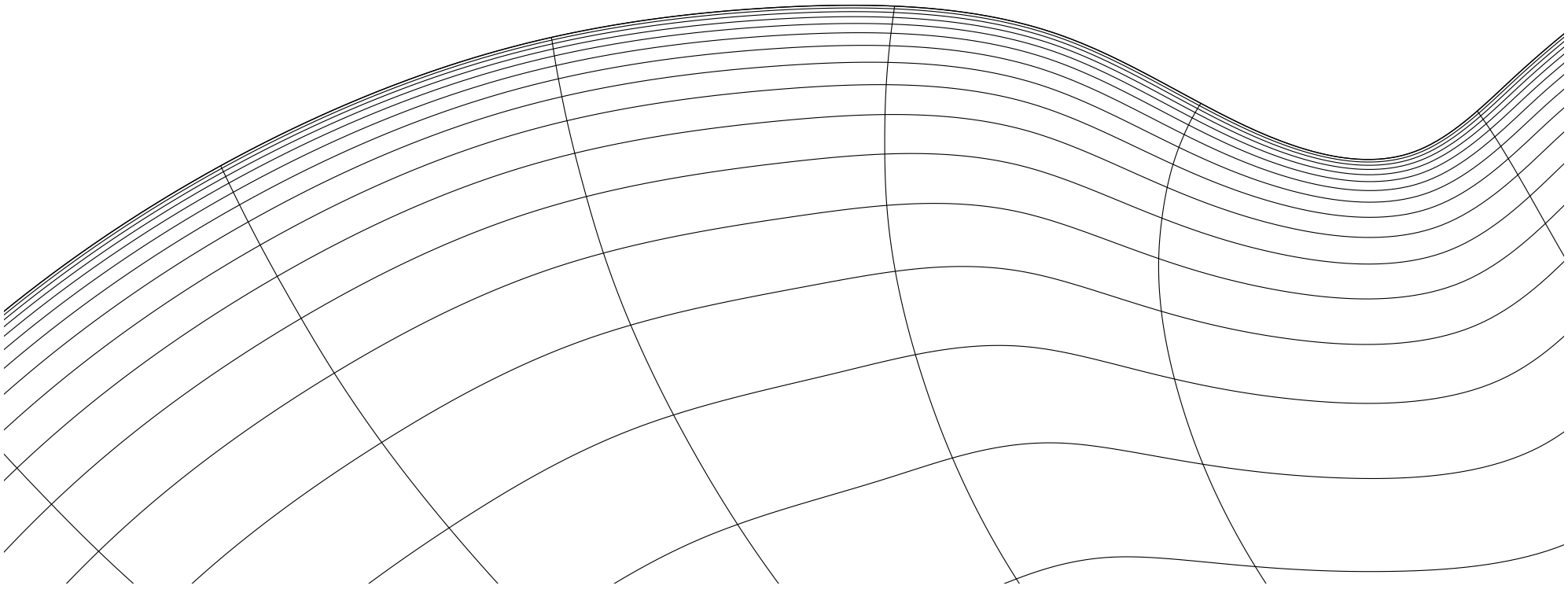}
\caption{A zoom-in on two different segments of the boundary showing the boundary orthogonalisation as well as the boundary layer that is largely maintained along the entire boundary.}
\label{fig:male_screw_boundary_layer_algebraic_zoom}
\end{figure}

\noindent We are considering the screw geometry with orthogonalised reference parameterisation from Figure \ref{fig:male_screw_boundary_orth_geom} (left) and associated controlmap $\bs: \hOm \rightarrow \hOm$ from Figure \ref{fig:male_screw_boundary_orth_param} (left). Using this reference parameterisation, we are creating a boundary layer in which value of $k > 0$ in~\eqref{eq:constant_transverse_deriv_opt} is given by $\overline{b} / 72$, where $\overline{b}$ is the average value of all $b_i$ over all $L_{k^i}$ in~\eqref{eq:boundary_layer_local_tangent_bundle_diagonal}. Figure \ref{fig:male_screw_boundary_layer_algebraic} shows the resulting parameterisation along with the controlmap $\tilde{\bs}: \hOm \rightarrow \hOm$, while Figure \ref{fig:male_screw_boundary_layer_algebraic_zoom} shows a zoom-in on two different segments on the boundary.
As can be seen, the methodology preserves the boundary orthogonalisation while now providing precise control over the density of the boundary layer, which is largely maintained along the entire boundary. As before, the boundary layer density is tuned by the value of $k > 0$. Since the method performs a sequence of algebraic operations on the boundary patches, unlike in the previous methodology, the interior patches remain unchanged. As such, the best results are obtained when the boundary patches cover the majority of $\partial \hOm$ with only a handful of small patches in the interior.

%% file: sec-conclusion.tex
\section{Conclusion}
We have presented a PDE-based parameterisation framework for planar multipatch domains based on the concept of harmonic maps. For this, we presented a total of four different numerical approaches capable of computing valid parameterisations for a wide range of piecewise smooth Lipschitz domains bounded by a collection of spline curves. We presented three different algorithms in nondivergence form, two of which are in mixed form and one based on $C^0$-DG. Furthermore, we presented one approach based on the inverse harmonicity requirement's weak form. We concluded that the NDF-discretisations in mixed form performed similarly well in the essayed benchmark problems while consistently exhibiting slightly better convergence rates than the $C^0$-DG approach. On the other hand, we concluded that the $C^0$-DG approach is the computationally least expensive approach that can be initialised with a degenerate initial iterate. The experiments demonstrated that the weak form discretisation converges reliably when initialised with the solution of one of the NDF-discretisations while performing only marginally worse than or on par with Winslow's original approach. Since the $C^0$-DG approach is usually sufficiently close to the discrete root, we concluded that a combination with the weak form discretisation constitutes a computationally feasible and effective means to compute a uniformly nondegenerate map for the geometries considered in this work. Hereby, the combination of the two approaches substantially reduces the need for a posteriori refinement in case the NDF-solution is degenerate, thanks to the weak form's barrier property. \\
We have augmented the parameterisation framework with mechanisms that allow for control over the map's parametric properties. Hereby, we presented techniques capable of incorporating many of the commonly-desired parametric features into the the computed maps, such as homogeneous cell sizes and boundary layers. For combining harmonic maps with parametric control, we mainly employed the weak form discretisation and concluded that its barrier property is an effective means of maintaining uniform nondegeneracy, even when confronted with coordinate transformations that induce extreme cell size heterogeneity, such as in Figures \ref{fig:male_screw_boundary_layer_algebraic} and \ref{fig:male_screw_boundary_layer_algebraic_zoom}. \\
Utilising an only essentially bounded diffusivity for the purpose of inducing a coordinate transformation via a controlmap, while enabling many novel ways of controlling the outcome, is associated with a number of potential robustness bottlenecks, such as the possibility of creating singularities in the interior of the domain. Here, we proposed a stabilisation via Gaussian blending functions on the quadrangulation's vertices. However, a more thorough investigation of the effect of the controlmap's reduced regularity on the computed maps constitutes a topic for future research. Furthermore, given that the parametric domain is typically given by a convex polygon, we see great potential in the use of computationally inexpensive algebraic methods to create a controlmap that builds desired features into the harmonic map. This constitutes another topic for further research.

\section*{Acknowledgements}
The authors gratefully acknowledge the support of the Swiss National Science Foundation through the project ‘‘Design-through-Analysis (of PDEs): the litmus test’’ n. 40B2-0 187094 (BRIDGE Discovery 2019). Jochen Hinz is grateful for the proof reading and feedback he received from Rafael V\'azquez. \\
Furthermore, the authors are grateful for the coding help they received from the Nutils core development team.